\documentclass{amsart}
\usepackage[
top=30truemm,
bottom=30truemm,
left=30truemm,
right=30truemm]{geometry}
\usepackage{amsfonts, amssymb, amsmath, verbatim, amsthm, mathrsfs, amscd, enumerate, ascmac, fancyhdr, caption, hyperref, framed, subfigure}
\usepackage{bbm}
\usepackage{eucal}
\usepackage{tikz}
\usetikzlibrary{shapes.geometric}
\usetikzlibrary{arrows.meta,arrows}
\usepackage{graphicx} 
\usepackage{pgf, tikz}
\usepackage{pgfplots}
\usepackage{float}
\usepackage{graphicx,color}
\usepackage{tcolorbox}
\usepackage{mathtools}

\def\XXint#1#2#3{{\setbox0=\hbox{$#1{#2#3}{\int}$}
    \vcenter{\hbox{$#2#3$}}\kern-.5\wd0}}

\usepackage[framemethod=tikz]{mdframed}
\mdfsetup{linecolor=blue,backgroundcolor=gray!10,roundcorner=10pt,leftmargin=2pt,innerleftmargin=2pt,innerrightmargin=2pt}
\usepackage{url}

\usepackage{xcolor}
\usepackage{tikz}
\usetikzlibrary{positioning,intersections, calc, arrows,decorations.markings}

\pgfplotsset{compat=1.18}

\newtheorem{theorem}{Theorem}[section]
\newtheorem{corollary}[theorem]{Corollary}

\newtheorem{lemma}[theorem]{Lemma}
\newtheorem{proposition}[theorem]{Proposition}
\theoremstyle{definition}

\newtheorem{remark}[theorem]{Remark}

\newtheorem{example}[theorem]{Example}

\numberwithin{equation}{section}

\usepackage{relsize}

\def\sgn{\mathop{\mathrm{sgn}}\nolimits}
\def\supp{{\rm supp\,}}

\def\diam{{\rm diam\,}}

\def\rmd{{\rm d}}
\def\1{\mathbbm{1}}
\def\Rmap{R_{\xi}\eta}

\usepackage[dvipsnames]{xcolor}

\definecolor{scisRed}{RGB}{190, 90, 90}

\definecolor{purple}{HTML}{645394}

\definecolor{aubergine}{HTML}{795367}

\definecolor{rose}{HTML}{AD718B}

\definecolor{lightgreen}{HTML}{A9C6B3}

\definecolor{green}{HTML}{587c64}

\definecolor{teal}{HTML}{075f6f}

\definecolor{blue}{HTML}{4670a8}

\makeatletter
\renewcommand{\l@section}{\@tocline{1}{0pt}{1pc}{}{\scshape}}
\renewcommand{\l@subsection}{\@tocline{2}{0pt}{2.3pc}{2.8pc}{}}
\makeatother
\mathtoolsset{showonlyrefs} 

\title[Wigner Representations of the Extension Operator for Curves]{Wigner Representations and
Weighted $L^2$ Inequalities for the Fourier Extension Operator for Curves}

\author[A. Tierney]{Amy Tierney}
\address[Amy Tierney]{School of Mathematics, The Watson Building, University of Birmingham, Edgbaston, Birmingham, B15 2TT, England.}
\email{axt944@bham.ac.uk; amyjtierney@gmail.com}

\begin{document}

\date{\today}


\begin{abstract}
For a non-degenerate smooth curve $S$ in $\mathbb{R}^n$, we prove a Wigner representation for its associated Fourier extension operator: $|Eg|^2 = X_{2k}^*W_{2k}g$. Here $X_{2k}$ is a restriction of the $X$-Ray transform via the $2k$-th Frenet vector, and $W_{2k}$ is a new curve-carried Wigner transform.
When $S$ is the helix, we additionally construct a representation in terms of a restriction of the $\lfloor n/2 \rfloor$-plane transform.
Consequently, we derive weighted $L^2$ inequalities for the extension operator for the curve. Explicitly, we prove Sobolev--Mizohata--Takeuchi inequalities and Mizohata--Takeuchi-type inequalities for orthonormal families of functions.
\end{abstract}

\maketitle
\vspace{-2em}

{\setlength{\parskip}{0.2pt}\tableofcontents}
\vspace{-2em}
\newpage
\section{Introduction}
A prominent object of interest in harmonic analysis is the \textit{Fourier extension operator} associated with a submanifold $S$ of $\mathbb{R}^n$, given by
\begin{equation}
    \widehat{g\rmd\sigma_S}(x) = \int_S g(\xi)e^{-2\pi i x \cdot \xi} \rmd\sigma_S(\xi), \qquad g \in L^1(\rmd\sigma_S),
\end{equation}
where $\rmd\sigma_S$ denotes the surface measure on $S$. 
Considerable attention has been devoted to its \textit{weighted $L^2$ inequalities}, which form the focus of this paper.

One strategy to make progress on such problems is to prove useful representations of 
$|\widehat{g\rmd\sigma_S}|^2$.
Some such representations have been in terms of a tomographic transform (see \cite{Planchon_Vega, BNS, BN_tomography} and the references therein). 
The tomographic object of interest, both in these representations and the current paper, is the \textit{$m$-plane transform} defined via
\begin{equation}
\label{k_plane_defn}
    T_{m, n}f(\pi, v) = \int_{\pi}f(y + v) \rmd y,
\end{equation}
where $\pi$ is an $m$-plane and $v \in \pi^{\perp}$ is a translation parameter, together forming the standard parametrisation of the Grassman manifold of $m$-planes, $\mathfrak{M}_{m, n}$. Moreover, $\rmd y$ denotes Lebesgue measure on $\pi$. In the cases of $m=1$ and $m=n-1$, these are the $X$-Ray and Radon transforms respectively, denoted by $X$ and $\mathcal{R}$.

Another wider benefit of such an expression in the case of the sphere comes from observing that $|\widehat{g\rmd \sigma_S}|^2$ represents a certain local intensity of a monochromatic optical wave field and therefore its representations hold importance in optics -- see, for example \cite{optics_alonso}.
\subsection{Wigner representations for hypersurfaces}
The present work is inspired by a recent paper of Bennett, Gutiérrez, Nakamura and Oliveira
\cite{BGNO}, wherein they establish a representation of $|\widehat{g\rmd\sigma_S}|^2$ for $S$ a smooth, strictly convex hypersurface. Their representation again features a tomographic transform (here the $X$-Ray transform) but also \textit{a Wigner transform}; the key strength and novelty of their identities comes from an appropriately constructed Wigner transform dependent on the underlying hypersurface.

The purpose of the present paper is to prove an analogous representation for curves in $\mathbb{R}^n$ by constructing appropriate Wigner transforms. To highlight the role of the geometrically defined Wigner transform, we refer to such expressions as \textit{Wigner representations}. Consequently, we derive weighted $L^2$ inequalities motivated by \cite{BGNO} and orthonormal weighted extension estimates motivated by \cite{BBGNO}.

\medspace
\noindent
\textbf{The paraboloid and the classical Wigner transform.}
Introduced by Wigner in 1932 \cite{Wigner}, the \textit{classical Wigner transform} is defined via
\begin{equation}
\label{classical_Wigner_intro}
    W(g_1, g_2)(x,v) = \int_{\mathbb{R}^d} g_1\left(x-\frac{y}{2}\right)\overline{g_2\left(x + \frac{y}{2}\right)}e^{-2\pi i v \cdot y} \rmd y, \qquad 
g_1, g_2 \in L^2(\mathbb{R}^d).
\end{equation}
Wigner observed that for a solution to the free Schrödinger equation $2 \pi i \frac{\partial u}{\partial t} =  \Delta u$ with initial data,
$u(x, 0) = u_0(x)$, the Wigner transform of the solution, $F(x,v,t) := W(u(\cdot, t), u(\cdot, t))(x, v)$, satisfies the transport equation
\begin{equation}
    \frac{\partial F}{\partial t}= 2v \cdot \nabla_x F,
\end{equation}
whose solution is $F(x, v, t) = F_0(x + 2tv, v)$, where $F_0:= W(u_0,u_0)$.
Another important property of the Wigner transform, referred to as the marginal property, states that
\begin{equation}
\label{first_marginal_classical_WT}
    |u(x,t)|^2 = \int_{\mathbb{R}^{d}} F(x,v,t) dv.
\end{equation}
Combining the above properties, observe that
\begin{equation}
\label{BGNO_parab_phase_space}
    |u(x,t)|^2 = \int_{\mathbb{R}^{d}} F(x,v,t) dv = \int_{\mathbb{R}^d} F_0(x + 2tv, v) dv :=\rho W(u_0, u_0)(x,t),
\end{equation}
where $\rho$ can be interpreted as an adjoint parametrised $X$-Ray transform with further details to follow.
Since the parametrised extension operator for the paraboloid,
\begin{equation}
   u(x, t)= \int_{\mathbb{R}^{d}} \widehat{u_0}(\xi)e^{-2\pi i( x \cdot \xi + t|\xi|^2)} d\xi, \qquad u_0 \in L^1(\mathbb{R}^d),
\end{equation}
solves the above Schrödinger equation with initial data $u_0$, \eqref{BGNO_parab_phase_space} can be interpreted as a phase-space representation for the extension operator for the paraboloid expressed in terms of the Wigner transform.

\medspace
\noindent
\textbf{{Hypersurfaces.}} Motivated by the above, 
Bennett, Gutiérrez, Nakamura and Oliveira \cite{BGNO} established a more general Wigner representation
\begin{equation}
|\widehat{g\rmd\sigma_S}|^2 = X^{*}_SW_Sg,
\end{equation}
for \textit{smooth, strictly convex hypersurfaces $S$}. See also Alonso \cite{optics_alonso} for an earlier discussion of Wigner representations arising in optics.
Here, they introduced a \textit{surface-carried Wigner transform} defined via
\begin{equation}
\label{intro:WT_hypersurface}
    W_S(g_1, g_2)(u,y) = \int_{S}g_1(u') \overline{g_2(R_uu')} e^{-2\pi iy\cdot (u' - R_u u')}J(u,u')\rmd\sigma(u'), \qquad g_1, g_2 \in L^2(\mathbb{R}^n),
\end{equation}
where $(u,y)\in TS$, the tangent bundle. For $u \neq u'$, $R_uu'$ is defined to be the unique point $u'' \in S$ with $u'' \neq u'$ such that
\begin{equation}
    N(u) \cdot (u' - u'')=0
\qquad \text{and} \qquad
    N(u) \wedge N(u') \wedge N(u'') =0. \qquad \footnote{It is important to note here that the second condition is imposed to guarantee uniqueness and in two dimensions this condition becomes redundant, as such the former condition is the one to keep in mind when comparing to curves later.}
\end{equation}
Moreover, $R_uu := u$ for all $u \in S$; $J(u,u')$ is the reciprocal of the Jacobian of the mapping $u \mapsto R_u u'$ expressible in terms of the shape operator of $S$; and $N: S \to \mathbb{S}^{n-1}$ is the Gauss map, with $\mathbb{S}^{n-1}$ denoting the $(n-1)$-dimensional sphere.
The \textit{$X$-ray transform} can also be expressed via
\begin{equation}
Xf(\omega, v) = \int_{\mathbb{R}}f(v + s\omega) \rmd s,
\end{equation}
for a direction $\omega \in \mathbb{S}^{n-1}$  and a translation parameter $v \in \langle \omega \rangle^{\perp}$. 
Furthermore, $X_S$ is a restricted $X$-Ray transform in the following sense:
\begin{equation}
\label{restricted_X-ray_BGNO}
    X_Sf(u, v) = Xf(N(u), v), \qquad u \in S, v \in TS,
\end{equation}
and $X_S^{*}$, appearing in the phase-space formulation, refers to its adjoint.
\subsection{Main results: Wigner representations for curves}
In the present work, we seek a Wigner representation of $|\widehat{gd\sigma_S}|^2$ when $S$ is a curve in $\mathbb{R}^n$.
In particular, we consider compact $C^{\alpha}$ curves in $\mathbb{R}^n$, with $\alpha \geq n+1$, parametrised by arclength so that
    $S = \{\phi(t): t \in I \}$, for some interval $I \subset \mathbb{R}$, and assume that $S$ is \textit{well-curved}, that is,
\begin{equation}
\label{well_curved_n}
    |\phi'(\xi) \wedge \phi''(\xi) \wedge \dots \wedge \phi^{(n)}(\xi)| \neq 0 \qquad \text{for all } \xi \in I.
\end{equation} 
We work with the parametrised extension operator for this curve, given by
\begin{equation}
    E_Sf(x) = \int_I f(t) e^{-2\pi i x \cdot \phi(t)} \rmd t, \qquad f\in L^1(I).
\end{equation}
Essential to the study of curves in this work are their associated Frenet frame and curvatures, denoted by $(e_1, \dots, e_n)$ and $\kappa_1, \dots, \kappa_{n-1}$. The first three vectors denote the tangent, normal and binormal to the curve, and we emphasise this in three dimensions by denoting $T = e_1$, $N= e_2$ and $B=e_3$. 
The final Frenet curvature $\kappa_{n-1}$ is sometimes referred to as the torsion. In three dimensions, this is highlighted by referring to $\kappa_1 = \kappa$ as the curvature and $\kappa_2 = \tau$ as the torsion.

\medskip
\noindent
\textbf{Wigner Transforms for Curves in $\mathbb{R}^3$.}
For the sake of exposition, we begin by presenting the Wigner transform and representation for $C^4$ curves in $\mathbb{R}^3$. We impose additional technical conditions that amount to restricting the length of the interval $I$ so that the Frenet frame cannot rotate too much in the interval, but defer discussion to Section \ref{Section_curvesR3}.

For $g_1, g_2 \in L^2(I)$, we define the \textit{curve-carried Wigner transform}, via
\begin{equation}
\label{intro:WT_curveR3}
    W_N(g_1, g_2)(\xi,v) = \int_{I}g_1(\eta) \overline{g_2(R_{\xi}\eta)} e^{-2\pi iv\cdot (\phi(\eta) - \phi(R_{\xi}\eta))}J(\xi,\eta)\rmd \eta,
\end{equation}
where $(\xi, v) \in \{ (\xi, v): \xi \in I, v \in \langle N(\xi) \rangle^{\perp}
    \}:= D_N.$ Despite the apparent similarity between the hypersurface and curve-carried Wigner transforms, the domain in the curve case is the (principal) normal line bundle rather than the tangent bundle.
For $\xi \neq \eta$, $R_{\xi}\eta$ is defined to be the unique point $\zeta \in I$ such that $\zeta \neq \eta$ and
\begin{equation}
     N(\xi) \cdot (\phi(\eta) - \phi(\zeta)) =0.
\end{equation}
Moreover, $R_{\xi}\xi := \xi$ for all $\xi \in I$ and in analogy with the surface-carried Wigner transform,
$J(\xi,\eta)$ is the reciprocal of the Jacobian of the mapping $\xi \mapsto R_{\xi} \eta$. Explicitly, for $\xi \neq \eta$,
\begin{equation}
        J(\xi, \eta) = \frac{|N'(\xi)\cdot (\phi(\eta) - \phi(\Rmap))|}{|N(\xi)\cdot T(\Rmap)|},
\end{equation}
where we note that the Frenet formulae state $N'(\xi) = - \kappa(\xi) T(\xi) + \tau(\xi) B(\xi)$, and the expression encompasses contributions from all of the Frenet vectors.
When the input functions are the same, we refer to this as the \textit{curve-carried Wigner distribution}, and sometimes denote $W_N(g, g) = W_Ng$.
Our Wigner representation and main result for curves in $\mathbb{R}^3$ is as follows:
\begin{theorem}[\textbf{Wigner representation for Curves in $\mathbb{R}^3$}]
\label{thm_intro:wigner_rep_3d}
Let $S$ be a compact $C^4$ curve in $\mathbb{R}^3$ as described above. Then,
\begin{equation}
    |E_Sg|^2 = X^{*}_NW_Ng,
\end{equation}
where $X_Nw(\xi, v) := Xw(N(\xi), v)$ for $(\xi, v) \in D_N$.
\end{theorem}
For functions $h:D_N \to \mathbb{C}$, we have the following explicit formula for the adjoint restricted $X$-Ray transform:
\begin{equation}
    X_{N}^* h(x) = \int_I h(\xi, P_{\langle N(\xi)\rangle^{\perp}}(x)) \rmd \xi,
\end{equation}
which may be interpreted as a superposition of dirac masses on lines $\ell_{\xi, v}:= \{ v + t N(\xi) : t \in \mathbb{R}\}$.

\medskip
\noindent
\textbf{Comparison to the wave-packet approach for curves.}
    A powerful technique in harmonic analysis involves expressing the extension operator for a submanifold of a function as a superposition of wavepackets. One such application of the wavepacket decomposition is in proving weighted $L^2$ inequalities referred to as Mizohata--Takeuchi-type inequalities: for hypersurfaces refer to \cite{CIW} and for curves to \cite{CLPY_MT_curves} (see also the references therein).

 When the underlying submanifold is a curve,
 the wavepackets are adapted to slabs. Explicitly, for a fixed point $\xi_i$ in a discrete set of $R^{-\frac{1}{n}}$ separated points on the curve, the slabs are parallelepipeds whose edges have directions given by the Frenet frame at $\xi_i$, with side length approximately $R^{\frac{j}{n}}$ in the $j$-th Frenet direction (see \cite{CLPY_MT_curves} for further detail). Taking into consideration the differences to be expected from the comparison in the setting of hypersurfaces (see \cite{CIW} and Remark 1.15 of \cite{BGNO}), the wavepacket decomposition for curves might suggest that our representation ought to be superpositions of hyperplanes with normals in the Frenet tangent direction. Alternatively, after rescaling, a superposition of lines with direction parallel to the final Frenet direction ($e_n(\xi)$ for some point $\xi \in I$) could also be suggested.

Therefore, one might wonder whether it is possible to construct a Wigner transform $W_R$, with the same structure as \eqref{intro:WT_hypersurface} and \eqref{intro:WT_curveR3} but different defining condition for $R_{\xi}\eta$, for which $|E_Sg|^2 = \mathcal{R}_T^* W_R g$. Here, the restriction of the Radon transform, 
\begin{equation}
   \mathcal{R}_Tf(\xi, v_t) := \mathcal{R}f(\langle T(\xi) \rangle^{\perp}, v_t), \qquad \xi \in I, v_t \in \langle T(\xi) \rangle,
\end{equation}
is chosen in accordance with the wavepacket considerations.

In terms of the defining condition on $\Rmap$, this question would be to ask for fixed $\xi, \eta \in I$, does there exist a third point $\zeta\in I$ with $\zeta \neq \eta$ for which both
\begin{align}
    N(\xi)\cdot(\phi(\eta) - \phi(\zeta)) = 0 \qquad \text{and} \qquad B(\xi) \cdot (\phi(\eta) - \phi(\zeta)) =0?
\end{align}
However, for a suitably small interval $I$, there \textit{does not exist} a third point $\zeta \in I$ such that $\zeta \neq \eta$ satisfying both equations. In fact, there \textit{does not exist} a third point $\zeta \in I$ with $\zeta \neq \eta$ satisfying the second condition. It is instructive to consider the example of the moment curve parametrised by $(t, t^2, t^3)$ for say $t \in [-1, 1]$ and $\xi =0$, in which case this statement is immediate.
The statement implies that not only can we not construct a representation via this Wigner transform as a superposition of hyperplanes, but we also cannot do so in terms of the final Frenet direction (the binormal direction), at least with the current approach to constructing Wigner transforms.

\medskip
\noindent
\textbf{Wigner transforms for curves in $\mathbb{R}^n$.}
In higher dimensions, interestingly, for a suitably small interval $I$, there do not exist three distinct points $\xi, \eta, \zeta \in I$ satisfying
\begin{equation}
    e_{2k-1}(\xi) \cdot (\phi(\eta) - \phi(\zeta))=0.
\end{equation} 
We defer to Proposition \ref{prop:no_third_pt} for its formal statement; its proof is illuminating in understanding why our methods invoke results with a strong dependence on the parity of the Frenet vectors.

That we cannot impose a defining condition for $\Rmap$ in terms of the odd Frenet vectors motivates the following definition of the Wigner transform in higher dimensions.
Fix $k \in \mathbb{N}$, $2k \leq n$ and let $S$ be a $C^{n+ 2k -1}$ curve in $\mathbb{R}^n$ satisfying certain technical conditions to be clarified in Section \ref{Section_curvesRn}.
For $g_1, g_2 \in L^2(I)$, we define the \textit{$k^{\text{th}}$ curve-carried Wigner transform},  via
\begin{equation}
    W_{2k}(g_1, g_2)(\xi,v) = \int_{I}g_1(\eta) \overline{g_2(R_{\xi}\eta)} e^{-2\pi iv\cdot (\phi(\eta) - \phi(R_{\xi}\eta))}J(\xi,\eta)\rmd \eta.
\end{equation}
Here, $(\xi, v) \in \{ (\xi, v): \xi \in I, v \in \langle e_{2k}(\xi) \rangle^{\perp}
    \}:= D_{2k}$, the $k^{\text{th}}$-Frenet line bundle.
For $\xi \neq \eta$, $R_{\xi}\eta$ is defined to be the unique point $\zeta \in I$ such that
\begin{equation}
    e_{2k}(\xi)\cdot(\phi(\eta) - \phi(\zeta)) =0.
\end{equation}
Moreover, $R_{\xi}\xi := \xi$ for all $\xi \in I$ and
  \begin{equation}
        J(\xi, \eta) := \frac{|e_{2k}'(\xi)\cdot (\phi(\eta) - \phi(\Rmap))|}{|e_{2k}(\xi)\cdot e_1(\Rmap)|},
    \end{equation}
    with $e_{2k}'(\xi)= -\kappa_{j-1}(\xi) e_{j-1}(\xi) + \kappa_j e_{j+1}(\xi)$.
\begin{theorem}[Wigner Representation for the extension operator for Curves in $\mathbb{R}^n$]
\label{intro_main_thm}
    Let $S$ be a well-curved, compact patch of a $C^{n+2k -1}$ curve in $\mathbb{R}^n$ as described above. Then,
\begin{equation}
    |E_Sg|^2 = X^{*}_{2k}W_{2k}g,
\end{equation} 
where for $(\xi, v) \in D_k$, we define a restriction of the $X$-Ray transform via
\begin{equation}
    X_{2k}f(\xi,v) := Xf(e_{2k}(\xi), v).
\end{equation}
\end{theorem} 
For $h:D_{2k} \to \mathbb{C}$, the adjoint restricted $X$-Ray transform $X_{2k}^*h$ is a superposition of dirac masses on lines $\ell_{\xi, v}:= \{v + t e_{2k}(\xi): t \in \mathbb{R}\}$.

\medskip
\textbf{An `optimal' Wigner representation for the helix.}
Proposition \ref{prop:no_third_pt} implies that via a Wigner transform with structure similar to those presented, but different defining condition for $\Rmap$, the best possible representation would be
\begin{equation}
\label{best_possible_WR}
   |E_Sg|^2 = T_{e}^*W_{\star}g,
\end{equation}
where $W_{\star}g$ is the `best possible' Wigner transform with $R_{\xi}\eta$ defined to satisfy
\begin{equation}
    e_{2k}(\xi) \cdot (\phi(\eta) - \phi(\zeta))=0
\end{equation}
simultaneously \textit{for all} $k \in \mathbb{N}$ with $2k \leq n$.
Moreover, $$T_ef(\xi, v) := T_{\lfloor \frac{n}{2}\rfloor, n}f(Y_e, v)$$ for $(\xi, v) \in D_o:= \{(\xi, v): \xi \in I, v \in Y_o\}$, where $Y_o$ and $Y_e$ are the spans of the odd and even Frenet vectors, respectively. 
The best possible representation \eqref{best_possible_WR} differs from the Wigner representation we prove for $n>3$.

However, \eqref{best_possible_WR} \textit{can be attained} in the special case of the \textit{helix} in $\mathbb{R}^n$ which is categorised by its constant Frenet curvatures $\kappa_1, ..., \kappa_{n-1}$.
The optimal representation is possible in this case since $\Rmap = 2\xi - \eta$ satisfies 
\begin{equation}
    e_{2k}(\xi)\cdot(\phi(\eta) - \phi(\Rmap)) =0
\end{equation}
for all $\xi, \eta \in I$ and \textbf{all} $k\in \mathbb{N}$ such that $2k \leq n$ when $\phi: I \to \mathbb{R}^n$ is the arclength parametrised helix.

\subsection{Applications of the Wigner representations}
Equipped with geometric identities for the extension operator for curves, we make two observations which motivate two applications of the Wigner representation. Let $w: \mathbb{R}^n \to \mathbb{C}$ be a weight function\footnote{This is to say, we consider a suitable class of weights such as those which are bounded and have compact support.}
 and let $g \in L^2(I)$. Then,
\begin{equation}
\label{Intro_obs1}
    \int_{\mathbb{R}^n}|E_S g(x)|^2 w(x) \rmd x = \int_I \int_{\langle e_{2k}(\xi) \rangle^{\perp}} W_{2k}g(\xi,v)X_{2k}w(\xi, v) \rmd v \rmd\xi.
\end{equation}
Moreover, for a sequence $(g_j)$ of $L^2$ functions on $I$,
\begin{equation}
\label{Intro_obs2}
    \int_{\mathbb{R}^n}\sum_j\lambda_j|E_S g_j(x)|^2 w(x) \rmd x = \int_I \int_{\langle e_{2k}(\xi) \rangle^{\perp}} \left(\sum_j \lambda_j W_{2k}g_j(\xi,v)\right)X_{2k}w(\xi, v) \rmd v \rmd\xi.
\end{equation}
The strength of the second observation \eqref{Intro_obs2} is in its ability to `linearise' $\sum_j\lambda_j|E_Sg_j|^2$ leading to applications in orthonormal weighted extension estimates for curves. Our application (Theorem \ref{thm_orthonormal_MT_Rn_curve}) is motivated by \cite{BBGNO} where the authors use the Wigner representation for hypersurfaces from \cite{BGNO} to prove orthonormal Mizohata--Takeuchi-type estimates for hypersurfaces.

Observation \eqref{Intro_obs1} suggests that the Wigner representation may be useful in proving weighted $L^2$ estimates related to the Mizohata--Takeuchi conjecture. For a smooth $m$-dimensional surface $S$ in $\mathbb{R}^n$ and a \textit{non-negative} locally integrable weight $w$, this asks for what $\alpha \geq 0$ does there exist a finite implicit constant such that
\begin{equation}
\label{MT_conj_dim_m}
\tag{MT}
    \int_{B_R} |\widehat{g\rmd\sigma_S}(x)|^2w(x) \rmd x \lesssim_{\alpha} R^{\alpha} \|g\|_{L^2(\rmd\sigma_S)}^{2}  \sup_{(u,v) \in TS}{T_{n-m, n}w((T_{u}S)^{\perp}, v)}?\footnote{In this work, we write $A \lesssim B$ to denote $A \leq cB$, and highlight the dependence on any parameters by writing say $A \lesssim_{\alpha} B$ to indicate that the constant $c$ depends on $\alpha$.}
\end{equation}
Much of the progress surrounding the Mizohata--Takeuchi conjecture has been concerned with the hypersurface setting ($m=n-1)$, with $T_{n-m, n}w((T_{u}S)^{\perp}, v) = X_S w(u, v)$, the restricted $X$-Ray transform of \eqref{restricted_X-ray_BGNO}. Refer to \cite{Takeuchi1974, Takeuchi1980, mizohata, Radial_weights_Barcelo_Ruiz_Vega, Radial_weights_Carbery_Soria} for its origins and \cite{ BBC_BARCELÓ_BENNETT_CARBERY_2008, CIW, Ortiz_MT_cone, MT_random_weights, Yang_MT_exponential_sums, triangles_BKNO} for further results pertaining to hypersurfaces.
We highlight the important recent development of Cairo's counterexample \cite{counterexample_MT} demonstrating that \eqref{MT_conj_dim_m} is false when $\alpha = 0$ for all hypersurfaces that are not part of a hyperplane. See the very recent paper \cite{MT_local_counterexample} demonstrating that the local conjecture is also false in this case.
For earlier developments in this direction, refer to \cite{Cairo_Zhang, Fenves_MT_counterexample, zhang_Stein}. 

Regarding curves in $n$-dimensions, we refer to \cite{CLPY_MT_curves} for the best known sufficient $\alpha$ for which \eqref{MT_conj_dim_m} holds and \cite{Cairo_Zhang} for a necessary lower bound on $\alpha$ for certain $C^k$ curves. Both papers also consider a variant of \eqref{MT_conj_dim_m} for curves with an $X$-Ray transform instead of a Radon transform, in which case $\alpha > n-2$ is required due to scaling considerations.

Motivated by the first observation of our Wigner representation \eqref{Intro_obs1}, we seek global bounds with a different norm on the restricted tomographic transform of the weight in the translation parameter. Such variants of \eqref{MT_conj_dim_m} are motivated by
\cite{DMV} for $L^2$-variants when $S$ is a quadratic submanifold; \cite{BGNO} for $L^2$-Sobolev variants when $S$ is a hypersurface and \cite{Ferrante_thesis} for $L^1$ Sobolev variants when $S$ is an $m$-dimensional submanifold.
We prove $L^2$-Sobolev variants for curves:
\begin{theorem}[$L^2$ Sobolev--Mizohata--Takeuchi Inequality for Curves]
\label{intro:Sob_MT}Let $S$ be a well-curved, compact patch of a $C^{n+2k -1}$ curve in $\mathbb{R}^n$ as described above. Moreover, let $s < \frac{1}{2}$. Then there exists a finite implicit constant depending only on $s$ and the diameter of $S$ such that
    \begin{equation}
    \int_{\mathbb{R}^n} |E_Sg(x)|^2w(x) \rmd x 
    \lesssim_s \: C_k(\phi) \: \left\|g\right\|_{L^2(I)}^2 \sup_{\xi \in \supp^*(g)} \left( \int_{\langle e_1(\xi), e_{2k}(\xi) \rangle^\perp} \left\| X_{2k}w(\xi, \cdot + v_b)) \right\|_{\dot{H}_{v_t}^s(\langle e_1(\xi) \rangle)} \rmd v_b \right),
\end{equation}
for all weight functions $w$.
Here, $\supp^*(g) := \{ \xi \in I: \Rmap \in \supp(g)  \text{ for some } \eta \in \supp(g)\}$ and $\dot{H}^s$ is the homogeneous Sobolev space.
\end{theorem}
We defer to Theorem \ref{corol2b_l2sobolev_n} for the full statement, including the simple explicit expression of $C_k(\phi)$.
\begin{remark}[Signed weights]
    Our results do not require non-negativity of the weight, a condition imposed in the Mizohata--Takeuchi conjecture. This was expected as the same permissibility of signed weights occurs in \cite{BGNO}.
\end{remark}

As suggested by the Wigner representations, we obtain stronger results for the helix, where the $X$-Ray transform is replaced by an $\lfloor \frac{n}{2}\rfloor$-plane transform:
\begin{theorem}[$L^2$ Sobolev--Mizohata--Takeuchi Inequality for the Helix in $\mathbb{R}^n$]
\label{intro:Sob_MT_helix}
Let $S$ be the helix with arclength parametrisation $\phi: I \to \mathbb{R}^n$. Moreover, let $s < \frac{1}{2}$. Then there exists a constant $c$ depending on $s$ such that
    \begin{equation}
    \int_{\mathbb{R}^n} |E_Sg(x)|^2w(x) \rmd x 
    \leq c \: \left\|g\right\|_{L^2(I)}^2 \sup_{\xi \in \supp^*(g)} \left( \int_{Y_b(\xi)} \left\| T_e w(\xi, \cdot + v_b)) \right\|_{\dot{H}_{v_t}^s(\langle e_1(\xi) \rangle)} \rmd v_b \right),
\end{equation}
where $Y_b(\xi) := \langle \{ e_{2k+1}(\xi): k \in \mathbb{N}, 2k\leq n\} \rangle$.
\end{theorem}

\textbf{Stein's conjecture.}
The Sobolev--Mizohata--Takeuchi inequalities follow from a Sobolev variant of a stronger conjectural inequality attributed to Stein (see Theorem \ref{corol2a_l2sobolev_n}). Presented by Stein in the 1978 Williamstown harmonic analysis conference \cite{Stein1978} in the context of the disk multiplier, the conjecture also has an extension operator form. It asks for what $\alpha$, do we have that 
\begin{equation}
\label{Stein_conj_dim_m}
\tag{ST}
    \int_{B_R} |\widehat{g\rmd\sigma_S}(x)|^2w(x) \rmd x \lesssim_{\alpha} R^{\alpha} \int_S |g(u)|^{2}  \sup_{v \in T_uS}{T_{n-m, n}w(u, v)} \rmd\sigma_S(u)?
\end{equation}
Since the verity of this inequality would imply the Mizohata--Takeuchi inequality, the counterexample \cite{MT_local_counterexample} again tells us that the local statement is false for all hypersurfaces that are not contained in a hyperplane.

\textbf{The controlling constant.} A distinctive feature of the Stein and Mizohata--Takeuchi conjectures compared with the classical Fourier restriction conjecture is that \eqref{MT_conj_dim_m} and \eqref{Stein_conj_dim_m} \textit{are true} for $\alpha =0$ when $S$ is (a patch of) an $m$-plane.
Therefore, it is natural to expect global Stein and Mizohata--Takeuchi-type inequalities for curves to \textit{not depend on a lower bound on any of the Frenet curvatures}. Indeed, the Sobolev--Mizohata--Takeuchi inequality for hypersurfaces of \cite{BGNO} did not depend on the curvature of $S$ in absolute terms. 

Both the Sobolev--Stein and Mizohata--Takeuchi inequalities are controlled by the constant $C_{k}(\phi)$ with explicit expression provided in Theorem \ref{corol2b_l2sobolev_n}.
In two dimensions, $C_1(\phi)$ is controlled by the same scale-invariant ratio of curvatures as in \cite{BGNO}. In three dimensions, we will later see (in Lemma \ref{lem:control_Cphi_3d}), that $C(\phi)$ may be controlled by the simpler quantity
\begin{equation}
    \tilde{C}(\phi) := \sup_{\xi, \eta \in I} \frac{\kappa(\xi) + |\tau(\xi)|}{\kappa(\eta)}.
\end{equation}
This quantity is scale invariant; we defer the details to Remark \ref{rem:scale_invar_const_3d}.

\medskip
\textbf{Further contextual remarks.} 
 The restricted Radon transform appearing on the right-hand side of \eqref{MT_conj_dim_m} for curves can be expressed as
\begin{equation}
\label{intro:Radon_Xray}
\mathcal{R}_Tw(\xi, v_t) = \int_{\langle e_1(\xi), e_{2k}(\xi)\rangle^{\perp}}  X_{2k}w(\xi, v_t + v_b)     \rmd v_b 
\end{equation}
for $\xi \in I$ and $v_t \in \langle e_1(\xi) \rangle$.
Observe that,
\begin{equation} 
\sup_{v_t} \mathcal{R}_Tw(\xi, v_t) \leq  \int_{\langle e_1(\xi), e_{2k}(\xi)\rangle^{\perp}} \sup_{v_t} X_{2k}w(\xi, v_t + v_b)     \rmd v_b.
\end{equation}
To compare the right-hand side of this inequality with the quantity appearing in Theorem \ref{intro:Sob_MT}, one may consider analogous Sobolev-embedding considerations to those of Remark 2.9 of \cite{BGNO}. 

The Sobolev variants are strongest when the weight is smooth.
For example, taking $w(x) = \psi(\frac{x}{R})$ for a smooth bump function $\psi$ in Theorem \ref{intro:Sob_MT}, recovers the Agmon--Hörmander inequality for curves \cite{Agmon_Hormander} up to an $\varepsilon$-loss, explicitly
\begin{equation}
    \int_{B(0, R)}|E_Sg|^2 \rmd x \lesssim R^{n-1 + \varepsilon} \|g\|_2^2.
\end{equation}
See also \cite{CLPY_MT_curves} for further context of this inequality in the setting of Mizohata--Takeuchi-type inequalities for curves
and Remark 1.13 of \cite{BGNO} for a relevant comparison of the strength of the Sobolev inequality for different weights.
For further clarification on the mixed norm of the tomographic transform on the right-hand side of Theorem \ref{intro:Sob_MT}, refer to Section \ref{subsec:wigner_heuristicsR3}, and for additional contextual remarks, refer to Section \ref{sec:orthonormal}.

\medskip
\textbf{Structure of the paper.}
Section \ref{Section_curvesR3} considers curves in $\mathbb{R}^3$, providing the framework in a more geometrically intuitive setting before considering curves in $\mathbb{R}^n$ in Section \ref{Section_curvesRn}. Both cover the Wigner representations, that is Theorems \ref{thm_intro:wigner_rep_3d} and \ref{intro_main_thm} respectively, and the first application of the representations, that is Theorem \ref{intro:Sob_MT}.
Section \ref{sec:orthonormal} then considers the second application to orthonormal variants of Mizohata--Takeuchi-type inequalities, being Theorem \ref{thm_orthonormal_MT_Rn_curve}. In Section \ref{Section_helix}, the special case of the helix is considered: the improved Wigner representation \eqref{best_possible_WR}, and its consequences, Theorem \ref{intro:Sob_MT_helix} and Theorem \ref{orthonormal_MT_helix}, are proven. Finally, in Section \ref{Section_second_marginal}, we prove an additional property of the curve-carried Wigner transform and in Section \ref{Section_bilinear_bounds} we prove the bilinear fractional integral estimates which are used in the proof of Theorem \ref{thm_orthonormal_MT_Rn_curve}.
 
 \section{Curves in \texorpdfstring{$\mathbb{R}^3$}{R3}}
\label{Section_curvesR3}
\subsection{A class of curves}
In this section, let $S$ denote a compact piece of a $C^{4}$ curve in $\mathbb{R}^3$ with arclength parametrisation given by
    $S = \{\phi(t): t \in I \}$, for some interval $I \subset \mathbb{R}$.
Assume that $S$ is well-curved, that is,
\begin{equation}
\label{well_curved}
    |\phi'(\xi) \wedge \phi''(\xi) \wedge \phi'''(\xi)| \neq 0 \qquad \text{for all } \xi \in I.
\end{equation}
Associated with the curve is its Serret-Frenet frame $(T, N, B)$, where $T(\xi), \: N(\xi), \: B(\xi)$ represent the unit tangent, normal, and binormal, respectively at some point $\phi(\xi)$ on the curve and are defined such that for all $\xi \in I$,
\begin{equation}
    T(\xi) = \phi'(\xi); \qquad N(\xi) = \frac{T'(\xi)}{\|T'(\xi)\|}; \qquad B(\xi)=T(\xi) \times N(\xi).
\end{equation}
They satisfy the following Serret-Frenet equations
\begin{align}
    T' = \kappa N \quad \quad
    N' = - \kappa T + \tau B \quad \quad
    B' = - \tau N,
\end{align}
where, 
$\kappa(\xi)$ is the curvature of $S$ at $\xi$ and
$\tau(\xi)$ is the torsion of $S$ at $\xi$. The well-curved assumption imposes that $\kappa(\xi) \neq 0$ and $\tau(\xi) \neq 0$ for all $\xi \in I$. Moreover, $\kappa(\xi) := \|T'(\xi)\| >0$.
We make the technical assumption that
\begin{equation}
\label{CONDN_3D_wedge_prod_different_pt}
\tag{VC}
   | T(\xi)\wedge N(\eta) \wedge B(\zeta)| \geq \frac{1}{2} \quad \text{for all $\xi, \eta, \zeta \in I$.}
\end{equation}
\begin{remark}
\label{rem:VC_implies_TC_NC}
The volume or frame-stability condition 
\eqref{CONDN_3D_wedge_prod_different_pt} implies the following two conditions important to the arguments of this section:
\begin{equation}
\label{eqn_tangents_condition}
\tag{TC}
    T(\xi)\cdot T(\eta) \geq \frac{1}{2} \quad \text{for all $\xi, \eta \in I$,}
\end{equation}
and
\begin{equation}
\tag{NC}
\label{eqn_normals_condition}
    N(\xi)\cdot N(\eta) \geq \frac{1}{2}\quad \text{for all $\xi, \eta \in I$.}
\end{equation}
While condition \eqref{CONDN_3D_wedge_prod_different_pt} also implies an analogous condition regarding the dot product of the binormals at different points, we will not explicitly make use of this.
We remark that the role of the constant $\frac{1}{2}$ is not important. Our applications of the Wigner representation discussed in the introduction behave suitably well under breaking the curve up into boundedly many pieces. Therefore, the importance here is that $\frac{1}{2}$ is some uniform constant strictly between $0$ and $1$.
\end{remark}
\noindent
We will later assume that 
\begin{equation}
\tag{IC}
\label{eqn_N'_dot_T_not_zero_condition}
    N'(\xi)\cdot T(\eta) \neq 0 \qquad \text{for all } \xi, \eta \in I, 
\end{equation}
which we refer to as the \textit{injectivity condition} as it is a sufficient condition to ensure the injectivity of a later defined map, see Remark \ref{rem:g_increas_bijective}. 
    We note that the Frenet formulae state $N'(\xi) \cdot T(\eta) = - \kappa(\xi) T(\xi) \cdot T(\eta) + \tau(\xi)B(\xi) \cdot T(\eta)$, hence $N'(\xi) \cdot T(\xi) = -\kappa(\xi)$. The importance here being that $N'(\xi)\cdot T(\eta)$ does not vanish when the points coincide.

\begin{remark}[Comparison with curves in two dimensions]
    Bennett, Gutiérrez, Nakamura and Oliveira \cite{BGNO} impose a similar technical restriction on the class of hypersurfaces they consider. In the two dimensional case, their condition stipulates that $N(\xi) \cdot N(\eta) \geq \frac{1}{2}$. Not only is this condition a two-dimensional analogue of condition \eqref{eqn_normals_condition}, but it also encompasses all of the two-dimensional analogues of the above conditions.
    Indeed, in two dimensions, the following two identities can be seen to hold by suitable rotations:
    \begin{equation}
        N(\xi) \cdot N(\eta) = T(\xi) \cdot T(\eta), \qquad T(\xi) \wedge N(\eta) = T(\xi) \cdot T(\eta).
    \end{equation}
    Hence, their condition is equivalent to a two dimensional analogue of condition \eqref{CONDN_3D_wedge_prod_different_pt}. Finally, in the plane, condition \eqref{eqn_N'_dot_T_not_zero_condition} is implied by condition \eqref{eqn_tangents_condition} and non-vanishing curvature (from the well-curved assumption), since 
    $N'(\xi) \cdot T(\eta) = - \kappa(\xi) T(\xi) \cdot T(\eta) \neq 0$.
\end{remark}
\subsection{A curve-carried Wigner transform}
For $g_1, g_2 \in L^2(I)$, we define the \textit{curve-carried Wigner transform}, via
\begin{equation}
    W_N(g_1, g_2)(\xi,v) = \int_{I}g_1(\eta) \overline{g_2(R_{\xi}\eta)} e^{-2\pi iv\cdot (\phi(\eta) - \phi(R_{\xi}\eta))}J(\xi,\eta)\rmd \eta,
\end{equation}
where $(\xi, v) \in \{ (\xi, v): \xi \in I, v \in \langle N(\xi) \rangle^{\perp}
    \}:= D_N.$
\noindent
For $\xi \neq \eta$, $R_{\xi}\eta$ is defined to be the unique point $\zeta \in I$ such that $\zeta \neq \eta$ and
\begin{equation}
\label{eqn_THE_defining_condition_R}
     N(\xi) \cdot (\phi(\eta) - \phi(\zeta)) =0.
\end{equation}
Moreover, $R_{\xi}\xi := \xi$ for all $\xi \in I$.
Away from $\xi = \eta$, the map $\xi \mapsto R_{\xi} \eta$ is continuously differentiable and
\begin{equation}
\label{explicit_J_formula}
        J(\xi, \eta) := \left| 
        \frac{\partial \Rmap}{\partial \xi}
        \right| = \frac{|N'(\xi)\cdot (\phi(\eta) - \phi(\zeta))|}{|N(\xi)\cdot T(\zeta)|}.
\end{equation}
This follows from a simple application of the implicit function theorem by noticing that condition \eqref{eqn_normals_condition} and the fundamental theorem of calculus together imply that $N(t) \cdot T(s) \neq 0$ for all $s, t \in I$ with $t \neq s$.

\begin{example}[The Helix]
\label{example_Helix_R3}
Let $S$ be the helix, with arclength parametrisation given by
    \begin{equation}
    \label{eqn_example_Helix_R3}
    \tag{He}
        \phi(s) = (R\sin(\lambda s), R\cos(\lambda s), m\lambda s), \qquad s \in I, \qquad \text{with } |I| \leq \frac{\pi}{3\lambda}, \qquad \text{and } \lambda = \frac{1}{(R^2 + m^2)^{1/2}}.
    \end{equation}
The Frenet normal for the helix is given by
$N(s) = (-\sin(\lambda s), -\cos(\lambda s), 0)$. 
Therefore, the defining condition \eqref{eqn_THE_defining_condition_R} becomes $R(\cos(\lambda(\xi - \eta)) = \cos(\lambda(\xi -\zeta)))$.
Since $|s - t| \leq \frac{\pi}{3\lambda}$ for all $s, t \in I$, and $\zeta \neq \eta$, we conclude that $\zeta = 2\xi - \eta$.
It is immediate that $J(\xi, \eta) = 2$ for all $\xi, \eta \in I$.
\end{example}
\noindent
The next lemma demonstrates that if there exists a $\zeta \in I$ with $\zeta \neq \eta$ satisfying condition \eqref{eqn_THE_defining_condition_R}, then $\xi$ must lie between $\eta$ and $\zeta$. We have already established this in the example of the helix where $\xi$ is the midpoint of $\eta$ and $\zeta$. This ordering underpins the well-definedness of the curve-carried Wigner transform and is central to proving many of its properties and consequently the main results of this paper. 
Though its proof is elementary, a similar argument demonstrates why other desired Wigner transforms are ill-defined, as we see in Remark \ref{eqn:the_wrong_condn_B}.
\begin{lemma}[The in-between lemma]
\label{xi_in_middle}
     Let $S$ satisfy condition \eqref{eqn_normals_condition}. Let $\xi, \eta \in I$ be distinct points and $\zeta \in I$ satisfy $\zeta \neq \eta$ and condition \eqref{eqn_THE_defining_condition_R}.
Then the function
\begin{equation}
\label{aux_function_xi_in_middle_lemma}
    f: I \to \mathbb{R}, \qquad f(t) = N(\xi)\cdot (\phi(\eta) - \phi(t))
\end{equation}
is strictly concave, with maximum attained at $t = \xi$ and zeroes $t = \eta$, $t = \zeta$.
In particular,     $$\min\{\eta,\zeta \}< \xi < \max\{\eta,\zeta \}.$$
\end{lemma}
\begin{proof}[Proof of Lemma \ref{xi_in_middle}]
Fix $\xi, \eta \in I$ and let $\zeta \in I$ satisfy $\zeta \neq \eta$ and condition \eqref{eqn_THE_defining_condition_R}.
Observe that
\begin{equation}
    f'(t)= -T(t)\cdot N(\xi), \qquad \text{and} \qquad f''(t)= -\kappa(t)N(t)\cdot N(\xi).
\end{equation}
Since for all $t \in I$ we have both $\kappa(t)>0$
and $N(t)\cdot N(\xi) \geq \frac{1}{2}$ for all $\xi, t \in I$ (from condition \eqref{eqn_normals_condition}), we also have that
\begin{equation}
    f''(t)= -\kappa(t)N(t)\cdot N(\xi) < 0.
\end{equation}
Thus $f'$ is strictly decreasing. This is to say that for $t_1, t_2 \in I$ such that $t_1 < \xi < t_2$,
\begin{equation}
    f'(t_2) < f'(\xi) < f'(t_1).
\end{equation}
Since also $f'(\xi)=0$, this implies 
\begin{equation}
    f'(t_2) < 0 < f'(t_1).
\end{equation}
Hence, $f$ is strictly increasing on $(- \infty, \xi)\cap I$ and strictly decreasing on $( \xi, \infty)\cap I$.
We note that $f(\eta)= f(\zeta)=0$, $f \in C^4$  and recall that $\eta \neq \zeta$. Hence, we conclude that $\min\{\eta,\zeta \}< \xi < \max\{\eta,\zeta \}$, as $f$ can have at most one zero on either side of $\xi$.
\end{proof}
\begin{remark}[Existence implies uniqueness]
\label{exist_implies_unique}
Let $\xi, \eta \in I$. We see that if there exists a $\zeta \in I $ with $\zeta \neq \eta$ such that
$N(\xi) \cdot (\phi(\eta) - \phi(\zeta)) =0$, then the proof above reveals it is unique.
Notice, we have used the normals condition \eqref{eqn_normals_condition} to draw this conclusion and this has the role of choosing a suitably small interval $I$ to ensure this uniqueness. As we saw in the example of the helix, there existed multiple solutions to $N(\xi)\cdot (\phi(\eta)- \phi(\zeta))=0$ for fixed $\xi, \eta \in I$ unless $I$ was suitably restricted. 
\end{remark}
\begin{remark}[Rationale for this choice of $\zeta$]
As discussed in the introduction, for fixed $\xi, \eta \in I$, there does not exist a third point $\zeta \in I$ with $\zeta \neq \eta$ such that  
\begin{equation}
\label{eqn:the_wrong_condn_B}
    B(\xi) \cdot (\phi(\eta) - \phi(\zeta)) =0.
\end{equation}
Defining an auxiliary function $f_B:I \to \mathbb{R}$ via $f_B(t) = B(\xi) \cdot (\phi(\eta) - \phi(t))$, one can run a similar argument to that of Lemma \eqref{xi_in_middle} to observe that $f_B$ is \textit{strictly monotonic} on $I \setminus \{\xi\}$ and hence has a unique zero. We refer forward to Proposition \ref{prop:no_third_pt} for the formal statement and proof.
\end{remark}
Two important properties that we expect a well-defined Wigner transform to satisfy are a certain marginal property and the conjugate-symmetric property. Governing the proofs of each of these properties is a map: $\xi \mapsto \Rmap$ and $\eta \mapsto \Rmap$ respectively.
To see that the Wigner transform is well-defined, we need to understand the domains of the maps $\xi \mapsto \Rmap$ and $\eta \mapsto \Rmap$. Explicitly, for fixed $\xi \in I$, define
       \begin{equation}
\label{actual_interval_Ixi}
       I_{\xi} :=\{\eta \in I:  N(\xi) \cdot (\phi(\eta) - \phi(\zeta))=0 \text{ for some } \zeta \neq \eta \}.
       \end{equation}
    Moreover, for fixed $\eta \in I$, define
\begin{equation}
\label{actual_interval_Ieta}
       \dot{I}_{\eta} :=\{\xi \in I:  N(\xi) \cdot (\phi(\eta) - \phi(\zeta))=0 \text{ for some } \zeta \neq \eta \}.
       \end{equation}
The following elementary lemma establishes that both domains are intervals.
\begin{lemma}[Domains of the $R$ map]
\label{lem:domains_Rmap}
Let $S$ satisfy \eqref{eqn_N'_dot_T_not_zero_condition}.
Denote $I=[a, b]$ and define $f$ as in \eqref{aux_function_xi_in_middle_lemma}. For fixed $\xi \in I$, then
\begin{equation}
\label{star_1}
\tag{$\star_1$}
    I_{\xi} = \begin{cases}
      [a', b] \qquad \text{if } f(a) < f(b) \\
      [a, b'] \qquad \text{if } f(a) > f(b) \\
      [a, b] \qquad \text{ if } f(a) = f(b),
  \end{cases}  
\end{equation}
where $a' \in [a, \xi]$ is the unique point such that $f(a') = f(b)$ and $b' \in [\xi, b]$ is the unique point such that $f(a) = f(b')$.

Now fix $\eta \in I$ and define $g: \dot{I}_{\eta} \to I$ via $g(t) = R_{t} \eta$. Then,
\begin{equation}
\label{star_2}
\tag{$\star_2$}
    \dot{I}_{\eta}= [a'', b''],
\end{equation}
where $a'', b'' \in I$ are defined such that $g(a'') = a$ and $g(b'') =b$. 
\end{lemma}
\begin{proof}We prove \eqref{star_1} and \eqref{star_2} in turn.\\
\textbf{\eqref{star_1}:}
Fix $\xi \in (a, b)$, noting that if $\xi = a$ or $b$, then $I_{\xi}= \{\xi \}$. 
    Assume that $f(a) < f(b)$ (the case $f(a) > f(b)$ is symmetric and the case $f(a) = f(b)$ is immediate).
    The existence and uniqueness of the point $a'$ follows from the concavity of $f$ (that is $f(a) < f(b) < f(\xi)$) and the intermediate value theorem.
    
To show that $[a',b] \subseteq I_{\xi}$, it suffices to consider $\eta \in (a', \xi)$, since the endpoints follow using $f(\xi)=\xi$ and $f(a')=f(b)$ and the remaining case follows by a symmetric argument. Since $f$ is strictly increasing on $(a', \xi)$, we have
    $f(b) = f(a') < f(\eta) < f(\xi).$
Hence, the intermediate value theorem implies that there exists a $\zeta \in [\xi, b]$ such that $f(\zeta) = f(\eta)$ and moreover, $\zeta \neq \eta$. Hence $\eta \in I_{\xi}$.

To see that $I_{\xi} \subseteq [a', b]$, it suffices to show that $I_{\xi}\cap [a, a') = \emptyset$. Assume for a contradiction that $\eta \in I_{\xi} \cap [a, a')$. Then there exists a $\zeta \in I$ with $\zeta \neq \eta$ such that $f(\eta) = f(\zeta)$.  Since $f$ is strictly increasing on $[a, \xi]$ and $\eta < a' \leq \xi$, observe
    $f(\zeta) = f(\eta) \leq f(a') = f(b).$
Again using $\eta < a' < \xi$, the in-between lemma dictates that $\zeta \in (\xi, b]$. Since $f$ is strictly decreasing on $(\xi, b]$, we have $f(\zeta) \geq f(b)$ and therefore a contradiction.

\textbf{\eqref{star_2}:}
We begin by showing that $g$ is strictly increasing on $I$.
By the fundamental theorem of calculus and the injectivity condition \eqref{eqn_N'_dot_T_not_zero_condition}, observe
\begin{equation}
\label{use_for_injectivity}
    N'(\xi)\cdot (\phi(\eta)- \phi(\zeta)) = \int_{\zeta}^{\eta} N'(\xi)\cdot T(t)\rmd t \neq 0.
\end{equation}
Moreover, $|N(t) \cdot T(s) |\leq 1$ for all $s, t \in I$. 
Thus, $J(\xi, \eta) \neq 0$ by 
\eqref{explicit_J_formula}. Hence, $g$ is strictly monotonic away from $\xi = \eta$. 

Now assume for a contradiction that $g$ is strictly decreasing away from $\xi = \eta$ and choose $\xi_1, \xi_2\in I$ such that $\xi_1 < g(\xi_2) \leq \xi_2 < \eta$. Then, $\xi_1 < g(\xi_2) < g(\xi_1)$. However, $g(\xi_1) \leq \xi_1$ by Lemma \ref{xi_in_middle}. Hence $g$ is strictly increasing away from $\xi= \eta$. 
Next observe that for $\xi_1, \xi_2 \in I$ satisfying $\xi_1 < \eta < \xi_2$, Lemma \ref{xi_in_middle} implies $g(\xi_1) < \xi_1 < \eta = g(\eta) = \eta < \xi_2 < g(\xi_2) $, and consequently $g$ is strictly increasing on $I$.

Recall that $g$ is continuous away from $\xi = \eta$. That $g$ is also continuous at $\xi = \eta$ is a simple application of the monotone convergence theorem and Lemma \ref{xi_in_middle}.

Since $g$ is continuous and strictly increasing on $I$ with $g(a'') = a$ and $g(b'') =b$, we conclude that $g([a'', b'']) = [a, b]$ and thus, 
$[a'', b''] \subseteq \dot{I}_{\eta}$. To conclude that $\dot{I}_{\eta} \subseteq [a'', b'']$,
assume that there exists a $\xi \in \dot{I}_{\eta}\setminus [a'', b'']$.
Consider $\xi \in [a, a'')$ as other case follows similarly. Then, Lemma \ref{xi_in_middle} implies
$g(\xi) \leq g(a'') = a$, contradicting that $g$ is strictly increasing on $I$.
\end{proof}

\begin{remark}[Bijectivity of the map $g: \dot{I}_{\eta} \to I$]
\label{rem:g_increas_bijective}
    In the above lemma, we established that the map $g:\dot{I}_{\eta} \to I$ is surjective. Injectivity follows quickly from \eqref{use_for_injectivity}.
    Indeed, let $\xi_1, \xi_2 \in I$ be distinct and assume $g(\xi_1) = g(\xi_2) = \zeta$. Then
    \begin{equation}
        N(\xi_1)\cdot (\phi(\eta)- \phi(\zeta)) = 0 \qquad \text{and} \qquad N(\xi_2)\cdot (\phi(\eta)- \phi(\zeta)) = 0,
    \end{equation}
    and $\zeta \neq \eta$.
    Now define $h: I \to \mathbb{R}$ to be such that $h(t) = N(t)\cdot (\phi(\eta)- \phi(\zeta))$.
    Since $h(\xi_1) = h(\xi_2)$, by the mean value theorem, there exists a $\xi \in I$ with $\xi_1 < \xi < \xi_2$ such that $h'(\xi)=0$.
   However, \eqref{use_for_injectivity} states that $h'(t) \neq 0$. Hence we conclude that $\xi_1 = \xi_2$. Note that we used condition \eqref{eqn_N'_dot_T_not_zero_condition} in arriving at \eqref{use_for_injectivity}; this provides the reason for its naming.
\end{remark}
\begin{proposition}[The Wigner Transform is Conjugate-Symmetric]
\label{prop:conj_sym_3D}
    The Wigner transform satisfies
    \begin{align}
        W_N(g_1, g_2)(\xi, v) = \overline{W_N(g_2, g_1)(\xi, v)},
    \end{align}
    for all $(\xi, v) \in D_N$.
    In particular, the Wigner distribution is real-valued.
\end{proposition}
\begin{proof}
Fix $\xi \in I$ and change variables via the map $\gamma_T: I \to \tilde{I}$. 
We notice that for all $\eta \in I_{\xi}$, $R_{\xi}(R_{\xi}\eta) = \eta$ due to the symmetry of the defining condition $N(\xi)\cdot (\phi(\eta) - \phi(\zeta))$ in $\eta$ and $\zeta$ and the uniqueness of $\zeta$ where it exists (see Remark \ref{exist_implies_unique}). In particular, $R_{\xi}$ is an involution and hence a bijection. 
Recall also that $N(s) \cdot T(t) \neq 0$ for all $s, t \in I$ with $s \neq t$. 
Then for all $\xi, \eta \in I$ with $\xi \neq \eta$, $\eta \mapsto \Rmap$ is continuously differentiable and
\begin{equation}
\label{second_Jacobian}
    \Delta(\xi, \eta) :=  \left| \frac{\partial \Rmap}{\partial \eta} \right| =\frac{|N(\xi)\cdot T(\eta)|}{|N(\xi)\cdot T(\Rmap)|}.
\end{equation}
Therefore, 
    \begin{align}
        W_N(g_1, g_2)(\xi, v) &=
        \int_{I_{\xi}}g_1(\eta) \overline{g_2(R_{\xi}\eta)} e^{-2\pi iv\cdot (\phi(\eta) - \phi(R_{\xi}\eta))}J(\xi,\eta)\rmd \eta \\
        &=
        \int_{I_{\xi}}g_1(R_{\xi}\zeta) \overline{g_2(\zeta)} e^{-2\pi iv\cdot (\phi(R_{\xi}\zeta) - \phi(\zeta))}J(\xi,R_{\xi}\zeta)\Delta(\xi, \zeta)\rmd \zeta.
    \end{align}
Using $R_{\xi}(R_{\xi}\eta) = \eta$, \eqref{explicit_J_formula} and \eqref{second_Jacobian} we conclude 
\begin{equation}
    \label{combine_Jac1_2}
     J(\xi, \eta) \Delta(\xi, \zeta) = J(\xi, \zeta).
\end{equation}
\end{proof}
\subsection{A Wigner Representation}
Equipped with a sensible curve-carried Wigner transform, our goal is to prove a Wigner representation for the extension operator for the curve, that is, to prove Theorem \ref{intro_main_thm} in three dimensions. 
The classical Wigner transform satisfies two marginal properties, we refer to \cite{Folland_harmonic_analysis_in_phase_space} for a full treatment. We begin with the first marginal property, given by \eqref{first_marginal_classical_WT} in the introduction, where we saw its relevance in proving a Wigner representation for the Fourier extension operator associated with the paraboloid. 

With a first marginal property for the curve-carried Wigner transform in mind, we begin by introducing an auxiliary function $F: I \times \mathbb{R}^3 \to \mathbb{R}$ defined as follows
\begin{equation}
    F(\xi, x) = \int_{I}g(\eta) \overline{g(R_{\xi}\eta)} e^{-2\pi ix\cdot (\phi(\eta) - \phi(R_{\xi}\eta))}J(\xi,\eta)\rmd \eta,
\end{equation}
so that $W_Ng$ is the restriction of $F$ to $D_N$. Explicitly, 
\begin{equation}
\label{transport_eqn_property}
    F(\xi, x) = W_Ng(\xi, P_{\langle N(\xi) \rangle^{\perp}}x),
\end{equation}
where
$P_{\langle N(\xi) \rangle^{\perp}}$ is the orthogonal projection onto $\langle N(\xi) \rangle^{\perp}$.

\begin{proposition}[The First Marginal Property]
\label{first_marg}
Let $S$ satisfy \eqref{eqn_normals_condition} and \eqref{eqn_N'_dot_T_not_zero_condition}. Then,
    \begin{equation}
    \int_{I} F(\xi,x) \rmd\xi = |E_Sg(x)|^2.
\end{equation}
\end{proposition}
\begin{proof}
As aforementioned, the true domain of integration for the Wigner transform is given by $I_{\xi}$, defined by \eqref{actual_interval_Ixi}.
Now defining
\begin{equation}
       U :=\{(\xi,\eta) \in I \times I:  N(\xi) \cdot (\phi(\eta) - \phi(\zeta))=0 \text{ for some } \zeta \neq \eta\},
       \end{equation}
we observe that
\begin{equation}
    U = \{(\xi,\eta) \in I \times I: \eta \in I_{\xi} \} = \{(\xi,\eta) \in I \times I: \xi \in \dot{I}_{\eta} \}.
\end{equation}
Hence, Fubini's theorem implies
\begin{align}
  \int_{I} F(\xi,x) \rmd \xi &=
  \int_{I} \int_{I_{\xi}} g(\eta) \overline{g(R_{\xi}\eta)} e^{-2\pi ix\cdot (\phi(\eta) - \phi(R_{\xi} \eta))}J(\xi,\eta)  \;\rmd\eta \; \rmd\xi \\
  &= \int_U g(\eta) \overline{g(R_{\xi}\eta)} e^{-2\pi ix\cdot (\phi(\eta) - \phi(R_{\xi} \eta))}J(\xi,\eta) \; \rmd\eta \; \rmd\xi \\
  &= \int_I \int_{\dot{I}_{\eta}} g(\eta) \overline{g(R_{\xi}\eta)} e^{-2\pi ix\cdot (\phi(\eta) - \phi(R_{\xi} \eta))}J(\xi,\eta) \; \rmd\xi \; \rmd\eta .
\end{align}
Changing variables via the map $\xi \mapsto \Rmap$, observe that
\begin{align}
\label{Wigner_integral_rep}
    \int_I \int_{\dot{I}_{\eta}} g(\eta) \overline{g(R_{\xi}\eta)} e^{-2\pi ix\cdot (\phi(\eta) - \phi(R_{\xi} \eta))}J(\xi,\eta) \; \rmd\xi \; \rmd\eta  &= \int_I \int_I g(\eta) \overline{g(\xi)} e^{-2\pi i x \cdot (\phi(\eta) - \phi(\xi))}  \; \rmd \xi \rmd \eta
    = |E_Sg(x)|^2.
\end{align}
Here we recall that $g: \dot{I}_{\eta} \to I$ defined by $g(t) = R_t \eta$ is bijective from Remark \ref{rem:g_increas_bijective}.
\end{proof}
We defer to Section \ref{subsec:wigner_heuristicsR3} for the statement of a second marginal property for the curve-carried Wigner transform.

Now, for functions $f: \mathbb{R}^3 \to \mathbb{C}$ and $(\xi, v) \in D_N$, define the following restricted $X$-Ray transform:
\begin{equation}
    X_N f(\xi, v) := Xf(N(\xi), v) = \int_{\langle N(\xi) \rangle} f(v + y) \rmd y.
\end{equation}
\begin{theorem}[Wigner representation for Curves in $\mathbb{R}^3$]
\label{thm:wigner_rep_3d}
Let $S$ satisfy \eqref{eqn_normals_condition} and \eqref{eqn_N'_dot_T_not_zero_condition}. Then,
\begin{equation}
    |E_Sg|^2 = X^{*}_NW_Ng,
\end{equation}
where for functions $h:D_N \to \mathbb{C}$,
\begin{equation}
    X_{N}^* h(x) = \int_I h(\xi, P_{\langle N(\xi)\rangle^{\perp}}(x)) \rmd \xi,
\end{equation}
and $P_{\langle N(\xi)\rangle^{\perp}}(x)$ denotes the orthogonal projection of $x$ onto $\langle N(\xi)\rangle^{\perp}$.
\end{theorem}
Observe that $X_N^*$ is the adjoint of $X_N$. Indeed by Fubini's theorem,
\begin{align*}
\langle X_N^*h, g \rangle &= \int_I \int_{\langle N(\xi) \rangle^{\perp}} h(\xi, v) X_Ng(\xi, v) \rmd v \rmd \xi \\
&=  \int_I \int_{\langle N(\xi) \rangle^{\perp}}\int_{\langle N(\xi) \rangle} h(\xi, v) g(v +y) \rmd y \rmd v \rmd \xi  \\
&=  \int_I \int_{\mathbb{R}^3} h(\xi, P_{\langle N(\xi) \rangle}(x)) g(x) \rmd x \rmd \xi \\
&= \int_{\mathbb{R}^3} \left(\int_I h(\xi, P_{\langle N(\xi) \rangle}(x)) \rmd \xi \right) g(x)   \rmd x\\
&= \langle h, X_Ng \rangle,
\end{align*}
for all functions $h \in L^2(D_N)$, $g \in L^2(\mathbb{R}^3)$.
\begin{proof}[Proof of Theorem \ref{thm:wigner_rep_3d}]
The Wigner representation follows quickly from the first marginal property, from \eqref{transport_eqn_property}, and from the definition of $X_N^*$:
\begin{equation}
    |E_Sg(x)|^2 = \int_I F(\xi, x) \rmd \xi = \int_I W_Ng(\xi, P_{\langle N(\xi) \rangle}(x)) \rmd \xi = X_N^*W_Ng(x),
\end{equation}
for all $x \in \mathbb{R}^3$.
\end{proof}

The surface-carried Wigner representation of \cite{BGNO} gave rise to a weighted integral representation for the Fourier extension operator associated with hypersurfaces. Strictly speaking, the analogous abstract $L^2$-duality argument does not hold here since the curve-carried Wigner distribution $W_Ng$ is not in $L^2(D_N)$, as will be seen in upcoming Remark \ref{W_N_notL2}. However $W_Ng$ is bounded by Cauchy--Schwarz and hence belongs to the natural domain of $X_N^*$.
Nevertheless, the following integral representation holds
\begin{equation}
\int_{\mathbb{R}^3}|E_Sg(x)|^2w(x) \rmd x
    =  \int_I \int_{\langle N(\xi) \rangle^{\perp}} W_Ng(\xi,v)X_Nw(\xi, v) \rmd v \rmd\xi,
\end{equation}
for all weight functions $w \colon \mathbb{R}^3 \to \mathbb{C}$.
To see this, we use Fubini's theorem, the first marginal property, the trivial decomposition $\mathbb{R}^3 = \langle N(\xi) \rangle \oplus \langle N(\xi) \rangle^{\perp}$, and the definition of the auxiliary function \eqref{transport_eqn_property} in turn to give
\begin{align}
\int_{\mathbb{R}^3}|E_Sg(x)|^2w(x) \rmd x
    &= \int_{\mathbb{R}^3} \left(\int_I F(\xi, x) \rmd \xi\right)w(x) \rmd x \\
    &=  \int_I  \int_{\langle N(\xi) \rangle^{\perp}} \int_{\langle N(\xi) \rangle}
    F(\xi, v+z) w(v + z)
     \rmd z \rmd v \rmd \xi \\
    &= \int_I  \int_{\langle N(\xi) \rangle^{\perp}}  
    W_{N}g(\xi, v) \left(\int_{\langle N(\xi) \rangle} w(v+z) \rmd z \right) \rmd v
    \rmd \xi \\
    &=  \int_I \int_{\langle N(\xi) \rangle^{\perp}} W_Ng(\xi,v)X_Nw(\xi, v) \rmd v \rmd\xi.
\end{align}
\subsection{Motivating applications via Wigner heuristics}
\label{subsec:wigner_heuristicsR3}
A further important property of the classical Wigner distribution is the second marginal property:
\begin{equation}
    \int_{\mathbb{R}^d}Wg(x, v) \rmd v = |g(x)|^2 \qquad \text{for all } x \in \mathbb{R}^n. \footnote{Again, we refer to \cite{Folland_harmonic_analysis_in_phase_space} for a discussion on such a property.}
\end{equation}
The curve-carried Wigner distribution satisfies a different marginal property:
\begin{equation}
\label{second_marginal_heuristics}
    \int_{\langle T(\xi) \rangle}W_Ng(\xi, v_t + v_b) \rmd v_t = |g(\xi)|^2, \qquad \text{for all }\xi \in I, v_b \in \langle B(\xi) \rangle.
\end{equation}
We defer to Proposition \ref{Second_marginal_property_curve_Rn} for the full statement.
One may think of the left-hand side as a `tangent marginal' in the sense that for fixed $\xi \in I$, the curve-carried Wigner distribution is integrated over $\langle T(\xi) \rangle$. Furthermore, the marginal may be expressed via another restricted $X$-Ray transform:
\begin{equation}
    \int_{\langle T(\xi) \rangle}W_Ng(\xi, v_t + v_b) \rmd v_t = X_TW_Ng(\xi, v_b) := X(W_Ng)(T(\xi), v_b), \qquad \xi \in I, v_b \in \langle B(\xi) \rangle.
\end{equation}
Importantly, we emphasise that \eqref{second_marginal_heuristics} states that this marginal is \textit{independent} of $v_b \in \langle B(\xi) \rangle$, indicating why we cannot expect a full second marginal property integrating over $\langle N(\xi) \rangle^{\perp}$.

A similar independence occurs when considering $L^2$ norms of the curve-carried Wigner distribution. Indeed, we will see shortly that
\begin{equation}
\label{L2_indep_vb}
    \int_{\langle T(\xi) \rangle}|W_Ng(\xi, v_t + v_b)|^2 \rmd v_t = \int_{I} |g(\eta)|^2|g(\Rmap)|^2 \frac{J(\xi, \eta)^2}{\tilde{J}(\xi, \eta)} \rmd \eta, \qquad \text{for all }\xi \in I, v_b \in \langle B(\xi) \rangle,
\end{equation}
with $\tilde{J}(\xi, \eta)$ to be later defined. The emphasis here is on the \textit{lack of dependence on the binormal direction}, which inspires the following remark:
\begin{remark}[The Curve-Carried Wigner Distribution is not in $L^2(D_N)$]
\label{W_N_notL2}
    We remark that \\
    ${W_Ng \notin L^2(D_N)}$. This follows quickly from \eqref{L2_indep_vb} using Tonelli's theorem.
\end{remark}

A further property, referred to as Moyal's identity and to be discussed in Section \ref{sec:orthonormal}, is again influenced by a similar lack of dependence on the binormal direction.

Now using the second marginal property, we motivate the first application of the Wigner integral representation to Mizohata--Takeuchi-type inequalities. Indeed, let $w: \mathbb{R}^3 \to \mathbb{C}$ be a non-negative weight function. Beginning from the integral representation \eqref{Wigner_integral_rep},
we express $v \in \langle T(\xi), B(\xi) \rangle$ in terms of its orthogonal projections onto $\langle T(\xi) \rangle$ and $\langle B(\xi) \rangle$. Explicitly, we write 
    $v = v_t + v_b$, where $v_t := P_{\langle T(\xi) \rangle}(v) \text{ and } v_b := P_{\langle B(\xi) \rangle}(v)$, giving
\begin{align}
\label{before_heuristics}
    \int_{\mathbb{R}^3} |E_Sg(x)|^2w(x) \rmd x 
    &=  \int_I \int_{\langle B(\xi) \rangle }\int_{\langle T(\xi) \rangle} W_Ng(\xi, v_t + v_b )X_Nw(\xi, v_t + v_b ) \rmd v_t \rmd v_b \rmd\xi.
\end{align}
Heuristically, \textit{if} the Wigner transform were non-negative, we could apply Hölder's inequality in $v_t$, implying
\begin{align}
    \int_{\mathbb{R}^3} |E_Sg(x)|^2w(x) \rmd x 
    &\leq \int_I \int_{\langle B(\xi) \rangle} \left( \int_{\langle T(\xi) \rangle} W_Ng(\xi, v_t + v_b) \rmd v_t  \right)\sup_{v_t} X_N w(\xi, v_t + v_b)     \rmd v_b\rmd\xi.
\end{align}
Consequently, the second marginal property \eqref{second_marginal_heuristics} and another application of Hölder's inequality would imply
\begin{align}
    \int_{\mathbb{R}^3} |E_Sg(x)|^2w(x) \rmd x 
    &\leq \int_I   |g(\xi)|^2 \left(\int_{\langle B(\xi) \rangle} \sup_{v_T} X_Nw(\xi, v_t + v_b)     \rmd v_b \right)\rmd\xi
    \\
    &\leq
\|g\|_{L^2(I)}^2 \sup_{\xi}\int_{\langle B(\xi) \rangle} \sup_{v_T} X_Nw(\xi, v_t + v_b)     \rmd v_b. 
\end{align}

However, we cannot expect the Wigner transform to be non-negative.
 Indeed, a well-known theorem of Hudson states that the classical Wigner transform is non-negative if and only if the input functions are Gaussians (see \cite{Folland_harmonic_analysis_in_phase_space} for further detail). However, the heuristics provide a measure of the proximity to the `conjectured' bounds as well as a guide and intuition for the following applications.

 \begin{remark}[Proximity to \eqref{MT_conj_dim_m}]
  Writing the Radon transform as in \eqref{intro:Radon_Xray}, we remark that
\begin{equation} 
\sup_{\xi} \sup_{v_t} \mathcal{R}_Tw(\xi, v_t) \leq \sup_{\xi} \int_{\langle T(\xi), N(\xi)\rangle^{\perp}} \sup_{v_t} X_N w(\xi, v_t + v_b)     \rmd v_b.
\end{equation}
Therefore, the above heuristics suggest a weaker inequality than the global Mizohata--Takeuchi inequality for curves (\eqref{MT_conj_dim_m} with $\alpha =0$). This reveals a distinction from the hypersurface case, wherein following similar heuristics leads to the false \eqref{MT_conj_dim_m} with $\alpha =0$. This difference arises from an $X$-Ray rather than Radon transform in the Wigner representation as well as the independence of the marginal on the binormal direction. Due to the presence of the same binormal independence in \eqref{L2_indep_vb}, a similar difference to the hypersurface setting can be seen in the Sobolev--Mizohata--Takeuchi results.
 \end{remark}

\textbf{The third map.}
As was the case for the conjugate-symmetric and first marginal properties, there is an important map that underpins the properties addressed here. For fixed $\xi \in I$, this third map $\gamma_T: I \to \tilde{I}$ is defined via
\begin{equation}
    \gamma_T(\eta) := T(\xi) \cdot (\phi(\eta) - \phi(\Rmap)).
\end{equation}

Due to its repeated relevance in the upcoming applications, we begin by collecting a few remarks and properties of the third map. 
First, we remark that the aforementioned tangent condition \eqref{eqn_tangents_condition} (which followed as a consequence of the condition \eqref{CONDN_3D_wedge_prod_different_pt}) will be important to ensure that this map is well-defined.
\begin{proposition}
\label{lemma_Jtilde_identity_3D}
Let $S$ satisfy \eqref{CONDN_3D_wedge_prod_different_pt}.
Fix $\xi \in I$.
Then $\gamma_T$ is continuously differentiable on $I_{\xi} \setminus \{\xi\}$ and
    \begin{equation}
        \gamma_T'(\eta) = T(\xi)\cdot \left(T(\eta) + \Delta(\xi, \eta) T(R_{\xi}\eta)\right) = \frac{B(\xi)\cdot (T(\eta) \wedge T(R_{\xi}\eta))}{N(\xi)\cdot T(R_{\xi}\eta)}.
    \end{equation}
\end{proposition}
\noindent
While we defer the proof of this proposition to Section \ref{subsubsec:understand_C}, for now we make the following remark: 
\begin{remark}[Bijectivity of the third map]
\label{bijective_third_map}
Fix $\xi \in I$. By the tangent condition \eqref{eqn_tangents_condition} and since $\Delta(\xi, \eta) \geq 0$, observe that
\begin{equation}
    \gamma'_T(\eta) = T(\xi)\cdot T(\eta) + \Delta(\xi, \eta) T(\xi) \cdot T(R_{\xi}\eta)
    \geq \frac{1 + \Delta(\xi, \eta)}{2} \geq \frac{1}{2}.
\end{equation}
This shows that $\gamma_T$ is strictly increasing on $I$ and hence bijective onto its image $\tilde{I}$.
\end{remark}
In what follows, we denote
\begin{equation}
\label{explicit_J_tilde}
    \tilde{J}(\xi, \eta) := |\gamma'_T(\eta)| = \frac{|B(\xi)\wedge T(\eta) \wedge T(R_{\xi}\eta)|}{|N(\xi)\cdot T(R_{\xi}\eta)|}.
\end{equation}

Essential to understanding the arguments of the second marginal property, Moyal's identity and the applications to come is the following lemma:
\begin{lemma}
\label{lem:gammaT_sim_a-b}
Let $S$ satisfy condition \eqref{eqn_tangents_condition}.  Then for all $\xi, \eta, \zeta \in I$, we have
    \begin{equation}
    \label{gammaT_sim_a-b}
       \frac{1}{2}|\eta - \zeta| \leq |T(\xi) \cdot (\phi(\eta)- \phi(\zeta)) | \leq |\eta - \zeta|,
    \end{equation}
    and
    \begin{equation}
    \label{phia-phib_sim_a-b}
       \frac{1}{\sqrt{2}}|\eta - \zeta| \leq |\phi(\eta) - \phi(\zeta)| \leq |\eta - \zeta|.
    \end{equation}
\end{lemma}
\begin{proof}
Fix $\xi, \eta \in I$. If $\eta = \zeta$, we are done, so assume $\eta \neq \zeta$.
First we prove the lower bound of \eqref{gammaT_sim_a-b}.
By the fundamental theorem of calculus,
    \begin{align*}
         T(\xi) \cdot(\phi(\eta) - \phi(\zeta)) =
    T(\xi) \cdot \left(\int_{\zeta}^{\eta}
    \phi'(t) \rmd t \right)=
    \int_{\zeta}^{\eta}
    T(\xi) \cdot T(t) \rmd t,
    \end{align*}
    Applying the tangent condition \eqref{eqn_tangents_condition} yields
    \begin{equation}
       \frac{1}{2}|\eta - \zeta| \leq |T(\xi) \cdot (\phi(\eta)- \phi(\zeta)) |.
    \end{equation}
Now we prove the lower bound of \eqref{phia-phib_sim_a-b} by a similar argument. Assume without loss of generality that $\xi < \eta$ and observe
\begin{align*}
    | \phi(\eta) - \phi(\zeta)|^2 &= (\phi(\eta) - \phi(\zeta))\cdot \int_{\xi}^{\eta}T(t) \rmd t =  \int_{\xi}^{\eta}(\phi(\eta) - \phi(\zeta))\cdot T(t) \rmd t \\
    &= \int_{\xi}^{\eta}\left(\int_{\xi}^{\eta}T(s) \rmd s\right)\cdot T(t) \rmd t = \int_{\xi}^{\eta}\int_{\xi}^{\eta} T(s) \cdot T(t) \rmd s\rmd t.
\end{align*}
The tangent condition \eqref{eqn_tangents_condition} implies that $| \phi(\zeta) - \phi(\eta)|^2 \geq \frac{1}{2}|\eta - \zeta|^2$.
    The reverse inequalities of both statements follow by using the trivial inequality $T(s)\cdot T(t) \leq 1$ for all $s,t \in I$ in place of \eqref{eqn_tangents_condition}.
\end{proof}

\subsection{Application 1: Weighted \texorpdfstring{$L^2$}{L2} inequalities}
In this section, we look to prove bounds reminiscent of Bennett, Gutiérrez, Nakamura and Oliveira's Sobolev Stein and Mizohata--Takeuchi inequalities for hypersurfaces \cite{BGNO}.

\begin{theorem}[$L^2$ Sobolev--Stein Inequality for Curves]
\label{corol2a_l2sobolev}
Let $S$ satisfy conditions \eqref{CONDN_3D_wedge_prod_different_pt} and \eqref{eqn_N'_dot_T_not_zero_condition}. Moreover, let $s < \frac{1}{2}$. Then, there exists a finite constant $c_s$ depending only on $s$ such that
    \begin{equation}
    \int_{\mathbb{R}^3} |E_Sg(x)|^2w(x) \rmd x 
    \leq c_s  C(\phi) \int_I \mathcal{B}_{S, 2s}(|g|^2, |g|^2)(\xi)^{\frac{1}{2}} \left( \int_{\langle B(\xi) \rangle} \left\| X_Nw(\xi, \cdot + v_b) \right\|_{\dot{H}_{v_t}^s} \rmd v_b \right)  \rmd\xi,
\end{equation}
where
\begin{equation}
    \mathcal{B}_{S, s}(g_1, g_2)(\xi) = \int_{I}\frac{g_1(\eta)\overline{g_2(R_{\xi}\eta)} J(\xi,\eta)}{|\phi(\eta)- \phi(\Rmap)|^{s}}\rmd \eta,
\end{equation}
and 
\begin{equation}
 C(\phi) := \sup_{\xi, \eta \in I}\left(\frac{|N'(\xi)\cdot (\phi(\eta) - \phi(R_{\xi}\eta))|}{|B(\xi)\wedge T(\eta) \wedge T(R_{\xi}\eta)|}\right)^{1/2}. 
 \end{equation}
\end{theorem}

\begin{proof}[Proof of Theorem \ref{corol2a_l2sobolev}]
We begin from the Wigner integral representation \eqref{Wigner_integral_rep} and decompose \\
$v = v_t + v_b$ for $v_t$ and $v_b$ the orthogonal projections onto $\langle T(\xi) \rangle$ and $\langle B(\xi) \rangle$ respectively:
\begin{align}
    \int_{\mathbb{R}^3} &|E_Sg(x)|^2w(x) \rmd x \\
    &=  \int_I \int_{\langle B(\xi) \rangle}\int_{\langle T(\xi) \rangle} W_Ng(\xi, v_t + v_b)X_Nw(\xi, v_t + v_b) \; \rmd v_t \rmd v_b \rmd\xi \\
    &= \int_I \int_{\langle B(\xi) \rangle}\int_{\langle T(\xi) \rangle} (-\Delta_{v_t})^{-s/2}(W_Ng(\xi,v_t + v_b)) \; (-\Delta_{v_t})^{s/2}(X_Nw(\xi, v_t + v_b)) \; \rmd v_t \rmd v_b \rmd\xi \\
    & \label{eqn_in_proof:sob_stein}
    \tag{$\star$}
    \leq \int_I \int_{\langle B(\xi) \rangle}
    \|W_Ng(\xi,\cdot + v_b)\|_{\dot{H}^{-s}_{v_t}} \|X_Nw(\xi, \cdot + v_b)\|_{\dot{H}^{s}_{v_t}}
    \rmd v_b
    \rmd \xi. 
\end{align}
The final line is an application of the Cauchy--Schwarz inequality, revealing the homogeneous Sobolev norms.
Fixing $\xi \in I$ and $v_b \in \langle B(\xi) \rangle$, we seek to analyse $\|W_Ng(\xi,\cdot + v_b)\|_{\dot{H}^{-s}_{v_t}}$.
First notice that for $x \in \mathbb{R}^n, v_t \in \langle T(\xi) \rangle$, we have $v_t \cdot x = v_t \cdot P_{\langle T(\xi) \rangle}(x)$ where the orthogonal projection onto $\langle T(\xi) \rangle$ is given by $P_{\langle T(\xi) \rangle}(x)= (x \cdot T(\xi))T(\xi)$. Now changing variables via the third map $\eta \mapsto T(\xi) \cdot (\phi(\eta) - \phi(\Rmap))$, observe that
\begin{align}
    W_N(g_1, g_2)(\xi, v_t + v_b) &= \int_I g_1(\eta) \overline{g_2(R_{\xi}\eta)} J(\xi,\eta) e^{-2\pi i v_b \cdot (\phi(\eta) - \phi(\Rmap))}  e^{-2\pi i v_t \cdot (\phi(\eta) - \phi(\Rmap))} \rmd \eta \\
    &= \int_{\tilde{I}} \frac{g_1(\eta(\tilde{\eta})) \overline{g_2(R_{\xi}\eta(\tilde{\eta}))} J(\xi,\eta(\tilde{\eta})) e^{-2\pi i v_b \cdot (\phi(\eta(\tilde{\eta})) - \phi(\Rmap(\tilde{\eta})))}}{\tilde{J}(\xi, \eta(\tilde{\eta}))}  e^{-2\pi i v_t \cdot \tilde{\eta}T(\xi)} \rmd \tilde{\eta} \\
    &:= \int_{\tilde{I}} G_{b}(\tilde{\eta})  e^{-2\pi i v_t \cdot \tilde{\eta}T(\xi)} \rmd \tilde{\eta}.
\end{align}

Here, for $\tilde{\eta} \in \mathbb{R}$, we have defined
\begin{equation}
\label{G_b_R3}
    G_{b}(\tilde{\eta}) := \frac{g_1(\eta(\tilde{\eta})) \overline{g_2(R_{\xi}\eta(\tilde{\eta}))} J(\xi,\eta(\tilde{\eta})) e^{-2\pi i v_b \cdot (\phi(\eta(\tilde{\eta})) - \phi(\Rmap(\tilde{\eta})))}}{\tilde{J}(\xi, \eta(\tilde{\eta}))}:= G(\tilde{\eta})e^{-2\pi i v_b \cdot (\phi(\eta(\tilde{\eta})) - \phi(\Rmap(\tilde{\eta})))},
\end{equation}
with the two notations of $G_b$ and $G$ introduced so we may simply use $|G_b(\tilde{\eta})| = |G(\tilde{\eta})|$, thereby highlighting where the independence of $v_b$ will arise.
Taking both input functions to be $g$, observe that
\begin{align*}
    \|W_Ng(\xi,\cdot + v_b)\|_{\dot{H}^{-s}_{v_t}(\langle T(\xi) \rangle)}^2 
    &= (2 \pi)^{-2s}\int_{\langle T(\xi) \rangle} \left| 
    \int_{\tilde{I}} |\tilde{\eta}|^{-s} G_{b}(\tilde{\eta})  e^{-2\pi i v_t \cdot \tilde{\eta}T(\xi)} \rmd \tilde{\eta}
    \right|^2 \rmd v_t.
\end{align*}
Now, parametrise $v_t = v_T T(\xi)$ for $v_T \in \mathbb{R}$ and notice that $v_t \cdot \tilde{\eta}T(\xi) 
= v_T \tilde{\eta}$. Then by an application of Plancherel's theorem, by \eqref{G_b_R3} and finally by a reversal of the earlier change of variables,
\begin{align*}
    \|W_Ng(\xi,\cdot + v_b)\|_{\dot{H}^{-s}(\langle T(\xi) \rangle)}^2 &=
    (2 \pi)^{-2s}\int_{\mathbb{R}} \left|
    \int_{\tilde{I}} |\tilde{\eta}|^{-s} G_{b}(\tilde{\eta})  e^{-2\pi i v_T \tilde{\eta}} \rmd \tilde{\eta} \right|^2 \rmd v_T \\
    &= (2 \pi)^{-2s} \int_{\tilde{I}} |\tilde{\eta}|^{-2s} \left| G(\tilde{\eta})
    \right|^2 \rmd \tilde{\eta} \\
     &= (2 \pi)^{-2s}
    \int_{I}\frac{|g(\eta)|^2  |g(R_{\xi}\eta)|^2 J(\xi,\eta)^2}{\tilde{J}(\xi, \eta)|T(\xi) \cdot (\phi(\eta) - \phi(\Rmap))|^{2s}}\rmd \eta.
\end{align*}
Using \eqref{explicit_J_formula} and the second identity of Lemma \ref{lemma_Jtilde_identity_3D}, yields
\begin{equation}
    \frac{J(\xi, \eta)}{\tilde{J}(\xi, \eta)}
    = \frac{|N'(\xi)\cdot (\phi(\eta) - \phi(R_{\xi}\eta))|}{|B(\xi)\wedge T(\eta) \wedge T(R_{\xi}\eta)|} 
    \leq \sup_{\xi, \eta \in I}\frac{|N'(\xi)\cdot (\phi(\eta) - \phi(R_{\xi}\eta))|}{|B(\xi)\wedge T(\eta) \wedge T(R_{\xi}\eta)|} := C(\phi)^2.
\end{equation}
 \noindent
Lemmas \ref{gammaT_sim_a-b} and \ref{phia-phib_sim_a-b} together imply that
\begin{equation}
\label{claim_application_2_proof}
    |T(\xi)\cdot (\phi(\eta)-\phi(R_{\xi}\eta)| \geq \frac{1}{2} |\eta - \Rmap| \geq \frac{1}{2}|\phi(\eta) - \phi(\Rmap)| \qquad \text{for all } \xi, \eta \in I.
\end{equation}
To summarise,
\begin{align}
    \| W_Ng(\xi,\cdot + v_b)\|_{\dot{H}^{s}_{v_t}}^2 \leq \pi^{-2s} C(\phi)^2 \mathcal{B}_{S, 2s}(|g|^2, |g|^2)(\xi).
\end{align}

Returning to \eqref{eqn_in_proof:sob_stein} and using the independence of the first homogeneous Sobolev norm on $v_b$, we conclude that
\begin{align}
    \int_{\mathbb{R}^3} &|E_Sg(x)|^2w(x) \rmd x \\
    &\leq \pi^{-s} C(\phi) \int_I \int_{\langle B(\xi) \rangle} \mathcal{B}_{S, 2s}(|g|^2, |g|^2)(\xi)^{\frac{1}{2}} \left\| X_Nw(\xi, \cdot  + v_b) \right\|_{\dot{H}_{v_t}^s}  \rmd v_b\rmd\xi \\
    &= \pi^{-s} C(\phi) \int_I \mathcal{B}_{S, 2s}(|g|^2, |g|^2)(\xi)^{\frac{1}{2}} \left( \int_{\langle B(\xi) \rangle} \left\| X_Nw(\xi, \cdot  + v_b) \right\|_{\dot{H}_{v_t}^s} \rmd v_b \right)  \rmd\xi.
\end{align}
\end{proof}

\begin{theorem}[$L^2$ Sobolev--Mizohata--Takeuchi Inequality for Curves]
\label{corol2b_l2sobolev}
Let $S$ satisfy \eqref{CONDN_3D_wedge_prod_different_pt} and \eqref{eqn_N'_dot_T_not_zero_condition}. Moreover, let $s < \frac{1}{2}$. Then there exists a finite implicit constant depending only on $s$ and the diameter of $S$ such that
    \begin{equation}
    \int_{\mathbb{R}^3} |E_Sg(x)|^2w(x) \rmd x 
    \lesssim C(\phi) \sup_{\xi \in \supp^*(g)} \left( \int_{\langle B(\xi) \rangle} \left\| X_Nw(\xi, \cdot + v_b) \right\|_{\dot{H}_{v_t}^s} \rmd v_b \right) \|g\|_{L^2(I)}^2.
\end{equation}
Here, $\supp^*(g) := \{ \xi \in I: \Rmap \in \supp(g)  \text{ for some } \eta \in \supp(g)\}$.
\end{theorem}
\begin{proof}[Proof of Theorem \ref{corol2b_l2sobolev}]
First note that if $\xi \notin \supp^*(g)$ then either $\eta \notin \supp^*(g)$ or $\Rmap \notin \supp^*(g)$, so $\mathcal{B}_{S, s}(|g|^2, |g|^2)(\xi)^{1/2}=0$.
Combining this with the above Sobolev--Stein inequality, we see that
\begin{align}
\int_{\mathbb{R}^2}|E_S g(x)|^2w(x)\rmd x
&\leq c
C(\phi)\int_{\supp^*(g)} \mathcal{B}_{S,2s}(|g|^2,|g|^2)(\xi)^{\frac{1}{2}} \left( \int_{\langle B(\xi) \rangle} \left\| X_Nw(\xi, \cdot +v_b) \right\|_{\dot{H}_{v_t}^s} \rmd v_b \right)  \rmd\xi.
\\
&\leq 
c C(\phi) \sup_{\xi \in \supp^*(g)} \left( \int_{\langle B(\xi) \rangle} \left\| X_Nw(\xi, \cdot + v_b) \right\|_{\dot{H}_{v_t}^s} \rmd v_b \right)   \|\mathcal{B}_{S,2s}(|g|^2,|g|^2)\|_{L^{1/2}(I)}^{1/2}
\end{align}
Now, we use the bilinear fractional integral estimate of Theorem \ref{bilinear_frac_higher_dim} to come, which states that for $s < 1$,
\begin{equation}
    \|\mathcal{B}_{S,s}(g_1,g_2)\|_{L^{1/2}(S)}\lesssim \|g_1\|_{L^1(S)}\|g_2\|_{L^1(S)}
\end{equation}
with uniform implicit constant.
Hence, 
\begin{equation}
    \|\mathcal{B}_{S,2s}(|g|^2,|g|^2)\|_{L^{1/2}(I)}\lesssim 
    \|g\|_{L^2(I)}^4,
\end{equation}
as required.

\end{proof}

\textbf{Understanding $C(\phi)$.}
\label{subsubsec:understand_C}
We wish to interpret the constant $C(\phi)$ which governs both the Sobolev--Stein and Sobolev--Mizohata--Takeuchi inequalities.
First, we prove Proposition \ref{lemma_Jtilde_identity_3D}, beginning with a remark on $\Delta(\xi, \eta)$.
\begin{remark}
\label{corol_injectivity_of_delta}
Fix $\xi\in I$. We claim 
\begin{equation}
 \frac{\partial R_{\xi}\eta}{\partial \eta} = - \Delta(\xi, \eta).
\end{equation}
Observe that for every $s, \xi \in I$,
\begin{equation}
     N(\xi)\cdot T(s) = \int^{s}_{\xi} \kappa(t)N(\xi)\cdot N(t) \rmd t,
\end{equation}
whose integrand is strictly positive on $I$ by the normals condition \eqref{eqn_normals_condition} and using $\kappa(\xi) >0$ from the well-curved hypothesis. In particular, writing
\begin{equation}
    m := \inf_{t \in I} \kappa(t) N(\xi) \cdot N(t) \qquad M := \sup_{t \in I} \kappa(t) N(\xi) \cdot N(t),
\end{equation}
we see that $0 < \frac{1}{2} \inf_{t \in I} \kappa(t) \leq m \leq M \leq \sup_{t \in I}\kappa(t)$.
Moreover,
\begin{equation}
\label{in_remark_bound_inj_delta}
    m |s- \xi| \leq |N(\xi) \cdot T(s)| \leq M|s - \xi| \qquad \text{and} \qquad \sgn\left(N(\xi) \cdot T(s)\right) = \sgn(s- \xi).
\end{equation}
In Lemma \ref{xi_in_middle} we saw that $\xi$ lies between $\eta$ and $\zeta= \Rmap$, so $\eta - \xi$ and $\zeta - \xi$ have opposite signs and hence $\frac{\eta - \xi}{\zeta - \xi}<0$. Therefore, using 
\eqref{in_remark_bound_inj_delta} and that $\frac{\eta - \xi}{\zeta - \xi}= - \left|\frac{\eta - \xi}{\zeta - \xi}\right|$, we obtain
\begin{equation}
   \frac{M(\eta-\xi)}{m(\zeta-\xi)} \leq \frac{N(\xi)\cdot T(\eta)}{N(\xi)\cdot T(\zeta)} \leq \frac{m(\eta-\xi)}{M(\zeta-\xi)} < 0.
\end{equation} 
\end{remark}
\begin{proof}[Proof of Proposition \ref{lemma_Jtilde_identity_3D}]
Fix $\xi \in I$.
By \eqref{second_Jacobian} we see that the function $\gamma_T$ is
 continuously differentiable away from $\xi = \eta$. Moreover,
    \begin{align}
     \tilde{J}(\xi, \eta) 
      &= 
     \left |T(\xi)\cdot \left(T(\eta) - \frac{\rmd R_{\xi}\eta}{\rmd \eta} T(R_{\xi}\eta)\right)\right| = \left|
     T(\xi)\cdot \left(T(\eta) + \Delta(\xi, \eta) T(R_{\xi}\eta)\right) \right|.
\end{align}

The final equality follows from Remark \ref{corol_injectivity_of_delta}.
    Seeking to prove the second identity, we use the explicit expression for $\Delta(\xi, \eta)$, that is \eqref{second_Jacobian}, to give 
    \begin{equation}
     T(\xi)\cdot \left(T(\eta) + \Delta(\xi, \eta) T(R_{\xi}\eta)\right) =
     \frac{T(\xi)\cdot \left((N(\xi)\cdot T(R_{\xi}\eta))T(\eta) - (N(\xi)\cdot T(\eta)) T(R_{\xi}\eta)\right)}{N(\xi)\cdot T(R_{\xi}\eta)}.
\end{equation}
Turning our attention to the numerator, we observe
\begin{align}
T(\xi)\cdot &((N(\xi)\cdot T(R_{\xi}\eta))T(\eta) - (N(\xi)\cdot T(\eta)) T(R_{\xi}\eta)) \\
&=
    (N(\xi)\cdot T(R_{\xi}\eta))(T(\xi)\cdot T(\eta)) - (N(\xi)\cdot T(\eta)) (T(\xi)\cdot T(R_{\xi}\eta)) \\
    &= (N(\xi) \times T(\xi))\cdot (T(R_{\xi}\eta) \times T(\eta)) \\
    &=B(\xi)\cdot (T(\eta) \times T(R_{\xi}\eta)),
\end{align}
thereby concluding the proof.
\end{proof}
Now, we turn to controlling $C(\phi)$, given by
\begin{equation}
    C(\phi) = \sup_{\xi, \eta} \frac{J(\xi, \eta)}{\tilde{J}(\xi, \eta)} = \sup_{\xi, \eta \in I}\frac{|N'(\xi)\cdot (\phi(\eta) - \phi(R_{\xi}\eta))|}{|B(\xi)\wedge (T(\eta) \wedge T(R_{\xi}\eta))|}.
\end{equation}

We endeavour to control $C(\phi)$ by a quantity which both removes the singularity at $\xi = \eta = \zeta$ and is 
more intuitive to interpret. This is precisely the lemma below.
\begin{lemma}[Controlling $C(\phi)$]
\label{lem:control_Cphi_3d}
Let $S$ satisfy \eqref{CONDN_3D_wedge_prod_different_pt} and \eqref{eqn_N'_dot_T_not_zero_condition}. Then
    \begin{equation}
        C(\phi) \leq 2\sup_{\xi, \eta \in I} \frac{\kappa(\xi) + |\tau(\xi)|}{\kappa(\eta)}.
    \end{equation}
\end{lemma}
\begin{proof}
    Fix $\xi, \eta \in I$ and let $\zeta = \Rmap$. Since $B(\xi) \cdot (T(\zeta) \wedge T(\zeta))=0$, observe that by an application of the generalised mean value theorem there exists a $c$ between $\eta$ and $\zeta$ such that
    \begin{align}\frac{N'(\xi)\cdot(\phi(\eta)- \phi(\zeta))}{B(\xi) \cdot (T(\eta) \wedge T(\zeta))} 
        := \frac{f(\eta) - f(\zeta)}{g(\eta)}= \frac{f(\eta) - f(\zeta)}{g(\eta)- g(\zeta)}
        = \frac{N'(\xi)\cdot T(c)}{\kappa(c)B(\xi) \cdot (N(c) \wedge T(\zeta))}.
    \end{align}
    Then using assumption \eqref{CONDN_3D_wedge_prod_different_pt}, the Frenet formulae, and finally the triangle inequality with $|T(\xi)\cdot T(c)| \leq 1$ and $|B(\xi)\cdot T(c)| \leq 1$ we see that
    \begin{align}
\frac{|N'(\xi)\cdot(\phi(\eta)- \phi(\zeta))|}{|B(\xi) \cdot (T(\eta) \wedge T(\zeta))|} 
        = \frac{|N'(\xi)\cdot T(c)|}{\kappa(c)|B(\xi) \wedge N(c) \wedge T(\zeta)|} \leq 2\frac{|- \kappa(\xi)T(\xi)\cdot T(c) + \tau(\xi)B(\xi)\cdot T(c)|}{\kappa(c)} \leq 2\frac{\kappa(\xi) + |\tau(\xi)|}{\kappa(c)}.
    \end{align}
    Taking suprema we conclude.
\end{proof}
\begin{example}[The Helix]
We return to the example in which $S$ is the arclength parametrised helix with $|I| \leq\frac{\pi}{3k}$. A short computation reveals that
\begin{equation}
    |I| \leq\frac{\pi}{3k} \iff |T(\xi) \wedge N(\eta) \wedge B(\zeta)| \geq 1/2 \text{ for all } \xi, \eta, \zeta \in I.
\end{equation}
Hence $T(\xi) \cdot T(\eta) \geq 1/2 \text{ for all } \xi, \eta, \in I$ as seen in \eqref{rem:VC_implies_TC_NC}. 
In the case of the helix, $\Delta(\xi, \eta) \equiv 1$, implying
\begin{equation}
        \tilde{J}(\xi, \eta) = T(\xi)\cdot T(\eta) + T(\xi) \cdot  T(R_{\xi}\eta) \geq 1.
    \end{equation}
Using the trivial inequality $T(s) \cdot T(t) \leq 1$ for all $t, s \in I$, we conclude that
\begin{equation}
    1 \leq \tilde{J}(\xi, \eta) \leq 2 \qquad \text{for all } \xi, \eta \in I.
\end{equation}
Recalling that $J(\xi, \eta) =2$ for all $\xi, \eta \in I$ for the helix, we conclude that 
\begin{equation}
    \frac{J(\xi, \eta)}{\tilde{J}(\xi, \eta)}\leq 2 \qquad \text{for all $\xi, \eta \in I$},
\end{equation}
and therefore $C(\phi) \leq 2$.
\end{example}

\begin{remark}[Scale invariance]
\label{rem:scale_invar_const_3d}
    We remark that the quantity 
    \begin{equation}
        \tilde{C}(\phi) = \sup_{\xi, \eta \in I}\frac{\kappa(\xi) + |\tau(\xi)|}{\kappa(\eta)}
    \end{equation}
    is scale invariant. To make this precise, let $\lambda S = \{ \lambda u: u \in S\}$. Then, $\phi_{\lambda}(\xi) = \lambda \phi(\frac{\xi}{\lambda})$ is the arclength parametrisation of $\lambda S$ since $|\phi_{\lambda}'(\xi)| = |\phi'(\xi/ \lambda)|=1$. Moreover, $\kappa_{\lambda}(\xi)= \frac{1}{\lambda}\kappa(\xi)$ and ${\tau_{\lambda}(\xi)= \frac{1}{\lambda}\tau(\xi)}$. Hence, $\tilde{C}(\phi_{\lambda})= \tilde{C}(\phi)$. This means that the Sobolev--Stein and Sobolev--Mizohata--Takeuchi inequalities for curves do not depend on any of the Frenet curvatures of $S$ in absolute terms, a property we sought since the Stein and Mizohata--Takeuchi inequalities hold true when $S$ is a straight line.
\end{remark}
  \section{Curves in \texorpdfstring{$\mathbb{R}^n$}{Rn}}
  \label{Section_curvesRn}
\subsection{Preliminaries and a class of curves}
Unless otherwise stated, in this section, assume that $S$ is a compact piece of a $C^{\alpha}$ curve in $\mathbb{R}^n$ for some $\alpha \geq n + 1$. Let $\phi: I \to \mathbb{R}^n$ be the arclength parametrisation of $S$ for some interval $I \subset \mathbb{R}$ and assume $S$ is well-curved (as defined by \eqref{well_curved_n}). 

\medskip
\textbf{The Frenet frame.}
Essential to the study of curves is the Frenet frame. We begin by briefly discussing the preliminaries and notation surrounding the Frenet frame and formulae here, referring to \cite{kuhnel} for a complete treatment.

The \textit{Frenet frame} of the curve $S$, denoted by $(e_1, \dots, e_n)$, is the collection of functions ${e_1, ..., e_n:I \to \mathbb{R}^n}$ that uniquely satisfy the following for every $t \in I$:
\begin{enumerate}
\item{\label{eqn:Frenet_frame_condn1}$e_1(t), e_2(t), ..., e_n(t)$ are orthonormal and positively oriented;}
\item{\label{eqn:Frenet_frame_condn2}$\langle e_1(t), ..., e_k(t)\rangle = \langle \phi'(t), ..., \phi^{(k)}(t)\rangle$ for each $k \in \{1, ..., n-1\}$;} 
    \item{\label{eqn:Frenet_frame_condn3} $ \phi^{(k)}(t) \cdot e_k(t) >0$ for each $k \in \{1, ..., n-1\}$.}
\end{enumerate}
Here, $e_1$ is the unit tangent, $e_2$, the unit normal, $e_3$ the unit binormal, and so on, where we adopted the standard notation of $T, N, B$ for curves in $\mathbb{R}^3$.
Moreover $e_j \in C^{\alpha-j}$ for all $j \in \{1, \dots, n-1\}$ and $e_n \in C^{\alpha - n+1}$.

Essential to this work are the \textit{Frenet equations} or \textit{formulae}.
\begin{theorem}[Frenet equations in $\mathbb{R}^n$]
\label{THM_Frenet_eqns}
    There exist functions ${\kappa_j\colon I \to \mathbb{R}}$ with $\kappa_{j}>0$ for ${j \in \{1, ..., n-2\}}$ and $\kappa_j \in C^{n-1 -j}$ such that for all $j \in \{1, \dots n\}$,
\begin{equation}
    e_j'(t) = -\kappa_{j-1}(t)e_{j-1}(t) + \kappa_{j}e_{j+1}(t).
\end{equation} 
Moreover, $\kappa_j \colon= \langle e_j', e_{j+1} \rangle = - \langle e_j, e_{j+1}' \rangle$.
\end{theorem}
\noindent
Here $\kappa_j$ is the \textit{$j$-th Frenet curvature} and the final Frenet-curvature $\kappa_{n-1}$ is sometimes referred to as the \textit{torsion}. The well-curved hypothesis also imposes that $\kappa_{n-1}(t) \neq 0$ for all $t \in I$.

\medskip
\textbf{Derivatives of Frenet vectors.}
Generalising the work of three dimensions to $n$ dimensions requires an understanding of higher order derivatives of Frenet vectors. In three dimensions, we needed at most first order derivatives, described by the Frenet equations. Now, we prove what we will require of a higher order analogue of the Frenet equations.

\begin{lemma}[Derivatives of the Frenet vectors]
\label{lemma_derivatives_frenet_vectors}
Fix $t \in I$.
    For $i, j \in \mathbb{N}$, with $j \leq n$ and $j + i \leq \alpha$,
    \begin{equation}
        e_{j}^{(i)}(t) \in \langle \{ e_m(t): m \in \mathbb{N}, \max \{ 1, j - i \}\leq m \leq  \min\{j+ i , n \}\} \rangle.
    \end{equation}
In particular, writing $e_{j}^{(i)}(t)$ in terms of its coordinates,
\begin{equation}
    e_{j}^{(i)}(t) = \sum_{m =j-i}^{j+i}c_m(t)e_m(t),
\end{equation}
we have
\begin{equation}
   c_{j-i}(t)= (-1)^{i}\prod_{l=j-i}^{j-1}\kappa_l(t) \qquad \text{and} \qquad
   c_{j+i}(t)= \prod_{l=j}^{j+i-1}\kappa_l(t).
\end{equation}
\end{lemma}
The way in which we will most often make use of this lemma is via the following corollary:
\begin{corollary}
\label{corol_derivatives_frenet_vectors}
Fix $i, j \in \mathbb{N}$, with $j \leq n$ and $j + i \leq \alpha$, and $t \in I$. Then for all $m \in \mathbb{N}$ such that $j-i < m < j+i$,
\begin{equation}
    e_j^{(i)}(t)\cdot e_{m}(t)=0.
\end{equation}
Moreover,
\begin{equation}
    e_j^{(i)}(t)\cdot e_{j-i}(t)= (-1)^{i}\prod_{m=j-i}^{j-1}\kappa_m(t), \qquad \text{and} \qquad
        e_j^{(i)}(t)\cdot e_{j+i}(t)= \prod_{m=j}^{j+i-1}\kappa_m(t).
\end{equation}
\end{corollary}
\noindent
    The corollary follows immediately from Lemma \ref{lemma_derivatives_frenet_vectors} by using the orthonormality of the Frenet vectors; explicitly it uses that $e_j^{(i)}(t)\cdot e_m(t) = c_m(t)$.
\begin{proof}[Proof of Lemma \ref{lemma_derivatives_frenet_vectors}]
Fix $j \in \{1, ..., n\}$.
    We prove Lemma \ref{lemma_derivatives_frenet_vectors} by induction on $i$.
    We begin with the base case of $i=1$. Observe that for all $t \in I$
    \begin{equation}
        e_{j}'(t) = -\kappa_{j-1}(t) e_{j-1}(t) + \kappa_{j}(t) e_{j+1}(t).
    \end{equation}
    Hence, writing 
    $e_j'(t)$ in terms of its coordinates, we see immediately from the above that
    $c_{j-1} = - \kappa_{j-1}(t)$ and $c_{j+1} = \kappa_{j}(t)$.
\\\\
Now suppose the claim holds for some $i \geq 1$, that is we assume that for all $t \in I$,
    \begin{equation}
        e_{j}^{(i)}(t) \in \langle \{ e_m(t): m \in \mathbb{N}, \max \{ 1, j - i \}\leq m \leq  \min\{j+ i , n \} \} \rangle.
    \end{equation}
Hence, we may write $e_{j}^{(i)}(t)$ in terms of its coordinates,
\begin{equation}
    e_{j}^{(i)}(t) = \sum_{m =j-i}^{j+i}c_m(t)e_m(t),
\end{equation}
and in particular, we assume that  
\begin{equation}
   c_{j-i}(t)= (-1)^{i}\prod_{l=j-i}^{j-1}\kappa_l(t) \qquad \text{and} \qquad
   c_{j+i}(t)= \prod_{l=j}^{j+i-1}\kappa_l(t).
\end{equation}
Let $c_m =0$ for all $m$ such that either $m> j+i$ or $m < j-i$.
The product rule and the Frenet equations give
\begin{align}
    e_{j}^{(i+1)}(t) 
    &= \sum_{m =j-i}^{j+i}c_m'(t)e_m(t) + c_m(t)e_m'(t) \\
    &= \sum_{m =j-i}^{j+i}c_m'(t)e_m(t) + (-\kappa_{m-1}(t)e_{m-1}(t)+ \kappa_{m}(t)e_{m+1}(t))c_m(t) \\
    &= \label{return_to_for_computing_deriv}\left(\sum_{m =j-i}^{j+i}(c_m'(t) -\kappa_{m}(t)c_{m+1}(t)+ \kappa_{m-1}(t)c_{m-1}(t))e_m(t)\right) 
    \\ & \qquad \qquad \qquad \qquad -\kappa_{j-i-1}(t)c_{j-i}(t)e_{j-i-1}(t)+ \kappa_{j+i}(t)c_{j+i}(t)e_{j+i+1}(t) \\
    &= \sum_{m =j-(i+1)}^{j+(i+1)}d_m(t)e_m(t). 
\end{align}
Here, we have defined
\begin{equation}
    \begin{cases}
        d_{j-i-1}(t) := - \kappa_{j-i -1}(t) c_{j-i}(t) = (-1)^{i+1}\prod_{l=j-(i+1)}^{j-1}\kappa_l(t); \\
        d_{m}(t) := c_m'(t) -\kappa_{m}(t)c_{m+1}(t)+ \kappa_{m-1}(t)c_{m-1}(t) \qquad \text{for all $j-i\leq m \leq j+i$};
        \\
        d_{j+i+1}(t) :=  \kappa_{j+i}(t) c_{j+i}(t) = \prod_{l=j}^{j+i}\kappa_l(t),
    \end{cases}
\end{equation}
where the second equalities in the first and third expressions use the assumption.
Thus, 
\begin{equation}
    e_j^{i+1}(t) \in 
    \langle \{ e_m(t): m \in \mathbb{N}, 1 \leq j - (i+1)\leq m \leq j+ (i+1) \leq n \} \rangle.
\end{equation}
This concludes the induction.
\end{proof}

\medskip
\textbf{Further assumptions on $S$.}
Using Lemma \ref{lemma_derivatives_frenet_vectors} and anticipating the work to come, we take $\alpha = n + 2k -1$ for the remainder of this section unless otherwise specified. This regularity assumption is to ensure the continuity of all of the derivatives of the Frenet vectors that will appear, and it aligns with three dimensions where $\alpha = 4$ was imposed.

As in three dimensions, $S$ will satisfy certain technical conditions that amount to restricting the length of the interval $I$. When $k=1$, sufficient technical conditions for the Wigner transform $W_2$ to be well-defined are precisely $n$-dimensional analogues of those imposed for curves in three dimensions. That is, \eqref{eqn_normals_condition} and \eqref{eqn_N'_dot_T_not_zero_condition} stated in $n$ dimensions:
\begin{equation}
\label{NC_ndim}
\tag{NCn}
    e_2(\xi)\cdot e_2(\eta) \geq \frac{1}{2} \qquad \text{for all } \xi, \eta \in I,
\end{equation}
and
\begin{equation}
\label{IC_ndim}
\tag{ICn}
    e_2'(\xi)\cdot e_1(\eta) \neq 0 \qquad \text{for all } \xi, \eta \in I.
\end{equation}
For the sake of exposition, in the case of $k>1$, we will state the technical conditions when they become relevant. 
\subsection{A curve-carried Wigner transform}
 For $g_1, g_2 \in L^2(I)$, and $k \in \mathbb{N}$ such that $2k \leq n$, we define the \textit{$k^{\text{th}}$ curve-carried Wigner transform}, via
\begin{equation}
    W_{2k}(g_1, g_2)(\xi,v) = \int_{I}g_1(\eta) \overline{g_2(R_{\xi}\eta)} e^{-2\pi iv\cdot (\phi(\eta) - \phi(R_{\xi}\eta))}J(\xi,\eta)\rmd \eta.
\end{equation}
Here, $(\xi, v) \in \{ (\xi, v): \xi \in I, v \in \langle e_{2k}(\xi) \rangle^{\perp}
    \}:= D_{2k}$, the \textit{$k^{\text{th}}$-Frenet line bundle}.
For $\xi \neq \eta$, $R_{\xi}\eta$ is defined to be the unique point $\zeta \in I$ such that $\zeta \neq \eta$ and
\begin{equation}
\label{eqn_THE_defining_condition_R_n}
    e_{2k}(\xi)\cdot(\phi(\eta) - \phi(\zeta)) =0.
\end{equation}
Moreover, $R_{\xi}\xi := \xi$ for all $\xi \in I$.
We impose the following technical assumption on $S$:
\begin{equation}
\tag{TA1}
\label{eqn_ej_dot_T_not_zero_condition}
    e_{2k}(\xi)\cdot e_1(\eta) \neq 0 \qquad \text{for all } \xi, \eta \in I \text{ such that } \xi \neq \eta.
\end{equation}
Therefore, as in the three dimensional case, we have an explicit formula for $J(\xi, \eta)$ away from $\xi = \eta$ via the implicit function theorem:
  \begin{equation}
    \label{explicit_J_formula_n}
        J(\xi, \eta) := \left| \frac{\partial \Rmap}{\partial \xi}\right|=\frac{|e_{2k}'(\xi)\cdot (\phi(\eta) - \phi(\zeta))|}{|e_{2k}(\xi)\cdot e_1(\zeta)|}.
    \end{equation}

    \begin{remark}[The first technical assumption \eqref{eqn_ej_dot_T_not_zero_condition}]
        In three dimensions, \eqref{eqn_ej_dot_T_not_zero_condition} was implied by condition \eqref{eqn_normals_condition}. The same is true in the case of $k=1$ with \eqref{eqn_ej_dot_T_not_zero_condition}. However this does not hold for $k>1$ so we have imposed \eqref{eqn_ej_dot_T_not_zero_condition} to ensure the finiteness of $J(\xi, \eta)$.
    \end{remark}

\begin{example}[The helix in $\mathbb{R}^3$]
    We delay details of the $n$ dimensional helix to Section \ref{Section_helix}, noting only that for a suitably small piece of the helix in $n$ dimensions, we have $R_{\xi}\eta = 2\xi - \eta$ for all $k \in \mathbb{N}$ with $2k \leq n$. Moreover, $J(\xi, \eta) =2$ for all $\xi, \eta \in I$.
\end{example}

We require a generalisation of the in-between lemma to arbitrary dimensions. 
The proof is similar to its three dimensional analogue when $k=1$, however differs for larger $k$. We present both proofs as this difference is a feature throughout.
\begin{lemma}[In-between Lemma when $k=1$]
\label{xi_in_middle_n_k1}
     Let $S$ satisfy \eqref{NC_ndim}. Moreover, let $\xi, \eta, \zeta \in I$ be distinct points such that $\zeta = \Rmap$. Define the function
\begin{equation}
\label{aux_function_xi_in_middle_lemma_n_k1}
    f\colon I \to \mathbb{R}, \qquad f(t) = e_2(\xi)\cdot (\phi(\eta) - \phi(t)).
\end{equation}
Then $f$ is strictly concave with maximum attained at $t = \xi$ and zeroes $t = \eta$, $t = \zeta$.
Moreover,
$$\min\{\eta,\zeta \}< \xi < \max\{\eta,\zeta \}.$$
\end{lemma}
\begin{proof}[Proof of Lemma \ref{xi_in_middle_n_k1}]
Fix $\xi, \eta \in I$, with $\zeta = \Rmap$.
Observe that
\begin{equation}
    f'(t)= -e_1(t)\cdot e_2(\xi), \qquad \text{and} \qquad f''(t)= -\kappa_1(t)e_2(t)\cdot e_2(\xi).
\end{equation}
Since for all $t \in I$ we have both $\kappa_1(t)>0$ and $e_2(t)\cdot e_2(\xi) \geq \frac{1}{2}$ (from the normal condition \eqref{NC_ndim}), we also have that
\begin{equation}
    f''(t)= -\kappa_1(t)e_2(t)\cdot e_2(\xi) < 0.
\end{equation}
Hence, we conclude as in three dimensions using that $f'(\xi) =0$.
\end{proof}

    For $k>1$, observe that by the Frenet formulae,
    \begin{equation}
        e_{2k}'(\xi)\cdot e_1(\xi)= (-\kappa_{2k-1}(\xi)e_{2k-1}(\xi)+ \kappa_{2k}(\xi)e_{2k+1}(\xi))\cdot e_1(\xi)=0.
    \end{equation}
    This holds by the orthonormality of the Frenet vectors since $2k-1 \neq 1$ and $2k+1 \neq 1$. In fact, as seen by Corollary \ref{corol_derivatives_frenet_vectors}, 
    \begin{equation}
        e_{2k}^{(i)}(\xi)\cdot e_1(\xi)=0
    \end{equation}
    for all $i < 2k-1$. Therefore, to run a similar argument for the general case, we require $2k-1$ derivatives followed by an analysis of the sign of the result. When $k=1$, the normal condition \eqref{NC_ndim} enables us to conclude, but for $k>1$ many more terms are accrued from the product rule and Frenet formulae. Therefore, a Taylor expansion about the point $\xi$ up to only the necessary terms, and the following assumption (replacing \eqref{NC_ndim}) appear to be the easiest way to conclude. Indeed, we impose that the length of the interval $I$ is sufficiently small, as we may, such that
\begin{equation}
\tag{TA2}
\label{eqn_taylor_condn}
   |\xi - \eta| \left| e_{2k}(\xi) \cdot e_1^{(2k)}(c)\right| < 2k \prod_{m=1}^{2k-1}|\kappa_m(\xi)| \qquad \text{for all } \xi, \eta, c \in I \text{ for which } |\xi - c|\leq |\xi - \eta|.
\end{equation}
\begin{remark}[Interpreting the second technical assumption \eqref{eqn_taylor_condn}]
    We emphasise that if the length of $I$ is sufficiently small, then condition \eqref{eqn_taylor_condn} holds.
    Indeed, the well-curved hypothesis and the continuity of the Frenet curvatures together imply that
    \begin{equation}
        \prod_{m=1}^{2k-1}|\kappa_m(\xi)| \geq C >0 \qquad \text{for all } \xi \in I.
    \end{equation}
    Moreover, Lemma \ref{lemma_derivatives_frenet_vectors} implies the continuity and hence boundedness of $\left| e_{2k}(\xi) \cdot e_1^{(2k)}(c)\right|$ for all $\xi, c \in I$.
    To further interpret the condition, we comment that when $c = \xi$, a careful iterative computation making use of \eqref{return_to_for_computing_deriv} and Corollary \ref{corol_derivatives_frenet_vectors} reveals that for all $\xi \in I$,
    \begin{equation}
        e_{2k}(\xi)\cdot e_1^{(2k)}(\xi) = \mathlarger{\sum_{j=1}^{2k-1}} (2k-j)\kappa_j'(\xi)\prod_{\substack{m=1 \\ m \neq j}}^{2k-1}\kappa_m(\xi).
    \end{equation}
     In the case of the helix, $e_{2k}(\xi)\cdot e_1^{(2k)}(\xi) =0$ for all $\xi \in I$ due to the constant Frenet curvatures. Moreover, due to Lemma \ref{corol_derivatives_frenet_vectors}, the condition may be immediately reinterpreted as
     \begin{equation}
         |\xi - \eta| \left| e_{2k}(\xi) \cdot e_1^{(2k)}(c)\right| < 2k |e_{2k}(\xi) \cdot e_1^{(2k-1)}(\xi)| \qquad \text{for all } \xi, \eta, c \in I \text{ for which } |\xi - c|\leq |\xi - \eta|.
     \end{equation}
    Condition \eqref{eqn_taylor_condn} is used in the subsequent lemma to ensure that the remainder of the Taylor expansion of $g_e(t) = e_{2k}(\xi)\cdot e_1(t)$ about the point $\xi \in I$ is suitably small. The condition quantifies how small the remainder must be in terms of the length of the interval $I$ and products of Frenet curvatures and their derivatives at different points on the curve.
\end{remark}

\begin{lemma}[In-between Lemma]
\label{xi_in_middle_n}
     Let $S$ satisfy \eqref{eqn_taylor_condn}. Moreover, let $\xi, \eta, \zeta \in I$ be distinct points such that $\zeta = \Rmap$. Define the function
\begin{equation}
\label{aux_function_xi_in_middle_lemma_n}
    f_e: I \to \mathbb{R}, \qquad f_e(t) = e_{2k}(\xi)\cdot (\phi(\eta) - \phi(t)).
\end{equation}
Then $f_e$ is strictly concave (or strictly convex) with maximum (or minimum) attained at $t = \xi$ and zeroes $t = \eta$, $t = \zeta$. \footnote{The distinction between concavity or convexity is determined by the sign of the final Frenet curvature.}
Furthermore,
$$\min\{\eta,\zeta \}< \xi < \max\{\eta,\zeta \}.$$
\end{lemma}
\begin{proof}[Proof of Lemma \ref{xi_in_middle_n}]
    Let $\xi, \eta, \zeta \in I$ be distinct points with that $\zeta = \Rmap$.
    Define the additional auxiliary function $g_e:I \to \mathbb{R}$ via
    \begin{equation}
    g_e(t):= - f_e'(t) = e_{2k}(\xi)\cdot e_1(t).
\end{equation}
Corollary \ref{corol_derivatives_frenet_vectors} implies that
\begin{equation}
g_e^{(i)}(\xi)= e_{2k}(\xi) \cdot e_{1}^{(i)}(\xi)=
    \begin{cases}
        0 \qquad \text{if $i <2k$}\\
        \prod_{m=1}^{2k}\kappa_m(\xi) \qquad \text{if $i =2k$}.
    \end{cases}
\end{equation}
Therefore, a Taylor expansion of $g_e$ about the point $\xi \in I$ gives
\begin{align*}
    g_e(t)&= \frac{(t- \xi)^{2k-1}}{(2k - 1)!} \prod_{m=1}^{2k-1}\kappa_m(\xi) + \frac{g_e^{(2k)}(c)(t-\xi)^{2k}}{(2k)!}\\
    &=(t- \xi)^{2k-1}\left(\frac{2k\prod_{m=1}^{2k-1}\kappa_m(\xi) + g_e^{(2k)}(c)(t- \xi)}{(2k)!}\right),
\end{align*}
for some $c$ between $t$ and $\xi$.
Condition \eqref{eqn_taylor_condn} implies
\begin{equation}
    |g_e^{(2k)}(c)| = |e_{2k}(\xi) \cdot e_1^{(2k)}(c)| <  \frac{2k\prod_{m=1}^{2k-1}|\kappa_m(\xi)|}{|t - \xi|}.
\end{equation}
Hence,
$\left|2k\prod_{m=1}^{2k-1}\kappa_m(\xi) + g_e^{(2k)}(c)(t- \xi)\right| >0$. 
Then, $(t - \xi)^{2k-1}>0$ for all $t \in (\xi, \infty) \cap I$ and $(t - \xi)^{2k-1}<0$ for all $t \in (-\infty, \xi) \cap I$.
Since also, $g_e(\xi)=0$, this implies that $f_e$ is strictly increasing on precisely one of $(-\infty, \xi) \cap I,$ and $ (\xi, \infty) \cap I$ and strictly decreasing on the other. 

Observe that ${f_e(\eta)= f_e(\zeta)=0}$ and ${f_e \in C^{n+ 2k - 1}}$. Since also $\eta \neq \zeta$, we conclude that $\min\{\eta,\zeta \}< \xi < \max\{\eta,\zeta \}$, as $f_e$ has at most one zero on each side of $\xi$.
\end{proof}

Using the same method, we clarify why the defining condition of $\zeta =\Rmap$ is 
    \begin{equation}
    e_{2k}(\xi)\cdot(\phi(\eta) - \phi(\zeta)) =0,
\end{equation}
instead of say 
    \begin{equation}
    e_{j}(\xi)\cdot(\phi(\eta) - \phi(\zeta)) =0,
\end{equation}
for any $j \in \{1, ..., n\}$. 

\begin{proposition}[No third point in odd Frenet directions]
\label{prop:no_third_pt}
There exists a non-empty interval $U \subseteq I$ for which there do not exist three distinct points $\xi, \eta, \zeta \in U$ satisfying
\begin{equation}
    e_{2k-1}(\xi) \cdot (\phi(\eta) - \phi(\zeta))=0.
\end{equation}
\end{proposition}
\begin{proof}
Let $U \subseteq I$ be some interval which we later choose.
    Fix $\xi, \eta \in U$ with $\xi \neq \eta$. Define two auxiliary functions $f_o: I \to \mathbb{R}$ and $g_o:I \to \mathbb{R}$ via
    \begin{equation}
        f_o(t):= e_{2k-1}(\xi)\cdot (\phi(\eta)- \phi(t)),
        \qquad \text{and} \qquad
    g_o(t):= - f_o'(t) = e_{2k-1}(\xi)\cdot e_1(t).
\end{equation}
An analogous argument to that of Lemma \ref{xi_in_middle_n} reveals
\begin{align*}
    g_o(t)
    &=(t- \xi)^{2k-2}\left(\frac{(2k-1)\prod_{m=1}^{2k-2}\kappa_m(\xi) + g_o^{(2k-1)}(c)(t- \xi)}{(2k-1)!}\right),
\end{align*}
for some $c$ between $t$ and $\xi$.
Then, choosing $U$ to satisfy
\begin{equation}
    \left|(2k-1)\prod_{m=1}^{2k-2}\kappa_m(\xi) + g_o^{(2k-1)}(c)(t- \xi) \right| > 0,
\end{equation}
and using $(t - \xi)^{2k-2}>0$ for all $t \neq \xi$, it is evident that $|g_o'(t)| > 0$ for all $t \in U$ with $t \neq \xi$. Hence, $f_o$ is strictly monotonic on $U \setminus \{\xi\}$. Since also $f_o(\eta) =0$, we conclude by continuity that there does not exist a $\zeta \in I$ such that $\zeta \neq \eta$ and $f_e(\zeta) =0$.
\end{proof}

As in three dimensions, the first marginal property and the conjugate-symmetry of the Wigner transform will rely on the maps $\xi \mapsto \Rmap$ and $\eta \mapsto \Rmap$ respectively.
    For fixed $\xi \in I$, define
       \begin{equation}
\label{actual_interval_Ixi_n}
       I_{\xi, k} :=\{\eta \in I:  e_{2k}(\xi) \cdot (\phi(\eta) - \phi(\zeta))=0 \text{ for some } \zeta \neq \eta \}.
       \end{equation}
    Moreover, for fixed $\eta \in I$, define
\begin{equation}
\label{actual_interval_Ieta_n}
       \dot{I}_{\eta, k} :=\{\xi \in I:  e_{2k}(\xi) \cdot (\phi(\eta) - \phi(\zeta))=0 \text{ for some } \zeta \neq \eta \}.
       \end{equation}
To generalise Lemma \ref{lem:domains_Rmap} to higher dimensions, that is to see that $I_{\xi, k}$ and $\dot{I}_{\eta, k}$ are intervals, we require a generalisation of the injectivity condition \eqref{eqn_N'_dot_T_not_zero_condition}. We have already seen that when $k=1$, the appropriate anologue is condition \eqref{IC_ndim}. The following assumption is precisely \eqref{IC_ndim} when $k=1$ and is sufficient to impose when $k>1$:
\begin{equation}
\tag{TA3}
\label{eqn_ej'_dot_T_not_zero_condition}
    e_{2k}^{(2k -1)}(\xi)\cdot e_1(\eta) \neq 0 \qquad \text{for all } \xi, \eta \in I.
\end{equation}
We defer to Remark \ref{rem:TA3} for further discussion, emphasising that \eqref{eqn_ej'_dot_T_not_zero_condition} holds for $I$ sufficiently small.

\begin{lemma}[Domains of the $R$ map]
\label{lem:domains_Rmap_ndim}
Let $S$ satisfy \eqref{eqn_ej'_dot_T_not_zero_condition}.
Denote $I=[a, b]$ and define $f_e$ as in \eqref{aux_function_xi_in_middle_lemma_n}. For all $\xi \in I$,
\begin{equation}
\label{star_1_n}
\tag{$\star_1$}
    I_{\xi, k} = \begin{cases}
      [a', b] \qquad \text{if } f_e(a) \leq f_e(b) \\
      [a, b'] \qquad \text{if } f_e(a) > f_e(b),
  \end{cases}  
\end{equation}
where $a' \in [a, \xi]$ is the unique point such that $f_e(a') = f_e(b)$ and $b' \in [\xi, b]$ is the unique point such that $f_e(a) = f_e(b')$.

Now fix $\eta \in I$ and define $g\colon \dot{I}_{\eta, k} \to I$ via $g(t) = R_{t} \eta$. Then,
\begin{equation}
\label{star_2_n}
\tag{$\star_2$}
    \dot{I}_{\eta, k}= [a'', b''],
\end{equation}
where $a'', b'' \in I$ are defined such that $g(a'') = a$ and $g(b'') =b$. 
\end{lemma}  
\begin{proof}
The proof follows similarly to that of Lemma \ref{lem:domains_Rmap} by replacing in-between Lemma \ref{xi_in_middle} with its generalisation -- Lemma \ref{xi_in_middle_n}. The difference is in the derivation of $J(\xi, \eta) \neq 0$ for all $\xi, \eta \in I$.
Hence, by \eqref{explicit_J_formula_n}, it remains to show that $G(\xi) := e_{2k}'(\xi)\cdot (\phi(\eta)- \phi(\zeta)) \neq 0$.
To this end, observe that
    \begin{equation}
       G(\xi) = e_{2k}'(\xi)\cdot (\phi(\eta)- \phi(\zeta)) = \int_{\zeta}^{\eta} e_{2k}'(\xi)\cdot e_1(t)\rmd t.
    \end{equation}
    In the case where $k =1$, we are done by assumption \eqref{eqn_ej'_dot_T_not_zero_condition} and the continuity of $e_{2k}'(\xi)\cdot e_1(t)$ in $t$.
     Now assume that $k >1$.
Corollary \ref{corol_derivatives_frenet_vectors} stipulates that $e_{i}^{(j)}(\xi)\cdot e_1(\xi)=0$ if and only if $1 <i -j$.
Therefore, 
\begin{equation}
    e_{2k}^{(j)}(s)\cdot e_1(t) = (e_{2k}^{(j)}(s) - e_{2k}^{(j)}(t))\cdot e_1(t) \qquad \text{for all } j \leq 2k -2.
\end{equation}
Applying this identity and the fundamental theorem of calculus $(2k-2)$ times, with each integral taken with respect to one-dimensional Lebesgue measure on the indicated integral, we obtain
   \begin{equation}
        G(\xi)= \int_{\zeta}^{\eta} \int_{t}^{\xi} \int_{t}^{s_1}\cdots \int_{t}^{s_{2k-3}} e_{2k}^{(2k-1)}(s_{2k-2})\cdot e_1(t) \; \rmd s_{2k-2} \cdots\rmd s_2 \rmd s_1 \rmd t.
    \end{equation}
 Since $\phi \in C^{n+ 2k - 1}$, $\gamma(s) = e_{2k}^{(2k-1)}(s)\cdot e_1(t) \in C^{n- 2k}$, and hence is continuous for all $s \in I$ since $2k \leq n$.
 Hence, by assumption \eqref{eqn_ej'_dot_T_not_zero_condition}, we conclude $G(\xi)\neq 0$ for all $\xi \in I$.
\end{proof}
\begin{remark}[The map ${g\colon \dot{I}_{\eta, k} \to I}$ is bijective]
\label{rem:g_increas_bijective_ndim}
    We conclude as in Remark \ref{rem:g_increas_bijective} that ${g\colon\dot{I}_{\eta, k} \to I}$ bijective -- the injectivity of the map follows from $G(\xi) \neq 0$ for all $\xi \in I$ proven above. This is where we made use of the $n$-dimensional analogue of the injectivity condition (that is \eqref{eqn_ej'_dot_T_not_zero_condition} in place of \eqref{eqn_N'_dot_T_not_zero_condition}), indicating the reason for its name.
\end{remark}

\begin{remark}[Understanding technical assumption 3 \eqref{eqn_ej'_dot_T_not_zero_condition}: the $n$-dimensional injectivity condition]
\label{rem:TA3}
  When $k=1$, the Frenet formulae imply that 
$$e_{2}'(\xi)\cdot e_1(\xi) = (-\kappa_{1}(\xi) e_{1}(\xi) + \kappa_{2}(\xi) e_{3}(\xi))\cdot e_1(\xi) = -\kappa_{1}(\xi) \neq 0.$$
Hence, as we saw in three dimensions, it is reasonable to impose that $e_{2}'(\xi)\cdot e_1(\eta) \neq 0$ for all $\xi, \eta \in I$.
However, for $k >1$, 
$$e_{2k}'(\xi)\cdot e_1(\xi) = (-{\kappa_{2k-1}(\xi) e_{2k-1}(\xi) + \kappa_{2k}(\xi) e_{2k+1}(\xi)})\cdot e_1(\xi) = 0.$$ 
Therefore we cannot assume that $e_{2k}'(\xi)\cdot e_1(\xi) = 0$ for all $\xi, \eta \in I$. 
Corollary \ref{corol_derivatives_frenet_vectors}
demonstrated that $e_{2k}^{(m)}(\xi)\cdot e_1(\xi) = 0$ for all $m < 2k-1$ and $\xi \in I$. Hence, the first condition of this form which does not vanish when the points collide is the one imposed here; indeed $|e_{2k}^{(2k-1)}(\xi)\cdot e_1(\xi)| = \prod_{m=1}^{2k-1} |\kappa_m(\xi)| > 0$ by the well-curved hypothesis.
\end{remark}

As in three dimensions, we may show that $W_{2k}(g_1, g_2)$ is conjugate-symmetric and hence the Wigner Distribution, $W_{2k}g$, is real-valued.
\begin{proposition}[The Wigner Transform is Conjugate-Symmetric]
Let $S$ satisfy conditions \eqref{eqn_ej_dot_T_not_zero_condition}, \eqref{eqn_taylor_condn} and \eqref{eqn_ej'_dot_T_not_zero_condition}. Moreover, let $\xi, \eta, \zeta \in I$ be distinct points such that $\zeta = \Rmap$. Then, the Wigner transform is conjugate-symmetric, in that
    \begin{align*}
        W_{2k}(g_1, g_2)(\xi, v) = \overline{W_{2k}(g_2, g_1)(\xi, v)}, \qquad \text{for all } (\xi, v) \in D_{2k}.
    \end{align*}
\end{proposition}
\begin{proof}
Fix $\xi \in I$ and change variables via the map $\eta \to R_{\xi}\zeta$, noting $R_{\xi}(R_{\xi}(\eta)) = \eta$ and the map $R_{\xi}:I_{\xi, k} \to I_{\xi, k}$ is bijective. Moreover, for all $\eta\in I_{\xi, k}$ with $\eta \neq \xi$, the map is continuously differentiable, and
       \begin{equation}
       \label{second_Jacobian_n}
    \Delta(\xi, \eta) :=  \left| \frac{\partial \Rmap}{\partial \eta} \right| = \frac{|e_{2k}(\xi)\cdot e_1(\eta)|}{|e_{2k}(\xi)\cdot e_1(\zeta)|}.
    \end{equation}
Hence
\begin{equation}
\label{combine_Jac1_2_n}
    J(\xi, \eta) \Delta(\xi, \zeta) = J(\xi, \zeta),
\end{equation}
and the result follows as in three dimensions.
\end{proof}

\subsection{Wigner Representation}
We now follow the three-dimensional methods to prove Wigner representations for a curve in $\mathbb{R}^n$, beginning by introducing the auxiliary function
\begin{equation}
\label{aux_fun_F_k}
F_k: I \times \mathbb{R}^n \to \mathbb{R}, \qquad
    F_k(\xi, x) = \int_{I}g(\eta) \overline{g(R_{\xi}\eta)} e^{-2\pi ix\cdot (\phi(\eta) - \phi(R_{\xi}\eta))}J(\xi,\eta)\rmd \eta,
\end{equation}
so that $W_{2k}g$ is the restriction of $F_k$ to $D_{2k}$. Explicitly,
\begin{equation}
\label{defn_F_k}
    F_k(\xi, x) = W_{2k}g(\xi, P_{\langle e_{2k}(\xi) \rangle^{\perp}}x),
\end{equation}
where $P_{\langle e_{2k}(\xi) \rangle^{\perp}}$ is the orthogonal projection onto $\langle e_{2k}(\xi) \rangle^{\perp}$. 

\begin{proposition}[The First Marginal Property]
\label{first_marg_n}
Let $S$ satisfy conditions \eqref{eqn_ej_dot_T_not_zero_condition}, \eqref{eqn_taylor_condn} and \eqref{eqn_ej'_dot_T_not_zero_condition}. Then,
    \begin{equation}
    \int_{I} F_k(\xi,x) \rmd\xi = |E_Sg(x)|^2.
\end{equation}
\end{proposition}
\begin{proof}
The proof mirrors that of Proposition \ref{first_marg}, changing variables via the map ${\xi \mapsto \Rmap}$ by using \eqref{explicit_J_formula_n}, Lemma \ref{lem:domains_Rmap_ndim} and Remark \ref{rem:g_increas_bijective_ndim}.
\end{proof}
We define another restriction of the $X$-Ray transform, now as an integral over the line with direction $e_{2k}(\xi)$ for fixed $\xi \in I$:
\begin{equation}
    X_{2k}f(\xi, v) := Xf(e_{2k}(\xi), v) = \int_{\langle e_{2k}(\xi) \rangle} f(v + y) \rmd y,
\end{equation}

\begin{theorem}[Wigner representation] 
\label{phase_space_rep_n}
Let $S$ satisfy conditions \eqref{eqn_ej_dot_T_not_zero_condition}, \eqref{eqn_taylor_condn} and \eqref{eqn_ej'_dot_T_not_zero_condition}. Then,
\begin{equation}
    |E_Sg|^2 = X^{*}_{2k}W_{2k}g,
\end{equation}
where for functions ${h \colon D_{2k} \to \mathbb{C}}$,
\begin{equation}
    X_{2k}^* h(x) = \int_I h(\xi, P_{\langle e_{2k}(\xi) \rangle^{\perp}}x) \rmd \xi.
\end{equation}
\end{theorem}
As in three dimensions, $X_{2k}^*$ may be seen to be the adjoint of the restricted $X$-Ray transform.
\begin{proof}[Proof of Proposition \ref{phase_space_rep_n}]
The statement follows from the first marginal property, \eqref{defn_F_k}, and the definition of $X_{2k}^*$.
\end{proof}
Again, we have the following integral representation for all weight functions $w: \mathbb{R}^n \to \mathbb{C}$:
\begin{align}
\int_{\mathbb{R}^n}|E_Sg(x)|^2w(x) \rmd x
    &=  \int_I \int_{\langle e_{2k}(\xi) \rangle^{\perp}} W_{2k}g(\xi,v)X_{2k}w(\xi, v) \rmd v \rmd\xi.
\end{align}   

\subsection{Motivating applications via Wigner heuristics}
As in three dimensions, the $k^{\text{th}}$ curve-carried Wigner distribution also satisfies a tangent marginal property:
\begin{equation}
    \int_{\langle e_1(\xi) \rangle} W_{2k}g(\xi, v_t + v_b) \rmd v_t = |g(\xi)|^2 \qquad \text{for all } \xi \in I, v_b \in \langle e_1(\xi), e_{2k}(\xi) \rangle^{\perp},
    \end{equation}
whose full statement is given by Theorem \ref{Second_marginal_property_curve_Rn}. 
Now, for fixed $\xi \in I$, the marginal is \textit{independent} of $v_b \in \langle e_1(\xi), e_{2k}(\xi) \rangle^{\perp}$.

The same independence appears when considering $L^2$ norms of the $k^{\text{th}}$ curve-carried Wigner distribution, where we will see that
\begin{equation}
    \int_{\langle e_1(\xi) \rangle}|W_{2k}g(\xi, v_t + v_b)|^2 \rmd v_t = \int_{I} |g(\eta)|^2|g(\Rmap)|^2 \frac{J(\xi, \eta)^2}{\tilde{J}(\xi, \eta)} \rmd \eta, \qquad \text{for all }\xi \in I, v_b \in \langle e_1(\xi), e_{2k}(\xi) \rangle^{\perp},
\end{equation}
with $\tilde{J}(\xi, \eta)$ defined for curves in $n$-dimensions shortly. This implies that $W_{2k}g \notin L^2(D_k)$, as was expected from the work in three dimensions.
   
The choice of notation, emulating that of three dimensions, expresses that $v_t$ and $v_b$ have similar roles to their three dimensional counterparts.
The difference is in the range of values that $v_b$ may take: now $v_b \in \langle e_1(\xi), e_{2k}(\xi) \rangle^{\perp}$ which has dimension $n-2$, in place of $\langle B(\xi) \rangle$ with dimension $1$. 

To observe the impact of this difference, we follow the Wigner heuristics seen in Section \ref{subsec:wigner_heuristicsR3} and compare the proximity to \eqref{MT_conj_dim_m}. Indeed, if $W_{2k}g$ were non-negative, similar considerations would suggest
   \begin{equation}
       \label{wigner_heuristics_nd}
       \int_{\mathbb{R}^n} |E_S g|^2 w \leq \|g\|_2^2 \sup_{\xi} \int_{\langle e_1(\xi), e_{2k}(\xi)\rangle^\perp} \sup_{v_t} X_{2k} w(\xi, v_t + v_b)\rmd v_b.
   \end{equation}
   \begin{remark}[Proximity to \eqref{MT_conj_dim_m}]
Recall that
\begin{equation} \mathcal{R}_Tw(\xi, v_t) = \int_{\langle e_1(\xi), e_{2k}(\xi)\rangle^{\perp}}  X_{2k}w(\xi, v_t + v_b)     \rmd v_b,
\end{equation}
and observe that
\begin{equation} 
\label{MT_heuristic_compare}
\sup_{\xi} \sup_{v_t} \mathcal{R}_Tw(\xi, v_t) \leq \sup_{\xi} \int_{\langle e_1(\xi), e_{2k}(\xi)\rangle^{\perp}} \sup_{v_t} X_{2k}w(\xi, v_t + v_b)     \rmd v_b,
\end{equation}
where the left-hand side appears in \eqref{MT_conj_dim_m} and the right-hand side in \eqref{wigner_heuristics_nd}. 
\end{remark}
As in three dimensions, we cannot impose that $W_{2k}$ is non-negative, however, the heuristic provides a measure of proximity to the conjectured bounds and anticipates the structure of the bounds we expect to prove; the Mizohata--Takeuchi-type inequalities proven in this work employ a different norm in $v_t$ on the right-hand side of \eqref{MT_heuristic_compare}.
   
\medskip
\textbf{The third map.}
Since $v_t$ has the same role as in three dimensions, the third map in higher dimensions is $\gamma_T: I \to \tilde{I}$ defined for fixed $\xi \in I$ via
\begin{equation}
\label{third_map_defn_n_dim}
    \gamma_T(\eta) := e_1(\xi)\cdot (\phi(\eta)- \phi(\Rmap)),
\end{equation}
for all $\eta \in I$.
This is precisely the third map seen in three dimensions, and we make analogous remarks before endeavouring to apply the Wigner representation.
First, we impose the following assumption on the curves considered:
\begin{equation}
\label{eqn_tangents_condition_n}
\tag{TCn}
    e_1(\xi)\cdot e_1(\eta) \geq \frac{1}{2} \quad \text{for all $\xi, \eta \in I$},
\end{equation} 
the immediate generalisation of the tangent condition \eqref{eqn_tangents_condition}.

\begin{lemma}
\label{lemma_Jtilde_identity}
Fix $\xi \in I$.
Then $\gamma_T$ is continuously differentiable on $I_{\xi} \setminus \{\xi\}$ and
    \begin{equation}
        \tilde{J}(\xi, \eta):= |\gamma_T'(\eta)| = |e_1(\xi) \cdot \left(e_1(\eta)  + \Delta(\xi, \eta) e_1(R_{\xi}\eta) \right)|
        =\frac{|\langle e_{2k}(\xi)  \wedge e_1(\xi) , e_1(R_{\xi}\eta)  \wedge e_1(\eta) \rangle|}{|e_{2k}(\xi)\cdot e_1(\Rmap)|}.
    \end{equation}
\end{lemma}
In the above, $\langle \cdot, \cdot \rangle$ denotes the standard inner product on exterior powers -- see for example, Section 1.6.5 of \cite{darling_book}. Explicitly, for $v_1, v_2, w_1, w_2 \in \mathbb{R}^n$,
\begin{equation}
    \label{inner_prod_exterior_algebra}
    \langle v_1 \wedge v_2,  w_1 \wedge w_2 \rangle  = 
    \det(\langle v_i, w_j \rangle),
\end{equation}
where the inner product on the right hand side is the standard inner product on $\mathbb{R}^n$.
\begin{proof}
The argument follows similarly to the proof of Lemma \ref{lemma_Jtilde_identity_3D}, using that $\Delta(\xi, \eta) = \frac{- \partial \Rmap}{\partial \eta}$ follows using a Taylor expansion argument similar to those seen in this section, and the in-between lemma. Then, for fixed $\xi \in I$,
    \begin{equation}
     \gamma'_T(\eta) = \frac{e_1(\xi) \cdot \left((e_{2k}(\xi) \cdot e_1(R_{\xi}\eta) )e_1(\eta)  - (e_{2k}(\xi) \cdot e_1(\eta) ) e_1(R_{\xi}\eta) \right)}{e_{2k}(\xi) \cdot e_1(R_{\xi}\eta) }.
\end{equation}
Using \eqref{inner_prod_exterior_algebra}, the numerator becomes
\begin{align*}
e_1(\xi) \cdot &((e_{2k}(\xi) \cdot e_1(R_{\xi}\eta) )e_1(\eta)  - (e_{2k}(\xi) \cdot e_1(\eta) ) e_1(R_{\xi}\eta) ) \\
&=
    (e_{2k}(\xi) \cdot e_1(R_{\xi}\eta) )(e_1(\xi) \cdot e_1(\eta) ) - (e_{2k}(\xi) \cdot e_1(\eta) ) (e_1(\xi) \cdot e_1(R_{\xi}\eta) ) \\
    &= \langle e_{2k}(\xi)  \wedge e_1(\xi) , e_1(R_{\xi}\eta)  \wedge e_1(\eta) \rangle.
\end{align*}
\end{proof}
\begin{remark}[Bijectivity of the third map]
Let $S$ satisfy conditions \eqref{eqn_taylor_condn} and \eqref{eqn_tangents_condition_n}.
\label{bijective_third_map_n}
By the same argument as in Remark \ref{bijective_third_map}, now using \eqref{eqn_tangents_condition_n}, the third map is  bijective.
\end{remark}
We have an immediate generalisation of Lemma \ref{gammaT_sim_a-b} to curves in $\mathbb{R}^n$.
\begin{lemma}
\label{lem:gammaT_sim_a-b_ndim}
Let $S$ satisfy condition \eqref{eqn_tangents_condition}.  Then for all $\xi, \eta, \zeta \in I$, we have
    \begin{equation}
    \label{gammaT_sim_a-b_ndim}
       \frac{1}{2}|\eta - \zeta| \leq |e_1(\xi) \cdot (\phi(\eta)- \phi(\zeta)) | \leq |\eta - \zeta|,
    \end{equation}
    and
    \begin{equation}
    \label{phia-phib_sim_a-b_ndim}
       \frac{1}{\sqrt{2}}|\eta - \zeta| \leq |\phi(\eta) - \phi(\zeta)| \leq |\eta - \zeta|.
    \end{equation}
\end{lemma}

\subsection{Application 1: \texorpdfstring{$L^2$}{L2} Weighted Inequalities}
Now we turn to prove the $L^2$ Sobolev--Stein and Sobolev--Mizohata--Takeuchi inequalities for curves in $\mathbb{R}^n$.

\begin{theorem}[$L^2$ Sobolev--Stein Inequality for Curves]
\label{corol2a_l2sobolev_n}
Let $S$ satisfy conditions \eqref{eqn_ej_dot_T_not_zero_condition}, \eqref{eqn_taylor_condn}, \eqref{eqn_ej'_dot_T_not_zero_condition}, and \eqref{eqn_tangents_condition_n}. Moreover, let $s < \frac{1}{2}$. Then, there exists a finite constant $c_s$ depending only on $s$ such that
    \begin{equation}
    \int_{\mathbb{R}^n} |E_Sg(x)|^2w(x) \rmd x 
    \leq c_s  C_k(\phi) \int_I \mathcal{B}_{S, 2s}(|g|^2, |g|^2)(\xi)^{\frac{1}{2}} \left( \int_{\langle B(\xi) \rangle} \left\| X_{2k}w(\xi, \cdot + v_b)) \right\|_{\dot{H}_{v_t}^s} \rmd v_b \right)  \rmd\xi,
\end{equation}
where
\begin{equation}
    \mathcal{B}_{S, s}(g_1, g_2)(\xi) = \int_{I}\frac{g_1(\eta)\overline{g_2(R_{\xi}\eta)} J(\xi,\eta)}{|\phi(\eta)- \phi(\Rmap)|^{s}}\rmd \eta,
\end{equation}
and
\begin{equation}
    C_k(\phi) :=\sup_{\xi, \eta \in I}
    \left(\frac{e_{2k}'(\xi)\cdot (\phi(\eta) - \phi(R_{\xi}\eta))}{\langle e_{2k}(\xi)  \wedge e_1(\xi) , e_1(R_{\xi}\eta)  \wedge e_1(\eta)\rangle} \right)^{\frac{1}{2}}.
\end{equation}
\end{theorem}
\begin{proof}[Proof of Theorem \ref{corol2a_l2sobolev_n}]
We begin with the Wigner integral representation and decompose $v = v_t + v_b$ where $v_t$ and $v_b$ are the orthogonal projections of $v$ onto $\langle e_1(\xi) \rangle$ and $\langle e_1(\xi), e_{2k}(\xi) \rangle^{\perp}$ respectively:
\begin{align*}
    \int_{\mathbb{R}^n} |E_Sg(x)|^2w(x) \rmd x 
    &=  \int_I \int_{\langle e_1(\xi), e_{2k}(\xi) \rangle^{\perp}}\int_{\langle e_{1}(\xi) \rangle} W_{2k}g(\xi, v_t + v_b)X_{2k}w(\xi, v_t + v_b) \rmd v_t \rmd v_b \rmd\xi \\
    &\leq \int_I \int_{\langle e_1(\xi), e_{2k}(\xi) \rangle^{\perp}}
    \|W_{2k}g(\xi,\cdot + v_b)\|_{\dot{H}^{-s}_{v_t}} \|X_{2k}w(\xi, \cdot + v_b)\|_{\dot{H}^{s}_{v_t}}
    \rmd v_b
    \rmd \xi. \qquad (\star)
\end{align*} 
Fix $\xi \in I$ and $v_b \in \langle e_1(\xi), e_{2k}(\xi) \rangle^{\perp}$. Arguing as in three dimensions, changing variables via the analogous map $\gamma_T: I \to \tilde{I}$ defined by \eqref{third_map_defn_n_dim} and now making use of Remark \ref{bijective_third_map_n}, we have
\begin{align*}
    \|W_{2k}g(\xi,\cdot + v_b)\|_{\dot{H}^{-s}(\langle e_1(\xi) \rangle)}^2 
     &= (2 \pi)^{-2s}
    \int_{I}\frac{|g(\eta)|^2  |g(R_{\xi}\eta)|^2 J(\xi,\eta)^2}{\tilde{J}(\xi, \eta)|e_1(\xi) \cdot (\phi(\eta) - \phi(\Rmap))|^{2s}}\rmd \eta.
\end{align*}
Identity \eqref{explicit_J_formula_n} and the second identity of Proposition \ref{lemma_Jtilde_identity}, together imply that
    \begin{equation}
        \frac{J(\xi, \eta)}{\tilde{J}(\xi, \eta)}
        = \frac{e_{2k}'(\xi)\cdot (\phi(\eta) - \phi(R_{\xi}\eta))}{\langle e_{2k}(\xi)  \wedge e_1(\xi) , e_1(R_{\xi}\eta)  \wedge e_1(\eta)\rangle}
        \leq C_k^2(\phi).
    \end{equation} 
Using also Lemma \ref{lem:gammaT_sim_a-b_ndim} yields
\begin{align}
    \| W_{2k}g(\xi,\cdot + v_b)\|_{\dot{H}^{s}_{v_t}} \leq \pi^{-s} C_k(\phi) \mathcal{B}_{S, 2s}(|g|^2, |g|^2)(\xi)^{1/2}.
\end{align}
Returning to $(\star)$, we conclude by using the independence of the first homogeneous Sobolev norm on $v_b \in \langle e_1(\xi), e_{2k}(\xi) \rangle^\perp$:
\begin{align}
    \int_{\mathbb{R}^n} &|E_Sg(x)|^2w(x) \rmd x
    &= \pi^{-s} C_k(\phi) \int_I \mathcal{B}_{S, 2s}(|g|^2, |g|^2)(\xi)^{\frac{1}{2}} \left( \int_{\langle e_1(\xi), e_{2k}(\xi) \rangle^{\perp}} \left\| X_{2k}w(\xi, \cdot  + v_b) \right\|_{\dot{H}_{v_t}^s} \rmd v_b \right)  \rmd\xi.
\end{align}
\end{proof}

\begin{theorem}[$L^2$ Sobolev--Mizohata--Takeuchi Inequality for Curves]
\label{corol2b_l2sobolev_n}Let $S$ satisfy conditions \eqref{eqn_ej_dot_T_not_zero_condition}, \eqref{eqn_taylor_condn}, \eqref{eqn_ej'_dot_T_not_zero_condition}, and \eqref{eqn_tangents_condition_n}. Moreover, let $s < \frac{1}{2}$. Then there exists a finite implicit constant depending only on $s$ and the diameter of $S$ such that
    \begin{equation}
    \int_{\mathbb{R}^n} |E_Sg(x)|^2w(x) \rmd x 
    \lesssim_s C_k(\phi) \sup_{\xi \in \supp^*(g)} \left( \int_{\langle e_1(\xi), e_{2k}(\xi) \rangle^\perp} \left\| X_{2k}w(\xi, \cdot + v_b)) \right\|_{\dot{H}_{v_t}^s} \rmd v_b \right) \|g\|_{L^2(I)}^2.
\end{equation}
\end{theorem}
\begin{proof}[Proof of Theorem \ref{corol2b_l2sobolev_n}]
The result follows from Theorem \ref{corol2a_l2sobolev_n} using an analogous argument to that of Theorem \ref{corol2b_l2sobolev}.
\end{proof}

\medskip
\textbf{Understanding $C_k(\phi)$.}
We now seek to control $C_k(\phi)$. The goal is to remove the singularity, and indicate a simple additional technical condition under which $C_k(\phi)$ can easily be seen to be finite. 

\begin{lemma}[Controlling $C_k(\phi)$]
\label{lem:Ck_control}
Let $S$ satisfy \eqref{eqn_ej_dot_T_not_zero_condition}, \eqref{eqn_taylor_condn}, \eqref{eqn_ej'_dot_T_not_zero_condition} and \eqref{eqn_tangents_condition_n}. Then,
    \begin{equation}
        C_k(\phi) = \sup_{\xi, \eta \in I}\frac{J(\xi, \eta)}{\tilde{J}(\xi, \eta)}\leq 2\sup_{\xi, s, t \in I} \frac{|e_{2k}'(\xi) \cdot e_1^{(2k-2)}(s)|}{|e_{2k}(\xi) \cdot e_1^{(2k-1)}(t)|}.
    \end{equation}
\end{lemma}
\begin{proof}
    First, we control $J(\xi, \eta)$ using \eqref{explicit_J_formula_n}:
    \begin{equation}
    \label{eqn:control_ck}
        J(\xi, \eta) = \frac{|e_{2k}'(\xi)\cdot (\phi(\eta) - \phi(\zeta))|}{|e_{2k}(\xi)\cdot e_1(\zeta)|} := \frac{|f(\eta) - f(\zeta)|}{|g(\zeta)|} = \frac{|(f(\eta) - f(\xi) )- (f(\zeta)-f(\xi))|}{|g(\zeta)|},
    \end{equation}
where we have defined
\begin{equation}
    f(t):= e_{2k}'(\xi)\cdot \phi(t) \qquad g(t):= e_{2k}(\xi)\cdot e_1(t).
\end{equation}
Observe that Corollary \ref{corol_derivatives_frenet_vectors} stipulates that
\begin{equation}
    f'(\xi) = f''(\xi) = \cdots = f^{(2k-2)}(\xi) =0, \qquad g(\xi) = g'(\xi) = g''(\xi) = \cdots = g^{(2k-2)}(\xi) =0.
\end{equation}
Hence, Taylor expansions of $f$ and $g$ up to order $m-1$ yield
\begin{equation}
\label{in_lem_C_phi:f_g_deriv}
    f(t) - f(\xi) = \frac{f^{(2k-1)}(a_t) (t-\xi)^{2k-1}}{(2k-1)!}, \qquad g(t) = g(t) - g(\xi) = \frac{g^{(2k-1)}(b_t) (t-\xi)^{2k-1}}{(2k-1)!}.
\end{equation}
Returning to \eqref{eqn:control_ck} and using the triangle inequality, we see that
\begin{align*}
    J(\xi, \eta) = 
    \left| 
    \frac{f^{(2k-1)}(a_{\eta}) (\eta-\xi)^{2k-1} - f^{(2k-1)}(a_{\zeta}) (\zeta-\xi)^{2k-1}}{g^{(2k-1)}(b_{\zeta}) (\zeta-\xi)^{2k-1}}
    \right|
    \leq 
    \frac{|f^{(2k-1)}(a_{\eta})|}{|g^{(2k-1)}(b_{\zeta})|} 
    \left|\frac{\eta-\xi}{\zeta-\xi}\right|^{2k-1}
    + \frac{|f^{(2k-1)}(a_{\zeta})|}{|g^{(2k-1)}(b_{\zeta})|}.
\end{align*}
Now observe that by \eqref{second_Jacobian_n} and \eqref{in_lem_C_phi:f_g_deriv}, we have
\begin{equation}
    \Delta(\xi, \eta) = \frac{|g(\eta)|}{|g(\zeta)|} = \frac{|g^{(2k-1)}(b_{\eta})|}{|g^{(2k-1)}(b_{\zeta})|} 
    \left|\frac{\eta-\xi}{\zeta-\xi}\right|^{2k-1}.
\end{equation}
Therefore,
\begin{align*}
    J(\xi, \eta)
    \leq 
    \frac{|f^{(2k-1)}(a_{\eta})|}{|g^{(2k-1)}(b_{\eta})|} 
    \Delta(\xi, \eta)
    + \frac{|f^{(2k-1)}(a_{\zeta})|}{|g^{(2k-1)}(b_{\zeta})|} = G(\xi, a_{\eta}, b_{\eta})\Delta(\xi, \eta) +  G(\xi, a_{\zeta}, b_{\zeta}) \leq (1 + \Delta(\xi, \eta))\sup_{a, b \in I}G(\xi, a, b),
\end{align*}
where we have defined $G(\xi, a, b) := \frac{|f^{(2k-1)}(a)|}{|g^{(2k-1)}(b)|}$. 

Second, we seek to control $ \tilde{J}(\xi, \eta)$ by using Lemma \ref{lemma_Jtilde_identity} and the tangent condition \eqref{eqn_tangents_condition_n}:
\begin{equation}
    \tilde{J}(\xi, \eta) = |e_1(\xi)\cdot e_1(\eta) + \Delta(\xi, \eta) e_1(\xi)\cdot e_1(\zeta)| \geq \frac{1}{2}(1+ \Delta(\xi, \eta)).
\end{equation}
Hence,
\begin{equation}
    \frac{J(\xi, \eta)}{\tilde{J}(\xi, \eta)} \leq 2 \frac{(1 + \Delta(\xi, \eta))\sup_{a, b \in I}G(\xi, a, b)}{1+ \Delta(\xi, \eta)} = 2\sup_{a, b \in I}G(\xi, a, b),
\end{equation}
and we conclude by taking suprema.
\end{proof}
\begin{remark}[Finiteness of $C_k(\phi)$]
    We observe that there is no longer a singularity at $\xi$ since 
\begin{equation}
    f^{2k-1}(\xi) = - \prod_{m=1}^{2k-1}\kappa_m(\xi), \qquad g^{2k-1}(\xi) = \prod_{m=1}^{2k-1}\kappa_m(\xi).
\end{equation}
In order to avoid an additional singularity at some other point $\eta \in I$, we impose that
\begin{equation}
\label{final_technical_condition}
\tag{TA4}
    e_{2k}(\xi) \cdot e_1^{(2k-1)}(\eta) \neq 0 \qquad \text{for all $\xi, \eta \in I$.}
\end{equation}
Under condition \eqref{final_technical_condition}, $C_k(\phi)$ is finite.
\end{remark}
\begin{remark}[Interpreting condition \eqref{final_technical_condition}]
    We comment that this condition is similar to condition \eqref{eqn_ej'_dot_T_not_zero_condition} in the sense that it may be interpreted via an understanding of the derivatives of the Frenet vectors -- see Remark \ref{rem:TA3}.
\end{remark}

\begin{remark}[When $k=1$]
In the case of $k=1$, we may assume 
condition \eqref{NC_ndim}, and consequently have the further upper bound which is more intuitive to interpret:
     \begin{equation}
     \label{simple_bound_C1}
        C_1(\phi) \leq 4\sup_{\xi, \eta \in I} \frac{|\kappa_1(\xi)| + |\kappa_2(\xi)|}{\kappa_1(\eta)},
    \end{equation}
    To see why \eqref{simple_bound_C1} holds, observe that Lemma \eqref{lem:Ck_control}, the Frenet formulae, and assumption \eqref{NC_ndim} imply
      \begin{align*}
      C_1(\phi) \leq 2\sup_{\xi, s, t \in I} \frac{|e_{2}'(\xi) \cdot e_1(s)|}{|e_{2}(\xi) \cdot e_1'(t)|} =  2\sup_{\xi, s, t \in I} \frac{|(-\kappa_1(\xi)e_1(\xi) + \kappa_2(\xi) e_3(\xi)) \cdot e_1(s)|}{\kappa_1(t)|e_{2}(\xi) \cdot e_2(t)|} \leq 4\sup_{\xi, t \in I} \frac{|\kappa_1(\xi)| + |\kappa_2(\xi)|}{|\kappa_1(t)|}.
    \end{align*}
This recovers the same expression as in three dimensions.
Moreover, we remark that
\begin{equation}
    \tilde{C}(\phi):= \sup_{\xi, t \in I}\frac{|\kappa_1(\xi)| + |\kappa_2(\xi)|}{|\kappa_1(t)|}
\end{equation}
can easily be seen to be scale invariant by the same argument as in Remark \ref{rem:scale_invar_const_3d} three dimensions.
\end{remark}

 \section{Application 2: Orthonormal Weighted \texorpdfstring{$L^2$}{L2} Inequalities}
 \label{sec:orthonormal}
\subsection{Background}
Thus far, we have considered weighted estimates for $|E_Sg|^2$. Now we explore what happens if we replace $|E_Sg|^2$ with
\begin{equation}
    \sum_j \lambda_j |E_Sg_j|^2,
\end{equation}
for an orthonormal sequence $(g_j)_{j \in \mathbb{N}}$ with each $g_j \in L^2(I)$.
Such a question was posed in the setting of Mizohata--Takeuchi inequalities for hypersurfaces by Bennett, Bez, Gutiérrez, Nakamura and Oliveira in \cite{BBGNO}, where they made use of the Wigner representation of \cite{BGNO} to provide one such answer.

The question is motivated by the study of Fourier restriction, where, for smooth, compact, hypersurfaces, Frank and Sabin \cite{Frank_Sabin} proved an orthonormal variant of the classical Stein-Tomas theorem. We refer there and to \cite{BBGNO} for further background in the setting of Strichartz estimates and in the study of fermions. The orthonormal result of Frank and Sabin is stronger than the Stein-Tomas estimate in the sense that an analogous orthonormal result obtained via Stein-Tomas and the triangle inequality is strictly weaker than their result.
A similar phenomenon has been considered in the setting of well-curved curves. In this setting, Fourier restriction has been completely resolved by Drury in \cite{Drury1985_curves}. Moreover, orthonormal estimates for restriction to curves have been studied by \cite{Ko_thesis} and recently by \cite{CDH_ONFR_curves}.

Returning to the setting of orthonormal weighted inequalities, we ask what might a sensible orthonormal variant of \eqref{MT_conj_dim_m} be for curves. Motivated by the wave-packet considerations of \cite{BBGNO} in the setting of hypersurfaces, a first suggestion might ask for what $1\leq p < \infty$ and what sequences of functions $(g_j)$, with $g_j:I \to \mathbb{C}$ for each $j \in \mathbb{N}$, does there exist a finite implicit constant for which
\begin{equation}
\label{ONMT_conj_guess1}
    \int_{B_R} \sum_j |E_Sg(x)|^2w(x) \rmd x \lesssim R^{\alpha} \|(\lambda_j)\|_{\ell^{\beta}}  \|{\mathcal{R}_Tw}\|_{L^p(\{(\xi, v): \xi \in E, v \in \langle \omega \rangle^{\perp})}.
\end{equation}
This inequality has the following equivalent form by duality:
\begin{equation}
\label{ONMT_conj_guess1_dual}
    \sum_j\left( \int_{B_R} |E_Sg(x)|^2w(x) \rmd x \right)^p \lesssim R^{\alpha p}  \|{\mathcal{R}_Tw}\|_{L^p(\{(\xi, v): \xi \in E, v \in \langle \omega \rangle^{\perp})}^p.
\end{equation}
In both inequalities, $E \subseteq I$ containing $K = \cup_j \supp(g_j)$.
When $p=2$, Bennett, Bez, Gutiérrez, Nakamura, and Oliveira provided two answers to this question for hypersurfaces via two different approaches: the first is a Schatten approach using the techniques of Frank and Sabin \cite{Frank_Sabin} and the second is a Wigner approach, making use of the Wigner representation of \cite{BGNO}.

When $p= \infty$, we expect this to be false in its global form as a result of \cite{counterexample_MT} and \cite{MT_local_counterexample}. However, we are motivated by the Wigner heuristics to adapt the question we pose, especially since we will aim to use the Wigner representation to obtain results in the $p=2$ case. Driven by such considerations, we instead ask for $1 \leq p < \infty$, for what sequences of functions $(g_j)$, with $g_j:I \to \mathbb{C}$ for each $j \in \mathbb{N}$, does there exist a finite constant $C_p$ such that
\begin{equation}
\label{orthonormal_lp_question_curve_rn}
\tag{ONMT}
 \sum_j \left(\int_{\mathbb{R}^n}|E_Sg_j(x)|^2 w(x)\right)^p \leq C_p \left\| \int_{\langle e_1(\xi), e_{2k}(\xi)\rangle^{\perp}}
        \|X_{2k}w(\xi, v_t + v_b)\|_{L^p_{v_t}(\langle e_1(\xi) \rangle)} \rmd v_b \right\|_{L^p_{\xi}(E)}^p.
\end{equation}

When $p= \infty$, this is precisely what the Wigner heuristic suggested in the single input setting. Again, while we do not expect it to be true when $p= \infty$, this inequality takes the form of something more structurally similar to what we proved for the Sobolev-variants motivated by the heuristics.

When $p=1$, \eqref{orthonormal_lp_question_curve_rn} becomes
\begin{equation}
\sum_j \int_{\mathbb{R}^n}|E_Sg_j(x)|^2 w(x) \leq C_1 
        \|X_{2k}w\|_{L^1(\{(\xi, v_t, v_b): \xi \in E, v_t \in \langle e_1(\xi) \rangle, v_b \in \langle e_1(\xi), e_{2k}(\xi) \rangle^{\perp} \})}.
\end{equation}
In this setting, we do not see the distinction between the two possible variants of the problem.
Moreover, by an analogous argument as in \cite{BBGNO}, there exists a finite $C_1$ such that this inequality holds when $(g_j)$ is an orthonormal sequence in $L^2(I)$ and $E=K$. Indeed, notice that
\begin{equation}
 \sum_j \int_{\mathbb{R}^n}|E_Sg_j (x)|^2 w(x) = 
  \int_{\mathbb{R}^n} \left(\sum_j|E_Sg_j (x)|^2 \right) w(x) \rmd x,
\end{equation}
and
\begin{equation}
    \|X_{2k}w\|_{L^1(\{(\xi, v_t, v_b): \xi \in K, v_t \in \langle e_1(\xi) \rangle, v_b \in \langle e_1(\xi), e_{2k}(\xi) \rangle^{\perp} \})}= |K|\|w\|_1.
\end{equation}
Therefore \eqref{orthonormal_lp_question_curve_rn} when $p=1$ is equivalent to
\begin{equation}
\label{ONMT1_niceform}
\int_{\mathbb{R}^n} \left(\sum_j|E_Sg_j (x)|^2 \right) w(x) \rmd x  \lesssim |K| \|w\|_1.
\end{equation}
To see that this holds true, observe that by an application of Bessel's inequality,
\begin{equation}
    \sum_j|E_Sg_j (x)|^2 = \sum_j\left(\int_{K} g_j(\xi) e^{2\pi i x \cdot \phi(\xi)} \rmd \xi\right)^2 = \sum_j| \langle g_j, e^{2\pi i x \cdot \phi(\cdot)} \1_K \rangle|^2 \leq \int_{K}|e^{2\pi i x \cdot \phi(\xi)}|^2 \rmd \xi = |K|.
\end{equation}

The case of $p=2$ is a natural setting in which to work, and we provide one such answer to this question motivated by the work of \cite{BBGNO} using a different notion of orthonormality and the Wigner representation. This forms the second application of the Wigner representation and is precisely upcoming Theorem \ref{thm_orthonormal_MT_Rn_curve}.
 
    When $p=1$, the norm of the tomographic transform is easier to interpret than that of our Sobolev--Mizohata--Takeuchi result in the single input setting, or our upcoming result when $p=2$ in the orthonormal setting. 
Since it is easier to interpret, we make some observations relevant in contextualising the single input setting. In particular, observe that
\begin{equation}
    \|w\|_{L^1(B_R)} \lesssim R \|\mathcal{R}w\|_{L^\infty}.
\end{equation}
In the single-input setting, this means \eqref{ONMT1_niceform} implies \eqref{MT_conj_dim_m} for curves with $\alpha =1$. 
This can be compared with the trivial bound for the Mizohata--Takeuchi conjecture implied by Agmon--Hörmander which shows that \eqref{MT_conj_dim_m} for curves holds with $\alpha = n-1$. It can also be compared with \cite{CLPY_MT_curves} and \cite{Cairo_Zhang}.
Moreover, 
\begin{equation}
    \|w\|_{L^1(B_R)} \lesssim R^{n-1} \|Xw\|_{L^\infty},
\end{equation}
implying $\alpha = n-1$ in a variant of \eqref{MT_conj_dim_m} with $\|Xw\|_{L^{\infty}}$ replacing $\|\mathcal{R}w\|_{L^{\infty}}$. We refer to \cite{CLPY_MT_curves} and \cite{Cairo_Zhang} for such variants of this problem.
\subsection{An Orthonormal Mizohata--Takeuchi inequality for \texorpdfstring{$p=2$}{p equals 2}}
\label{subsec:ONMT2}
\begin{theorem}
\label{thm_orthonormal_MT_Rn_curve}
Let $S$ satisfy conditions \eqref{eqn_ej_dot_T_not_zero_condition},  \eqref{eqn_taylor_condn}, 
\eqref{eqn_ej'_dot_T_not_zero_condition} and \eqref{eqn_tangents_condition_n}. Moreover, let $(g_j)$ be a sequence of functions on $I$ for which 
\begin{equation}
    L(i,j) := \int_I \int_{\langle e_1(\xi) \rangle} W_{2k}g_i(\xi, v) \overline{W_{2k}g_j(\xi, v)} \: \rmd \mathcal{H}^1(v) \:\rmd \xi,
\end{equation}
is the (real-valued) kernel of an $l^2$-bounded operator.\footnote{Here, $\rmd \mathcal{H}^1$ denotes one dimensional Hausdorff measure on $\langle e_1(\xi) \rangle$.}
    Then,
    \begin{equation}
        \sum_j \left(
        \int_{\mathbb{R}^n}|E_Sg_j(x)|^2w(x) \rmd x
        \right)^2 \lesssim 
    \|X_{2k}w\|_{L^2_{\xi}(K^*)L^1(\langle e_1(\xi), e_{2k}(\xi)\rangle^{\perp})L^2(\langle e_1(\xi)\rangle)}^2
    \end{equation}
    for all signed weight functions $w$. Here
    \begin{equation}
        K^*:= \bigcup_j \supp^* (g_j).
    \end{equation}
\end{theorem}
Explicitly,
\begin{equation}
    \|X_{2k}w\|_{L^2_{\xi}(K^*)L^1(\langle e_1(\xi), e_{2k}(\xi)\rangle^{\perp})L^2(\langle e_1(\xi)\rangle)}^2 := \int_{K^*} \left( 
    \int_{\langle e_1(\xi), e_{2k}(\xi) \rangle^{\perp}}
    \left(
    \int_{\langle e_1(\xi) \rangle} |X_{2k}w(\xi, v_t + v_b)|^2 \rmd v_t
    \right)^{1/2}
    \rmd v_b
    \right)^2 \rmd \xi 
\end{equation}

To interpret the almost-orthogonality condition that $L(i, j)$ is the kernel of an $l^2$-bounded operator, we turn to another key property of the Wigner transform called \textit{Moyal's identity}. For the classical Wigner transform, this states that
\begin{equation}
    \langle W(f_1, f_2), W(g_1, g_2) \rangle = \langle f_1, g_1 \rangle \overline{\langle f_2, g_2 \rangle}.
\end{equation}
We wish to prove an analogous Moyal identity for the curve-carried Wigner transform both to interpret the assumption on $L(i,j)$ and for completeness in understanding the Wigner transform we have defined. However, as previously established, $W_{2k}(g_1, g_2) \notin L^2(D_{2k})$, and due to previous considerations regarding the independence on the `binormal' direction (more explicitly on $\langle e_1(\xi), e_{2k}(\xi) \rangle^{\perp}$), we instead seek a Moyal-type identity for \begin{equation}
    \int_I \int_{\langle e_1(\xi) \rangle} W_{2k}(f_1, f_2)(\xi, v) \overline{W_{2k}(g_1, g_2)(\xi, v)} \: \rmd \mathcal{H}^1(v) \:\rmd \xi.
\end{equation}
\begin{proposition}[Moyal's Identity for the Curve-Carried Wigner Transform]
\label{prop_Moyal_curve_Rn}
Let $S$ satisfy conditions \eqref{eqn_ej_dot_T_not_zero_condition}, \eqref{eqn_taylor_condn}, \eqref{eqn_ej'_dot_T_not_zero_condition}, and \eqref{eqn_tangents_condition_n}. 
Then, for all $\xi \in I$ and $v_b \in \langle e_1(\xi), e_{2k}(\xi) \rangle^\perp$,
\begin{align}
\label{Moyal_indep}
    \int_{\langle e_1(\xi) \rangle} W_{2k}(f_1, f_2)(\xi, v_t + v_b) &\overline{W_{2k}(g_1, g_2)(\xi, v_t + v_b)} \: \rmd v_t \\
    &= \int_{I}f_1(\eta) \overline{f_2(R_{\xi}\eta)} \; \overline{g_1(\eta)}   g_2(R_{\xi}\eta)
    \frac{ J(\xi,\eta)^2}{\tilde{J}(\xi, \eta)} \rmd \eta.
\end{align}
Consequently,
    \begin{align}
        \int_I \int_{\langle e_1(\xi) \rangle} W_{2k}(f_1, f_2)(\xi, v) \overline{W_{2k}(g_1, g_2)(\xi, v)} \: \rmd \mathcal{H}^1(v) \:\rmd \xi=
        \int_I f_1(\eta)\overline{g_1(\eta)} \left(\int_I \overline{f_2(\zeta)} g_2(\zeta) M(\eta, \zeta) \rmd \zeta \right) \rmd \eta,
    \end{align}
    where 
    \begin{equation}
        M(\eta, \zeta) = \frac{|e_{2k}'(\xi)\cdot (\phi(\eta) - \phi(\zeta))|}{|\langle e_{2k}(\xi)  \wedge e_1(\xi) , e_1(\zeta)  \wedge e_1(\eta)\rangle|},
    \end{equation}
and we have written $\xi= \xi(\eta, \zeta)$.
\end{proposition}
\noindent

We highlight the lack of dependence on $v_b \in \langle e_1(\xi), e_{2k}(\xi) \rangle$ in the first statement \eqref{Moyal_indep}.
\begin{proof}[Proof of Proposition \ref{prop_Moyal_curve_Rn}, Moyal's Identity for Curves]
Fix $\xi \in I$, and $v_b \in \langle e_1(\xi), e_{2k}(\xi) \rangle^{\perp}$.
As in the proof of the Sobolev--Stein inequality, changing variables via the map $\gamma_T: I \to \tilde{I}$, we have
\begin{align}
    W_{2k}(g_1, g_2)(\xi, v_t + v_b) 
    &= \int_{\tilde{I}} G_{b}(\eta(\tilde{\eta}))  e^{-2\pi i v_t \cdot \tilde{\eta}e_1(\xi)} \rmd \tilde{\eta},
\end{align}
where
\begin{equation}
\label{G_b_Rn}
    G_{b}(\tilde{\eta}) := \frac{g_1(\eta(\tilde{\eta})) \overline{g_2(R_{\xi}\eta(\tilde{\eta}))} J(\xi,\eta(\tilde{\eta}))) e^{-2\pi i v_b \cdot (\phi(\eta(\tilde{\eta}))) - \phi(\Rmap(\tilde{\eta})))}}{\tilde{J}(\xi, \eta(\tilde{\eta})))}=: G(\tilde{\eta})e^{-2\pi i v_b \cdot (\phi(t)) - \phi(\Rmap(\tilde{\eta})))}.
\end{equation}
so that $|G_b(\tilde{\eta})| = |G(\tilde{\eta})|$ for all $\tilde{\eta} \in \tilde{I}$.
Moreover, parametrising $v_t = v_T e_1(\xi)$ for $v_T \in \mathbb{R}$,
\begin{align}
    W_{2k}(g_1, g_2)(\xi, v_T e_1(\xi) + v_b) 
    &= \int_{\tilde{I}} G_{b}(\tilde{\eta})  e^{-2\pi i v_T \tilde{\eta}} \rmd \tilde{\eta}
    = \widehat{G_b}(v_T).
\end{align}
Defining similarly $F_b$ and $F$ by replacing $g_1, g_2$ with $f_1, f_2$, we proceed by denoting $v_t = P_{\langle e_1(\xi) \rangle}(v)$, ${v_b = P_{\langle e_1(\xi), e_{2k}(\xi) \rangle^{\perp}}(v)}$, using Parseval's theorem, and reversing the third change of variables:
\begin{align}
    \int_{\langle e_1(\xi) \rangle} W_{2k}(f_1, f_2)(\xi, v_t + v_b) \overline{W_{2k}(g_1, g_2)(\xi, v_t + v_b)} \rmd v_t
    &= \int_{\mathbb{R}} \widehat{F_b}(v_T) \overline{\widehat{G_b}(v_T)} \rmd v_T \\
     &= \int_{\tilde{I}}F(\tilde{\eta}) \overline{G(\tilde{\eta})} \rmd \tilde{\eta}  \\
    &= \int_{I}f_1(\eta) \overline{f_2(R_{\xi}\eta)} \; \overline{g_1(\eta)}   g_2(R_{\xi}\eta)
    \frac{ J(\xi,\eta)^2}{\tilde{J}(\xi, \eta)} \rmd \eta,
\end{align}
which is precisely \eqref{Moyal_indep}.
To prove the second equality of the theorem, we use Fubini's theorem, the change of variables $\xi \mapsto \Rmap$, expression \eqref{explicit_J_formula_n} and Lemma \ref{lemma_Jtilde_identity} as follows
\begin{align*}
   \int_I \int_{\langle e_1(\xi) \rangle} W_{2k}&(f_1, f_2)(\xi, v) \overline{W_{2k}(g_1, g_2)(\xi, v)} \rmd \mathcal{H}^1(v) \rmd \xi\\
     &= \int_I f_1(\eta) \overline{g_1(\eta)} \left(\int_{\tilde{I}}\overline{f_2(R_{\xi}\eta)}  g_2(R_{\xi}\eta)
    \frac{ J(\xi,\eta)^2}{\tilde{J}(\xi, \eta)} \rmd \xi \right)\rmd \eta  \\
    &=\int_I f_1(\eta)\overline{g_1(\eta)} \left(\int_I \overline{f_2(\zeta)} g_2(\zeta) M(\eta, \zeta) \rmd \zeta \right) \rmd \eta.
\end{align*}
\end{proof}
\noindent
Now, we prove the orthonormal Mizohata--Takeuchi inequality for curves via a Wigner approach.
\begin{proof}[Proof of Theorem \ref{thm_orthonormal_MT_Rn_curve}]
    By duality, it is equivalent to show
    \begin{equation}
    \int_{\mathbb{R}^n} \sum_j \lambda_j
        |E_Sg_j(x)|^2w(x) \rmd x
        \lesssim \|(\lambda_j)\|_{l^2}\left\|\left(\int_{\langle e_1(\xi), e_{2k}(\xi) \rangle^{\perp}}\|X_{2k}w(\xi, \cdot + v_b)\|_{L^2_{v_t}} \rmd v_b \right) \right\|_{L^2_{\xi}(K^*)}.
\end{equation}
We turn our attention to the left-hand side of the desired inequality, use the Wigner representation, and decompose $v = v_t + v_b$ with the usual notation to give
\begin{align*}
    \int_{\mathbb{R}^n} \sum_j \lambda_j
        |E_Sg_j(x)|^2w(x) \rmd x
         &=\int_I \int_{\langle e_1(\xi), e_{2k}(\xi) \rangle^{\perp}}\int_{\langle e_1(\xi) \rangle} \sum_j \lambda_j W_{2k}g_j(\xi, v_t + v_b)X_{2k}w(\xi, v_t + v_b) \rmd v_t \rmd v_b \rmd\xi.
\end{align*}
Hence by the Cauchy--Schwarz inequality, we see that
\begin{align*}
\label{proof_orthonormal_MT_Wigner_eqn2}
\int_{\mathbb{R}^n} \sum_j \lambda_j
|E_Sg_j(x)|^2w(x) \rmd x 
        &\leq \int_{I}\int_{\langle e_1(\xi), e_{2k}(\xi) \rangle^{\perp}} \left\|\sum_j \lambda_j W_{2k}g_j(\xi, \cdot + v_b)\right\|_{L^2_{v_t}}
    \|X_{2k}w(\xi, \cdot + v_b)\|_{L^2_{v_t}} \rmd v_b \rmd \xi 
\end{align*}
Now, recall that by \eqref{Moyal_indep},
\begin{align}
    \left\|\sum_j \lambda_j W_{2k}g_j(\xi, \cdot + v_b)\right\|_{L^2_{v_t}} 
    &= \left(\sum_i \sum_j \left(\lambda_i \overline{\lambda_j}\int_{I}g_i(\eta) \overline{g_i(R_{\xi}\eta)} \; \overline{g_j(\eta)}   g_j(R_{\xi}\eta)
    \frac{ J(\xi,\eta)^2}{\tilde{J}(\xi, \eta)} \rmd \eta \right)\right)^{1/2} := H(\xi).
\end{align}
Hence, using this highlighted lack of dependence on $v_b \in \langle e_1(\xi), e_{2k}(\xi) \rangle^\perp$; noticing that $W_{2k}g_j(\xi,v)=0$ for all $\xi \notin \supp^*(g_j)$; and by another application of the Cauchy--Schwarz inequality,
\begin{align*}
\int_{\mathbb{R}^n} \sum_j \lambda_j
|E_S&g_j(x)|^2w(x) \rmd x
    \leq \int_{I} H(\xi)
    \left(\int_{\langle e_1(\xi), e_{2k}(\xi) \rangle^{\perp}}\|X_{2k}w(\xi, \cdot + v_b)\|_{L^2_{v_t}} \rmd v_b \right) \rmd \xi \\
    &\leq  \|H\|_{2}
    \left\|\int_{\langle e_1(\xi), e_{2k}(\xi) \rangle^{\perp}}\|X_{2k}w(\xi, \cdot + v_b)\|_{L^2_{v_t}} \rmd v_b \right\|_{L^2_{\xi}(K^*)}.
\end{align*}
Noticing that by \eqref{Moyal_indep}, and the assumption that $L(i, j)$ is the kernel of an $l^2$ bounded operator
\begin{align}
    \|H\|_{L^2(I)}^2 &= \sum_i \sum_j \lambda_i \overline{\lambda_j} \int_I \int_{\langle e_1(\xi) \rangle}  W_{2k}g_i(\xi, v_t + v_b) \overline{W_{2k}g_j(\xi, v_t + v_b)} \: \rmd \mathcal{H}^1(v) \rmd \xi
    = \sum_i \sum_j \lambda_i \overline{\lambda_j}L(i, j)
    \leq \|(\lambda_j)\|_{l^2}^2.
\end{align}
Hence,
\begin{align*}
    \int_{\mathbb{R}^n} \sum_j \lambda_j
        |E_Sg_j(x)|^2w(x) \rmd x
        \leq \left\|(\lambda_j)\right\|_{l^2}
    \left\|\int_{\langle e_1(\xi), e_{2k}(\xi) \rangle^{\perp}}\|X_{2k}w(\xi, \cdot + v_b)\|_{L^2_{v_t}} \rmd v_b \right\|_{L^2_{\xi}(K^*)}.
\end{align*}
\end{proof}

\section{The Helix in \texorpdfstring{$\mathbb{R}^n$}{Rn}}
\label{Section_helix}
When all of the Frenet curvatures, $\kappa_1, ..., \kappa_{n-1}$ are constant functions, there is a simple way to find the Frenet vectors using the Frenet formulae. This nicely constructed example, categorised by its constant Frenet curvatures, is \textit{a helix in $\mathbb{R}^n$}.

Let $a_i, \alpha_i, b \in \mathbb{R}$ for $i \in \{1, ..., m\}$, be such that each $\alpha_i$ is distinct and $\sum_{i=1}^m a_i^2\alpha_i^2 + b^2 =1$ (where $b=0$ when the ambient dimension is even).
When the ambient dimension is even, the arclength parametrisation $\phi_e: I \to \mathbb{R}^{2m}$ of a helix satisfies
\begin{equation}
\label{helix_even}
    \phi_e(\xi) = (a_1 \sin(\alpha_1 \xi), a_1\cos(\alpha_1 \xi), ..., a_m \sin(\alpha_m \xi), a_m\cos(\alpha_m \xi)),
\end{equation}
for all $\xi \in I$.
On the other hand, when the ambient dimension is odd, the arclength parametrisation $\phi_o: I \to \mathbb{R}^{2m+1}$ of the helix is such that
\begin{equation}
\label{helix_odd}
    \phi_o(\xi) = \left(a_1 \sin(\alpha_1 \xi), a_1\cos(\alpha_1 \xi), ..., a_m \sin(\alpha_m \xi), a_m\cos(\alpha_m \xi), b\xi\right),
\end{equation}
for all $\xi \in I$.
In the former case, the $j$-th Frenet vector, $e_j: I \to \mathbb{R}^{2m}$ of a helix in $\mathbb{R}^{2m}$ takes the form
\begin{equation}
\label{helix_ej_n_even}
e_j(\xi)=
  \begin{cases} 
    (c_{1j} \cos(\alpha_1 \xi), -c_{1j}\sin(\alpha_1 \xi), ..., c_{mj} \cos(\alpha_m \xi), -c_{mj}\sin(\alpha_m \xi)) \qquad \text{ for $j$ odd}
    \\
    (c_{1j} \sin(\alpha_1 \xi), c_{1j}\cos(\alpha_1 \xi), ..., c_{mj} \sin(\alpha_m \xi), c_{mj}\cos(\alpha_m \xi)) \qquad \quad \text{ for $j$ even},
  \end{cases}
 \end{equation}
 for all $\xi \in I$.
Moreover, the $j$-th Frenet vector, $e_j: I \to \mathbb{R}^{2m +1}$ of the helix in $\mathbb{R}^{2m+1}$ satisfies
\begin{equation}
\label{helix_ej_n_odd}
e_j(\xi)=
  \begin{cases} 
    (c_{1j} \cos(\alpha_1 \xi), -c_{1j}\sin(\alpha_1 \xi), ..., c_{mj} \cos(\alpha_m \xi), -c_{mj}\sin(\alpha_m \xi), d_{j})\qquad \text{ for $j$ odd}    \\
    (c_{1j} \sin(\alpha_1 \xi), c_{1j}\cos(\alpha_1 \xi), ...,c_{mj} \sin(\alpha_m \xi), c_{mj}\cos(\alpha_m \xi) , 0 )\qquad \quad \text{ for $j$ even},
  \end{cases}
 \end{equation}
 for all $\xi \in I$.
In both of these expressions, $c_{ij} \in \mathbb{R}$ are constants depending on $a_i, \alpha_i$ and the first $(j -1)$ Frenet curvatures, $\kappa_1, ..., \kappa_{j-1}$.
We comment that $c_{i1}= a_i\alpha_i$ for each $i \in \{ 1,..., \lfloor \frac{n}{2}\rfloor \}$. Moreover, $d_{1}=b$, $d_{2k}=0$, and $d_{2k+1}= \frac{\kappa_1 \kappa_3 \cdots \kappa_{2k-1}}{\kappa_2 \kappa_4 \cdots \kappa_{2k}}b$. Refer to example \eqref{example_Helix_R3} for the helix in $\mathbb{R}^3$ and to the book \cite{kuhnel} for further detail.

\subsection{Another helical Wigner transform}
For the remainder of this section, let $S$ denote the helix in $\mathbb{R}^n$ with arclength parametrisation $\phi:I \to \mathbb{R}^n$ and $I \subset \mathbb{R}$, satisfying $\phi = \phi_e$ for even ambient dimension and $\phi_o$ for odd ambient dimension for some fixed appropriate $a_i, \alpha_i, b \in \mathbb{R}$. This has proven to be an important example in Section \ref{Section_curvesR3}, and Section \ref{Section_curvesRn} provides a suitable Wigner transform associated to this curve in arbitrary dimensions. 
\begin{lemma}
\label{system_of_eqn_helix_prop}
For fixed $\xi, \eta \in I$,
    \begin{equation}
    e_{2k}(\xi) \cdot(\phi(\eta)- \phi(2\xi - \eta))=0 \qquad \text{for all $k \in \mathbb{N}$ with $2k \leq n$.}
 \end{equation}
\end{lemma}
\begin{proof}
Let $k \in \mathbb{N}$ be such that $2k \leq n$. Let also $\xi, \eta, \zeta \in I$ and observe that
\begin{equation}
    e_{2k}(\xi)\cdot (\phi(\eta)- \phi(\zeta))= -2\sum_{i=1}^{m} a_i c_{ij}\sin\left(\frac{\alpha_i}{2}(\eta - \zeta)\right)\sin\left(\frac{\alpha_i}{2}(\eta + \zeta - 2\xi)\right),
\end{equation}
Taking $\zeta = 2\xi - \eta$, it is immediate that
    \begin{equation}
   e_j(\xi)\cdot (\phi(\eta)- \phi(2\xi - \eta))= 0.
\end{equation}
\end{proof}

Therefore, fixing $k \in \mathbb{N}$ to be such that $2k \leq n$, and letting $g_1, g_2 \in L^2(I)$, the \textit{$k^{\text{th}}$ helical Wigner transform} is given by
\begin{equation}
    W_{2k}(g_1, g_2)(\xi,v) = 2\int_{I}g_1(\eta) \overline{g_2(2\xi - \eta)} e^{-2\pi iv\cdot (\phi(\eta) - \phi(2\xi - \eta))}\rmd \eta,
\end{equation}
where $(\xi, v) \in D_{2k}$.
 
Motivated by Lemma \ref{system_of_eqn_helix_prop}, we will see that it is also possible to construct another natural and useful Wigner transform associated to the helix. We will prove another Wigner representation for the extension operator for the helix which results in a stronger Sobolev--Mizohata--Takeuchi inequality. To do so, we introduce some notation.
Fix $\xi \in I$ and define
\begin{align}
     Y_e:= Y_e(\xi) &:= \langle\{ e_{2k}(\xi): k \in \mathbb{N}, 2k \leq n\}\rangle; \\
     Y_o := Y_o(\xi) &:= \langle\{ e_{2k - 1}(\xi): k \in \mathbb{N}, 2k \leq n\}\rangle; \\
     Y_b:= Y_b(\xi) &:= \langle\{ e_{2k + 1}(\xi): k \in \mathbb{N}, 2k \leq n\}\rangle.
\end{align}
Then, $Y_e(\xi)$ is the span of the collection of all even Frenet vectors at $\phi(\xi)$, $Y_o(\xi)$ is the span of the collection of all odd Frenet vectors, and $Y_b(\xi)$ is the span of the collection of all odd Frenet vectors except for the tangent vector $e_1(\xi)$.

For $g_1, g_2 \in L^2(I)$, we define the \textit{helical Wigner transform}, via
\begin{equation}
    W_{e}(g_1, g_2)(\xi,v_o) = 2\int_{I}g_1(\eta) \overline{g_2(2\xi - \eta)} e^{-2\pi iv_o\cdot (\phi(\eta) - \phi(2\xi - \eta))}\rmd \eta,
\end{equation}
where $(\xi, v_o) \in \{(\xi, v_o): \xi \in I, v_o \in Y_o(\xi)\}:=D_o$.

The new helical Wigner transform also satisfies the typical properties, some of which it inherits from the $k^{\text{th}}$ helical Wigner transform and others which follow similarly.
\begin{proposition}
\label{conj_symm_helix_prop}
    The Helical Wigner Transform is conjugate symmetric in that it satisfies
    \begin{align*}
        W_e(g_1, g_2)(\xi, v_o) = \overline{W_e(g_2, g_1)(\xi, v_o)} \qquad \text{for all } (\xi, v_o) \in D_o.
    \end{align*}
\end{proposition}
\begin{proof} 
Fix $\xi \in I$ and change variables via the map $\eta \mapsto 2\xi - \eta$:
    \begin{align*}
        W_e(g_1, g_2)(\xi, v_o) =
        2\int_{I}g_1(2\xi - \zeta) \overline{g_2(\zeta)} e^{-2\pi iv_o\cdot (\phi(2\xi - \zeta) - \phi(\zeta))}\rmd \zeta = \overline{W_e(g_2, g_1)(\xi, v_o)}.
    \end{align*}
\end{proof}
Define the auxiliary function $F: I \times \mathbb{R}^n \to \mathbb{R}$ as follows
\begin{equation}
    F(\xi, x) = 2\int_{I}g(\eta) \overline{g(2\xi - \eta)} e^{-2\pi ix\cdot (\phi(\eta) - \phi(2\xi - \eta))}\rmd \eta,
\end{equation}
so that $W_eg$ is the restriction of $F$ to $D_o$. Observe that this is precisely the auxiliary function defined by \eqref{aux_fun_F_k} in the case of the helix. On this note, we also mention that $W_eg$ is the restriction of $W_{2k}g$ to $D_o$. Explicitly, for fixed $\xi \in I$,
\begin{equation}
\label{helix:aux_fun_We}
    F(\xi, x) = W_eg(\xi, P_{Y_o(\xi)}(x)) = W_{2k}g(\xi, P_{\langle e_{2k(\xi)\rangle}}(x)), \qquad x \in \mathbb{R}^n,
\end{equation}
where the second equality is precisely \eqref{defn_F_k}.
The relationship between the two Wigner transforms for the helix is further emphasised by the following equality
\begin{equation}
\label{helix:W2k_We}
    W_{2k}g(\xi, v) = W_eg(\xi, P_{Y_o(\xi)}(v)), \qquad v \in \langle e_1(\xi), e_{2k}(\xi) \rangle^\perp.
\end{equation}
Since $F$ is the same auxiliary function as \eqref{aux_fun_F_k}, the first marginal property is precisely Proposition \ref{first_marg} in the case of the helix (that is, with $J(\xi, \eta) =2)$.

Now, define a restriction of $Tf$, the $\lfloor n/2 \rfloor$-plane transform (see \eqref{k_plane_defn}). Indeed, for $(\xi, v_o) \in D_o$,
\begin{equation}
    T_ef(\xi, v_o) := Tf(Y_e(\xi), v_o) = \int_{Y_e(\xi)}f(v_e + v_o) \rmd v_e.
\end{equation}

\begin{theorem}[\textbf{Wigner representation}]  Let $S$ be the helix in $\mathbb{R}^n$ with arclength parametrisation $\phi: I \to \mathbb{R}^n$, and $I \subset \mathbb{R}$. Moreover, let $g \in L^2(I)$. Then,
\begin{equation}
    |E_Sg|^2 = T^{*}_e W_eg,
\end{equation}
where for functions $h: D_o \to \mathbb{C}$, 
\begin{equation}
\label{defn_Te*}
    T^{*}_e h(x) := \int_I h(\xi, P_{Y_o(\xi)}(x)) \rmd \xi.
\end{equation}
\end{theorem}
Again, it may be observed that $T^{*}_e$ is the adjoint of $T_e$, and the theorem follows from the first marginal property, \eqref{helix:aux_fun_We}, and \eqref{defn_Te*}.
Moreover, we have the following integral representation for all weight functions $w: \mathbb{R}^n \to \mathbb{C}$,
\begin{align*}
\int_{\mathbb{R}^n}|E_Sg(x)|^2w(x) \rmd x =  \int_I \int_{ Y_o(\xi) } W_eg(\xi,v_o) \:
T_ew(\xi, v_o) \: \rmd v_o \: \rmd\xi.
\end{align*}

\subsection{Motivating applications via Wigner heuristics}
The helical Wigner transform satisfies a tangent marginal property:
 \begin{equation}
        \int_{\langle e_1(\xi) \rangle} W_e(g_1, g_2)(\xi, v_t + v_o') \rmd v_t = g_1(\xi)\overline{g_2(\xi)} \qquad \text{for all } \xi \in I, v_o' \in Y_b(\xi).
    \end{equation}
This follows from the tangent marginal property for the $k^{\text{th}}$-Wigner transform by noticing that 
\begin{equation}
    W_e(g_1, g_2)(\xi, v_t + P_{Y_b(\xi)}(v_b)) = W_{2k}(g_1, g_2)(\xi, v_t + v_b), \qquad \xi \in I, v_b \in \langle e_1(\xi), e_{2k}(\xi) \rangle^\perp.
\end{equation}
Similarly, 
\begin{equation}
\label{helix:sobolev_equiv}
    \|W_e(g_1, g_2)(\xi, \cdot + v_o')\|_{\dot{H}^{-s}_{v_t}}= \|W_{2k}(g_1, g_2)(\xi, \cdot + v_b)\|_{\dot{H}^{-s}_{v_t}},
\end{equation}
for all $\xi \in I, v_b \in \langle e_1(\xi), e_{2k}(\xi)\rangle^{\perp}, v_o' \in Y_b(\xi)$.

Both the tangent marginal and the homogeneous Sobolev norm over the tangent direction exhibit an independence of $v_o' \in Y_b(\xi)$, inherited from the work for general curves in $n$-dimensions.

Consequently, the Wigner heuristics would suggest that if $W_eg$ were non-negative then,
\begin{align*}
    \int_{\mathbb{R}^n} |E_Sg(x)|^2w(x) \rmd x 
    \leq
\|g\|_{L^2(I)}^2 \sup_{\xi}\int_{Y_b(\xi)} \sup_{v_t} T_ew(\xi, v_t + v_o') \rmd v_o' .
\end{align*}

\begin{remark}[Proximity to \eqref{MT_conj_dim_m}]
Observe that
\begin{equation} \mathcal{R}_Tw(\xi, v_t) = \int_{Y_b(\xi)}\int_{Y_e(\xi)}  w(v_t + v_o' + v_e)   \rmd v_e  \rmd v_o'= \int_{Y_b(\xi)}  T_ew(\xi, v_t + v_o')     \rmd v_o' .
\end{equation}
   Notice
   \begin{equation}
       \sup_{\xi} \sup_{v_t}\int_{Y_b(\xi)}  T_ew(\xi, v_t + v_o')     \rmd v_o'  \leq \sup_{\xi}\int_{Y_b(\xi)} \sup_{v_t} T_ew(\xi, v_t + v_o')     \rmd v_o' \leq \sup_{\xi}\int_{\langle e_1(\xi), e_{2k}(\xi)\rangle^{\perp}} \sup_{v_t} X_{2k} w(\xi, v_t + v_o')     \rmd v_o'.
   \end{equation}
   This demonstrates that for the helix, we may aim to prove something stronger than for a general curve and closer to  \eqref{MT_conj_dim_m}. Again, we cannot assume that $W_{e}g$ is non-negative, but we will be motivated by the heuristics to prove $L^2$-Sobolev variants of this form.
\end{remark}
For the upcoming applications of the helical Wigner representation, it will simplify computations to impose the tangents condition \eqref{eqn_tangents_condition_n}. This can be expected to hold for a suitably chosen interval $I$, for example for the helix in three dimensions (see Example \ref{example_Helix_R3}) taking $|I| < \frac{\pi}{3\lambda}$ suffices.
\subsection{Application 1: Weighted \texorpdfstring{$L^2$}{L2} inequalities}
\begin{theorem}[$L^2$ Sobolev--Stein Inequality for the Helix in $\mathbb{R}^n$]
\label{appliation2a_l2sobolev_helix}
Let $S$ also satisfy \eqref{eqn_tangents_condition_n}. Moreover, let $s < \frac{1}{2}$ and $g \in L^2(I)$. Then,
    \begin{equation}
    \int_{\mathbb{R}^n} |E_Sg(x)|^2w(x) \rmd x 
    \leq 4(2\pi)^{2s} \int_I \mathcal{B}_{2s}(|g|^2, |g|^2)(\xi)^{\frac{1}{2}} \left( \int_{\langle Y_b(\xi)\rangle} \left\| T_e w(\xi, \cdot + v_o')) \right\|_{\dot{H}_{v_t}^s} \rmd v_o' \right)  \rmd\xi,
\end{equation}
where
\begin{equation}
    \mathcal{B}_{s}(g_1, g_2)(\xi) = \int_{I}\frac{g_1(\eta)\overline{g_2(2\xi - \eta)}}{|\eta - \xi|^{s}}\rmd \eta.
\end{equation}
\begin{remark}
\label{same_as_Kenig_Stein}
Observe that via the change of variables $\eta \mapsto \eta + \xi$,
    \begin{equation}
    \mathcal{B}_{s}(g_1, g_2)(\xi) = \int_{I'}\frac{g_1(\xi + \eta)\overline{g_2(\xi - \eta)}}{|\eta|^{s}}\rmd \eta.
\end{equation}
This is a local version of the classical bilinear fractional integral studied by Kenig and Stein in Theorem $2$ of \cite{KenigStein}.
\end{remark}
\end{theorem}
\begin{proof}[Proof of Theorem \ref{appliation2a_l2sobolev_helix}]
We again begin by applying the Wigner representation, and by writing ${v_o = v_t + v_o'}$, where $v_t = P_{\langle e_1(\xi)\rangle}(v_o)$ and $v_o' = P_{\langle Y_{b}(\xi) \rangle}(v_o)$:
\begin{align*}
    \int_{\mathbb{R}^n} |E_Sg(x)|^2w(x) \rmd x
    &=  \int_I \int_{Y_b(\xi)}\int_{\langle e_{1}(\xi) \rangle} W_eg(\xi, v_t + v_o')T_ew(\xi, v_t + v_o') \rmd v_t \rmd v_o' \rmd\xi \\
    &\leq \int_I \int_{Y_b(\xi)}
    \|W_eg(\xi,\cdot + v_o')\|_{\dot{H}^{-s}_{v_t}} \|T_ew(\xi, \cdot + v_o')\|_{\dot{H}^{s}_{v_t}}
    \rmd v_o'
    \rmd \xi. \qquad (\star)
\end{align*} 
\noindent
Fixing $\xi \in I$ and $v_o' \in Y_b(\xi)$, we seek to analyse the first homogeneous Sobolev norm.
By \eqref{helix:sobolev_equiv},
\begin{align*}
    \| W_eg(\xi,\cdot + v_o')\|_{\dot{H}^{s}_{v_t}}^2 
     &=4(2\pi)^{-2s}
    \int_{I}\frac{|g(\eta)|^2  |g(2\xi - \eta)|^2 }{\tilde{J}(\xi, \eta)|e_1(\xi) \cdot(\phi(\eta) - \phi(2\xi - \eta))|^{2s}}\rmd \eta.
\end{align*}
By Lemma \ref{lem:gammaT_sim_a-b_ndim},
\begin{equation}
    |e_1(\xi) \cdot(\phi(\eta) - \phi(2\xi - \eta))| \geq |\eta - \xi|.
\end{equation}
Moreover, \eqref{eqn_tangents_condition_n} implies
   \begin{equation}
   \label{J_tilde_helix}
        \tilde{J}(\xi, \eta)= e_1(\xi)\cdot \left(e_1(\eta) + e_1(2\xi - \eta)\right) \geq 1
    \end{equation}
for all $\xi, \eta \in I$.
Hence,
\begin{align*}
    \| W_eg(\xi,\cdot + v_o'))\|_{\dot{H}^{s}_{v_t}}^2
     &\leq 4(2\pi)^{-2s}
    \int_{I}\frac{|g(\eta)|^2  |g(2\xi - \eta)|^2 }{|\eta - \xi|^{2s}}\rmd \eta.
\end{align*}
Returning to $(\star)$, we conclude by using the lack of dependence on $v_o' \in Y_b(\xi)$: 
\begin{align*}
    \int_{\mathbb{R}^n} |E_Sg(x)|^2w(x) \rmd x
    &= 4(2\pi)^{-2s} \int_I \left( \int_{I}\frac{|g(\eta)|^2  |g(2\xi - \eta)|^2 }{|\eta - \xi|^{2s}}\rmd \eta \right)^{\frac{1}{2}} \left( \int_{Y_b(\xi)} \left\| T_ew(\xi, \cdot + v_o')) \right\|_{\dot{H}_{v_t}^s} \rmd v_o' \right)  \rmd\xi.
\end{align*}

\end{proof}
\begin{theorem}[$L^2$ Sobolev--Mizohata--Takeuchi Inequality for the Helix in $\mathbb{R}^n$]
\label{corol2b_l2sobolev_helix}
Assume $S$ satisfies \eqref{eqn_tangents_condition_n}, and let $s < \frac{1}{2}$ and $g \in L^2(I)$. Then there exists a constant $c_s$ only depending on $s$ such that
    \begin{equation}
    \int_{\mathbb{R}^n} |E_Sg(x)|^2w(x) \rmd x 
    \leq c_s \|g\|_{L^2(I)}^2 \sup_{\xi \in \supp^*(g)} \left( \int_{Y_b(\xi)} \left\| T_ew(\xi, \cdot + v_o')) \right\|_{\dot{H}_{v_t}^s} \rmd v_o' \right).
\end{equation}
Here, $\supp^*(g) 
$ is the set of all midpoints of pairs of points from $\supp(g)$.
\end{theorem}
\begin{proof}[Proof of Theorem \ref{corol2b_l2sobolev_helix}]
The statement follows from the Sobolev--Stein inequality, by using Remark \ref{same_as_Kenig_Stein} and the bilinear fractional integral estimates of Kenig and Stein in \cite{KenigStein}. 
\end{proof}

\subsection{Application 2: Orthonormal Mizohata--Takeuchi inequality for the helix}
We also may prove an orthonormal Mizohata--Takeuchi estimate and a Moyal identity for the helix with the improved Wigner representation. We state the results here, but omit their proofs since they follow similarly to those of Section \ref{sec:orthonormal} with analogous changes to those in the proof of the Sobolev--Stein and Sobolev--Mizohata--Takeuchi inequalities.
\begin{theorem}
\label{orthonormal_MT_helix}
Let $S$ satisfy \eqref{eqn_tangents_condition_n}. Moreover, let $(g_j)$ be a sequence of functions on $I$ for which 
\begin{equation}
    L(i, j) := \int_I \int_{\langle e_1(\xi) \rangle} W_eg_i(\xi, v_o) \overline{W_eg_j(\xi, v_o)} \rmd \mathcal{H}^1(v_o) \rmd \xi
\end{equation}
is the (real-valued) kernel of an $l^2$-bounded operator. Then, for all weight functions $w: \mathbb{R}^3 \to \mathbb{C}$,
    \begin{equation}
        \sum_j \left(
        \int_{\mathbb{R}^n}|Eg_j(x)|^2w(x) \rmd x
        \right)^2 \lesssim 
    \|T_ew(\xi, v_t + v_o')\|_{L^2_{\xi}(K^*)L^1(Y_b)L^2(\langle e_1(\xi) \rangle)}^2.
    \end{equation}
\end{theorem}
\noindent
To understand what it means for $L(i, j)$ to be the kernel of an $l^2$-bounded operator, we turn to Moyal's identity for the helical Wigner transform.
\begin{proposition}[Moyal's Identity for the Helical Wigner Transform]
\label{helical_Moyal}
$S$ satisfy \eqref{eqn_tangents_condition_n}. Moreover, let $f_1, f_2, g_1, g_2 \in L^2(I)$. Then,
    \begin{align}
    \int_I \int_{\langle e_1(\xi) \rangle} W_e(f_1, f_2)(\xi, v_o) \overline{W_e(g_1, g_2)(\xi, v_o)} \rmd \mathcal{H}^1(v_o) \rmd \xi =
        \int_I f_1(\eta)\overline{g_1(\eta)} \left(\int_I \frac{\overline{f_2(\zeta)} g_2(\zeta) }{\tilde{J}\left(\frac{\eta + \zeta}{2}, \eta\right)} \rmd \zeta \right) \rmd \eta.
    \end{align} 
\end{proposition}
Recall $\tilde{J}(\xi, \eta)$ is given by \eqref{J_tilde_helix} and is bounded below by $1$ for the helix under the hypothesis \eqref{eqn_tangents_condition_n}.

\section{The Second Marginal Property for curves in \texorpdfstring{$\mathbb{R}^n$}{Rn}}
\label{Section_second_marginal}
In this section, we formalise and prove the second marginal property of the $k^{\text{th}}$-curve-carried Wigner transform for the sake of completeness. Its proof is similar to that of the second marginal property of the surface-carried Wigner transform in \cite{BGNO}.
\begin{proposition}[The Second Marginal Property]
\label{Second_marginal_property_curve_Rn}
Let $S$ be a well-curved, compact patch of a $C^{n+2k -1}$ curve in $\mathbb{R}^n$ with arclength parametrisation, $\phi: I \to \mathbb{R}^n$ satisfying conditions \eqref{eqn_ej_dot_T_not_zero_condition},  \eqref{eqn_taylor_condn}, 
\eqref{eqn_ej'_dot_T_not_zero_condition},
and \eqref{eqn_tangents_condition_n}.
    Furthermore, let $\chi$ be a Schwartz function on $\langle e_1(\xi) \rangle$ for which $\chi(0)=1$ and $\chi_r(v_t)=\chi(\frac{v_t}{r})$ for each $r >0$. Then for continuous, compactly supported $g_1, g_2 \in L^2(I)$,
    \begin{equation}
    \int_{\langle e_1(\xi) \rangle} W_{2k}(g_1, g_2)(\xi, v_t + v_b) \chi_r(v_t) \rmd v_t \to g_1(\xi)\overline{g_2(\xi)}
    \end{equation}
    as $r \to \infty$ for all $\xi \in I$ and $v_b \in \langle e_1(\xi), e_{2k}(\xi) \rangle^\perp$.
\end{proposition}
\begin{proof}[Proof of the Second Marginal Property, Proposition \ref{Second_marginal_property_curve_Rn}]
    Fix $\xi \in I$, $v_b \in \langle e_1(\xi), e_{2k}(\xi) \rangle^{\perp}$. Parametrise $v_t = v_T e_1(\xi)$, for $v_T\in \mathbb{R}$ and observe that
    \begin{align*}
        \int_{\langle e_1(\xi) \rangle}& W_{2k}(g_1, g_2)(\xi, v_t + v_b) \chi_r(v_t) \rmd v_t \\
        &=
        \int_{\langle e_1(\xi) \rangle} \left(
        \int_{I}g_1(\eta) \overline{g_2(\Rmap)} J(\xi, \eta) e^{-2\pi i(v_t + v_b)\cdot (\phi(\eta) - \phi(\Rmap))}\rmd \eta  \right) \chi_r(v_t)
        \rmd v_t \\
        &=
        \int_{\mathbb{R}} \left(
        \int_{I}g_1(\eta) \overline{g_2(\Rmap)} J(\xi, \eta) e^{-2\pi iv_b\cdot (\phi(\eta) - \phi(\Rmap))}\rmd \eta\right)e^{-2\pi iv_Te_1(\xi)\cdot (\phi(\eta) - \phi(\Rmap))} \chi_r(v_T)
        \rmd v_T \\
        &= 
        \int_{I}g_1(\eta) \overline{g_2(\Rmap)} J(\xi, \eta) 
        e^{-2\pi iv_b\cdot (\phi(\eta) - \phi(\Rmap))}
        \hat{\chi}_r(e_1(\xi)\cdot (\phi(\eta) - \phi(\Rmap)))\rmd \eta \\
        &:= \mathcal{A}_{S, r}(g_1, g_2)(\xi).
    \end{align*} 
Since $g_1, g_2$ are continuous functions, it suffices to prove that $\mathcal{A}_{S, r}(1, 1)(\xi) \to 1$.
To this end, we change variables via the third map and then again via the map $z \mapsto r^{-1}x$:
    \begin{align*}
        \mathcal{A}_{S, r}(1, 1)(\xi)
        &= r 
        \int_{I} J(\xi, \eta)
        e^{-2\pi iv_b\cdot (\phi(\eta) - \phi(\Rmap))}
        \hat{\chi}(re_1(\xi)\cdot( \phi(\eta) - \phi(\Rmap)))\rmd \eta \\
        &= r 
        \int_{\tilde{I}}
        e^{-2\pi iv_b\cdot \gamma_{\xi}(\eta(z))} \frac{J(\xi,\eta(z))}{\tilde{J}(\xi,\eta(z))}
        \hat{\chi}(rz)\rmd z \\
        &=\int_{U}
        e^{-2\pi iv_b\cdot \gamma_{\xi}(\eta(r^{-1}x))} \frac{J(\xi,\eta(r^{-1}x))}{\tilde{J}(\xi,\eta(r^{-1}x))}
        \hat{\chi}(x)\rmd x.
    \end{align*}
As $r \to \infty$, $e_1(\xi) \cdot(\phi(\eta) - \phi(\Rmap)) = r^{-1}x \to 0$.
Lemmas \ref{gammaT_sim_a-b_ndim} and \ref{phia-phib_sim_a-b_ndim} together imply that 
$$e_1(\xi) \cdot(\phi(\eta) - \phi(\Rmap)) \lesssim \eta - \xi \lesssim |\phi(\eta)- \phi(\Rmap)|.$$
Hence as $e_1(\xi)\cdot (\phi(\eta) - \phi(\Rmap)) \to 0$, $|\phi(\eta) - \phi(\Rmap)| \to 0$ and $\eta \to \xi$. This in turn implies that $v_b \cdot(\phi(\eta) - \phi(\Rmap)) \to 0$. 
To conclude, it remains to understand the limits of $J(\xi, \eta)$ and $\tilde{J}(\xi, \eta)$ as $\eta \to \xi$, which is precisely Lemma \ref{jacobian_limits_n} below.
\end{proof}
\begin{remark}
    The requirement for $g_1, g_2$ to be continuous can be removed by proving a maximal function estimate. We expect that this would follow in a similar way to the maximal function estimate proven by \cite{BGNO} with suitable extensions made, provided one could prove further distance estimates.
\end{remark}
\begin{lemma}
For each $\xi \in I$, we have 
\label{jacobian_limits_n}
    \begin{align}
    \label{limJ}
        \lim_{\eta \to \xi}J(\xi, \eta) &= 2, \\
    \label{limDelta}
        \lim_{\eta \to \xi}\Delta(\xi, \eta) &= 1, \\
    \label{limJtilde}
        \lim_{\eta \to \xi}\tilde{J}(\xi, \eta) &= 2.
    \end{align}
\end{lemma}
\begin{remark}
    Recall that 
    \begin{align}
        J(\xi, \eta) &=  \left|\frac{e_{2k}'(\xi))\cdot (\phi(\eta) - \phi(R_{\xi}\eta)}{e_{2k}(\xi)\cdot e_{1}(R_{\xi}\eta)}\right|;\\
        \Delta(\xi, \eta) &= \left|\frac{e_{2k}(\xi)\cdot e_1(\eta)}{e_{2k}(\xi)\cdot e_{1}(R_{\xi}\eta)}\right|;\\
        \tilde{J}(\xi, \eta) & = |e_{1}(\xi)\cdot \left(e_{1}(\eta) + \Delta(\xi, \eta) e_{1}(R_{\xi}\eta)\right)|.
    \end{align}
\end{remark}
\begin{proof}
Fix $\xi \in I$, and let $\zeta = R_{\xi}(\eta)$. 
Define $f: I \to \mathbb{R}$ to be such that
    \begin{equation}
        f(t) = e_{2k}(\xi)\cdot\phi(t),
    \end{equation}
    so that \eqref{eqn_THE_defining_condition_R_n} stipulates that $f(\eta) = f(\zeta)$,
Hence, by Corollary \ref{corol_derivatives_frenet_vectors}, the Taylor expansion of $f$ about the point $\xi$ is
\begin{align*}
    f(t) - f(\xi)
    &=(t- \xi)^{2k}\left(\frac{(2k+1)A(\xi) + f^{(2k+1)}(a)(t- \xi)}{(2k+1)!}\right),
\end{align*}
for some $a$ between $t$ and $\xi$. For convenience in what follows, we have denoted
\begin{equation}
    A(\xi) := \prod_{m=1}^{2k-1}\kappa_{m}(\xi).
\end{equation}
Therefore, $f(\eta) = f(\zeta)$ implies
\begin{align*}
(\eta- \xi)^{2k}\left((2k+1)A(\xi) + f^{(2k+1)}(a_{\eta})(\eta- \xi)\right) = (\zeta- \xi)^{2k}\left((2k+1)A(\xi) + f^{(2k+1)}(a_{\zeta})(\zeta- \xi)\right),
\end{align*}
for $a_{\eta}$ between $\xi$ and $\eta$ and $a_{\zeta}$ between $\xi$ and $\zeta$.
Moreover, assumption \eqref{eqn_taylor_condn} implies that \\ $|f^{(2k+1)}(a)||t- \xi| < 2k |A(\xi)|$ for all $a, t, \xi \in I$. Therefore,
\begin{equation}
    |(2k+1)A(\xi) + f^{(2k+1)}(a)(t- \xi)| \geq (2k+1)|A(\xi)| - |f^{(2k+1)}(a)(t-\xi)| > |A(\xi)| \neq 0.
\end{equation}
Hence,
\begin{equation}
    \left(\frac{\eta- \xi}{\zeta - \xi}\right)^{2k} = \frac{(2k+1)A(\xi) + f^{(2k+1)}(a_{\zeta})(\zeta- \xi)}{(2k+1)A(\xi) + f^{(2k+1)}(a_{\eta})(\eta- \xi)}.
\end{equation}
Thus, as $\eta \to \xi$, we have $a_{\eta} \to \xi$, $\zeta \to \xi$, and $a_{\zeta} \to \xi$ by the in-between lemma. Hence, by assumption \eqref{eqn_taylor_condn}, we conclude that as $\eta \to \xi$,
\begin{equation}
    \left((2k+1)A(\xi) + f^{(2k+1)}(a_{\eta})(\eta- \xi)\right) \to (2k+1)A(\xi)  \text{  and  } \left((2k+1)A(\xi) + f^{(2k+1)}(a_{\zeta})(\zeta- \xi)\right) \to (2k+1)A(\xi).
\end{equation}
Since Lemma \ref{xi_in_middle_n} stipulates that $\xi$ lies between $\eta$ and $\zeta$, we conclude
\begin{align}
\label{lim_eta_xi_over_zeta_xi}
    \lim_{\eta \to \xi} \frac{\eta - \xi}{\zeta - \xi} =-1.
\end{align}
Similarly, a Taylor expansion of $f'$ gives
\begin{align*}
\label{in_limits:here1}
    f'(t) = e_{2k}(\xi) \cdot e_1(t)
    &=(t- \xi)^{2k-1}\left(\frac{2kA(\xi) + f^{(2k+1)}(b)(t- \xi)}{(2k)!}\right),
\end{align*}
for some $b$ between $t$ and $\xi$.
Thus,
\begin{align}
    \frac{e_{2k}(\xi)\cdot e_1(\eta)}{e_{2k}(\xi)\cdot e_{1}(R_{\xi}\eta)}=\frac{(\eta- \xi)^{2k-1}\left(2kA(\xi) + f^{(2k+1)}(b_{\eta})(\eta- \xi)\right)}{(\zeta- \xi)^{2k-1}\left(2kA(\xi) + f^{(2k+1)}(b_{\zeta})(\zeta- \xi)\right)}
\end{align}
for $b_{\eta}$ between $\eta$ and $\xi$ and $b_{\zeta}$ between $\zeta$ and $\xi$. Then, again as $\eta \to \xi$, $\zeta \to \xi$, $b_{\eta} \to \xi$, $b_{\zeta} \to \xi$ and
\begin{equation}
\label{in_limits:here2}
    \left(2kA(\xi) + f^{(2k+1)}(b_{\eta})(\eta- \xi)\right) \to 2kA(\xi) \text{  and  } \left(2kA(\xi) + f^{(2k+1)}(b_{\zeta})(\zeta- \xi)\right) \to 2k A(\xi).
\end{equation}
Hence,
\begin{equation}
    \lim_{\eta \to \xi}\frac{e_{2k}(\xi)\cdot e_1(\eta)}{e_{2k}(\xi)\cdot e_{1}(R_{\xi}\eta)}=\lim_{\eta \to \xi}\frac{(\eta- \xi)^{2k-1}}{(\zeta- \xi)^{2k-1}}
    = -1,
\end{equation}
and
\begin{equation}
    \lim_{\eta \to \xi}\Delta(\xi,\eta)
    = 1.
\end{equation}
Consequently \eqref{limJtilde} follows immediately by the algebra of limits:
\begin{equation}
    \lim_{\eta \to \xi} \tilde{J}(\xi, \eta) = |e_{1}(\xi)\cdot \left(e_{1}(\xi) + e_{1}(\xi)\right)| =2.
\end{equation}
To prove \eqref{limJ}, we define $g(t) := e_{2k}'(\xi) \cdot \phi(t)$ and use Corollary \ref{corol_derivatives_frenet_vectors} to give its Taylor expansion
\begin{align*}
    g(t)-g(\xi)
    &=(t- \xi)^{2k-1}\left(\frac{-2kA(\xi) + g^{(2k)}(c)(t- \xi)}{(2k)!}\right),
\end{align*}
where $c$ is between $\xi$ and $t$.
Hence, for some $c_{\eta}$ between $\eta$ and $\xi$ and some $c_{\zeta}$ between $\zeta$ and $\xi$, we have
\begin{equation}
    (e_{2k}'(\xi))\cdot (\phi(\eta) - \phi(\zeta)) = (\eta- \xi)^{2k-1}\left(\frac{-2kA(\xi) + g^{(2k)}(c_{\eta})(\eta- \xi)}{(2k)!}\right) - (\zeta- \xi)^{2k-1}\left(\frac{-2kA(\xi) + g^{(2k)}(c_{\zeta})(\zeta- \xi)}{(2k)!}\right).
\end{equation}
Hence, by \eqref{in_limits:here1},
\begin{equation}
    J(\xi, \eta) = \left(\frac{\eta - \xi}{\zeta - \xi}\right)^{2k-1}\left( \frac{-2kA(\xi) + g^{(2k)}(c_{\eta})(\eta- \xi)}{2kA(\xi) + f^{(2k+1)}(b_{\zeta})(\zeta-\xi)} 
    \right)
    -
     \frac{-2kA(\xi) + g^{(2k)}(c_{\zeta})(\zeta- \xi)}{2kA(\xi) + f^{(2k+1)}(b_{\zeta})(\zeta-\xi)}.
\end{equation}

Again, as $\eta \to \xi$, we have $c_{\eta} \to \xi$, and $c_{\zeta} \to \xi$ and
\begin{equation}
    \frac{-2kA(\xi) + g^{(2k)}(c_{\eta})(\eta- \xi)}{(2k)!}\to \frac{-A(\xi)}{(2k-1)!}, \qquad \text{and} \qquad \frac{-2kA(\xi) + g^{(2k)}(c_{\zeta})(\zeta- \xi)}{(2k)!} \to \frac{-A(\xi)}{(2k-1)!}.
\end{equation}
Also using \eqref{in_limits:here2}, \eqref{lim_eta_xi_over_zeta_xi}, and the algebra of limits, we conclude that
\begin{align*}
\lim_{\eta \to \xi} J(\xi, \eta)
    = 2.
\end{align*}

\end{proof}

\section{Curve-Carried Bilinear Fractional Integral Bounds}
\label{Section_bilinear_bounds}
\begin{theorem}[Surface-carried bilinear fractional integral bound]
\label{bilinear_frac_higher_dim}
Let $S$ be a well-curved, compact patch of a $C^{n+2k -1}$ curve in $\mathbb{R}^n$ with arclength parametrisation, $\phi: I \to \mathbb{R}^n$ satisfying conditions \eqref{eqn_ej_dot_T_not_zero_condition}, \eqref{eqn_taylor_condn}, \eqref{eqn_ej'_dot_T_not_zero_condition} and \eqref{eqn_tangents_condition_n}.
    Additionally, take $0<s<1$, and $\alpha \in [\frac{1}{2}, 1]$. Define
    \begin{equation}
        \mathcal{B}_{S,s}(g_1, g_2) = \int_I \frac{g_1(\eta)g_2(R_{\xi}\eta)}{|\phi(\eta)-\phi(R_{\xi}\eta)|^s}J(\xi, \eta)\rmd\eta.
    \end{equation}
    Then,
    \begin{equation}
        \|\mathcal{B}_{S, s}(g_1, g_2)\|_{L^{\alpha}(I)}\lesssim \|g_1\|_{L^{2\alpha}(I)}\|g_2\|_{L^{2\alpha}(I)},
    \end{equation}
    with implicit constant depending on $\alpha, s$ and the diameter of $S$.
\end{theorem}
The below proof follows in a similar way to the surface-carried bilinear-fractional integral bound proven by Bennett, Gutiérrez, Nakamura and Oliveira in \cite{BGNO} inspired also by Kenig and Stein \cite{KenigStein}.
\begin{proof}[Proof of Theorem \ref{bilinear_frac_higher_dim}]
We begin by noting that for the dyadic scales $\lambda \lesssim \diam(S) = \sup_{s, t \in I}|\phi(s) - \phi(t)|$, we have,
\begin{align}
    \|\mathcal{B}_{S,s}(g_1, g_2)\|_{L^{\alpha}(I)}^{\alpha} &= \int_I \left(\int_I \frac{g_1(\eta)g_2(R_{\xi}\eta)}{|\phi(\eta)-\phi(R_{\xi}\eta)|^s}J(\xi, \eta)\rmd\eta \right)^{\alpha} \rmd\xi
    \\
    &\leq \sum_{0 < \lambda 
    \lesssim \diam(S)}\int_I \left(
    \int_{\eta \in I: |\phi(\eta)-\phi(R_{\xi}\eta)| \approx  \lambda} \frac{g_1(\eta)g_2(R_{\xi}\eta)}{|\phi(\eta) - \phi(R_{\xi}\eta)|^s}J(\xi, \eta)\rmd \eta  \right)^{\alpha} \rmd \xi\\
    &\approx 
    \sum_{0 < \lambda \lesssim \diam(S)}
    \int_I \left(
    \int_{\eta \in I: |\phi(\eta) - \phi(R_{\xi}\eta)| \approx  \lambda} \frac{g_1(\eta)g_2(R_{\xi}\eta)}{\lambda^s}J(\xi, \eta)\rmd \eta\right)^{\alpha} \rmd \xi \\
    &= \sum_{0 < \lambda \lesssim \diam(S)}
    \lambda^{-{\alpha}s}\int_I \left(\int_{\eta \in I: |\phi(\eta) - \phi(R_{\xi}\eta)| \approx  \lambda} g_1(\eta)g_2(R_{\xi}\eta)J(\xi, \eta)\rmd \eta \right)^{\alpha} \rmd \xi\\
    &:= \sum_{0 < \lambda \lesssim \diam(S)}
    \lambda^{-{\alpha}s} \|F_{\lambda}\|_{L^{\alpha}(I)}^{\alpha}.
\end{align}
Here we have defined $$F_{\lambda}(\xi) := \int_{\eta \in I: |\phi(\eta)-\phi(R_{\xi}\eta)| \approx  \lambda} g_1(\eta)g_2(R_{\xi}\eta)J(\xi, \eta)\rmd \eta.$$
\noindent
Let ${\alpha} \in [\frac{1}{2}, 1]$ and fix a dyadic scale, $\lambda \lesssim \diam(S)$. Decomposing $I$ into a collection, $\Theta_\lambda$, of caps $\theta$ such that $|\theta| \approx \lambda$ followed by an application of Hölder's inequality yields
\begin{align}
\|F_{\lambda}\|_{L^{\alpha}}^{\alpha} 
&= \int_I \left( \int_{|\phi(\xi) - \phi(R_{\xi}\eta)| \approx \lambda} g_1(\eta)g_2(R_{\xi}\eta)J(\xi, \eta)\rmd \eta \right)^{\alpha} \rmd \xi\\
&= \sum_{\theta \in \Theta_\lambda}\int_\theta 
\left( \int_{|\phi(\eta) - \phi(R_{\xi}\eta)| \approx \lambda} g_1(\eta)g_2(R_{\xi}\eta)J(\xi, \eta)\rmd \eta
\right)^{\alpha}
\rmd \xi \\
&\lesssim \lambda^{1-{\alpha}} \sum_{\theta \in \Theta_\lambda}\left(\int_\theta  \int_{|\phi(\eta) - \phi(R_{\xi}\eta)| \approx \lambda} g_1(\eta)g_2(R_{\xi}\eta)J(\xi, \eta)\rmd \eta
\rmd \xi
\right)^{\alpha}.
\end{align}
Now observe that for all $\xi, \eta \in I$, \eqref{phia-phib_sim_a-b} and the in-between Lemma \ref{xi_in_middle_n} imply
 \begin{align}
        | \phi(\xi) - \phi(\eta)| \leq |\eta - \xi| \leq |\eta - \Rmap| \leq \sqrt 2 |\phi(\eta) - \phi(\Rmap)|. \label{distance_1}
    \end{align}
Therefore, if $|\phi(\eta) - \phi(R_{\xi}\eta))| \approx \lambda$, then  $|\phi(\xi) - \phi(\eta)| \lesssim \lambda$. 
Since also $\xi \in \theta$ with $|\theta| \approx \lambda$, we conclude $\eta \in \theta^*$ for $\theta^*$ a dilation of $\theta$ such that $2|\theta| = |\theta^*|$.
Similarly,
\begin{equation}
\label{distance_2}
        | \phi(\xi) - \phi(\Rmap)| \leq |\phi(\eta)-\phi(\Rmap)|,
\end{equation}
and hence $R_{\xi}\eta \in \theta^*$.
Therefore,
\begin{align}
\label{in_proof:F_lambda_norm}
\|F_{\lambda}\|_{L^{\alpha}}^{\alpha}
    &\lesssim \lambda^{1-{\alpha}} \sum_{\theta \in \Theta_\lambda}\left(\int_I  \int_{I} g_1\1_{\theta^*}(\eta)g_2\1_{\theta^*}(R_{\xi}\eta)J(\xi, \eta)\rmd \eta
\rmd \xi
\right)^{\alpha}.
\end{align}
Now, 
 \begin{equation}
     \int_I \int_I f(\eta)g(R_{\xi}\eta)J(\xi,\eta)\rmd \xi\rmd \eta = \|f\|_{L^1(I)}\|g\|_{L^1(I)}.
 \end{equation}
 Applying this identity to the right-hand side of \eqref{in_proof:F_lambda_norm}, followed by applications of Cauchy--Schwarz and Hölder's inequalities in turn produces
\begin{align}
\|F_{\lambda}\|_{L^{\alpha}}^{\alpha}
    &\lesssim \lambda^{1-{\alpha}} \sum_{\theta \in \Theta_\lambda}
\|g_1\1_{\theta^*}\|_{L^1(I)}^{\alpha}\|g_2\1_{\theta^*}\|_{L^1(I)}^{\alpha} \\
&\leq \lambda^{1-{\alpha}} \left(\sum_{\theta \in \Theta_\lambda}
\|g_1\1_{\theta^*}\|_{L^1(I)}^{2\alpha}\right)^{\frac{1}{2}}\left(\sum_{\theta \in \Theta_\lambda}
\|g_2\1_{\theta^*}\|_{L^1(I)}^{2\alpha}\right)^{\frac{1}{2}}
\\
&\leq \lambda^{1-{\alpha}} (2\lambda)^{\frac{2{\alpha}}{{\beta}'}} \left(\sum_{\theta \in \Theta_\lambda}
\|g_1\1_{\theta^*}\|_{L^{\beta}(I)}^{2{\alpha}}\right)^{\frac{1}{2}}\left(\sum_{\theta \in \Theta_\lambda}
\|g_2\1_{\theta^*}\|_{L^{\beta}(I)}^{2{\alpha}}\right)^{\frac{1}{2}}
\end{align}
for ${\beta} \geq 1$. In particular, we take $\beta = 2 \alpha$ since $\alpha \geq \frac{1}{2}$, implying that $\frac{2{\alpha}}{{\beta}'} = 2 \alpha -1$.
Observe that the collection, $\Theta_{\lambda}^*$, of $\theta^*$ covers $I$ with an overlap factor of $2$.
This implies that
\begin{equation}
    \left(\sum_{\theta \in \Theta_\lambda}
\|f\1_{\theta^*}\|_{L^{2\alpha}(I)}^{2{\alpha}}\right)^{\frac{1}{2}} \leq (2\|f\|_{L^{2\alpha}(I)}^{2{\alpha}})^{1/2} = \sqrt{2}\|f\|^{\alpha}_{L^{2\alpha}(I)},
\end{equation}
and hence
\begin{equation}
\label{Thm7.2b_Lq_result}
    \|F_{\lambda}\|_{L^{\alpha}}^{\alpha} \lesssim \lambda^{\alpha}\|g_1\|_{L^{\beta}(I)}^{\alpha} \|g_2\|_{L^{\beta}(I)}^{\alpha}.
\end{equation}
Combining the above, we have that
\begin{equation}
    \|\mathcal{B}_{S,s}(g_1, g_2)\|^{\alpha}_{L^{\alpha}(I)} = \sum_{0 < \lambda \lesssim \diam(S)}
    \lambda^{-{\alpha}s} \|F_{\lambda}\|_{L^{\alpha}(I)}^{\alpha} 
    \lesssim \sum_{0 < \lambda \lesssim \diam(S)}
    \lambda^{(1-s)\alpha}\|g_1\|_{L^{\beta}(I)}^{\alpha} \|g_2\|_{L^{\beta}(I)}^{\alpha}
\end{equation}
which converges since $1/2 < {\alpha}<1$ and $s<1$.
\end{proof}

\begin{remark}
    We note that we could have instead proven a bilinear fractional integral estimate for
    \begin{equation}
        \int_I \frac{g_1(\eta)\overline{g_2(R_{\xi}\eta)}}{|e_1(\xi)\cdot (\phi(\eta)-\phi(R_{\xi}\eta))|^s}J(\xi, \eta)\rmd\eta
        \qquad
        \text{or}
        \qquad
        \int_I \frac{g_1(\eta)\overline{g_2(R_{\xi}\eta)}}{|\eta-R_{\xi}\eta|^s}J(\xi, \eta)\rmd\eta.
    \end{equation}
The respective bilinear fractional integral estimates would follow in the same way, but make use of
\begin{equation}
\label{all_sim}
    |e_1(\xi)\cdot (\phi(\eta)-\phi(R_{\xi}\eta)| \sim |\eta - \Rmap| \sim |\phi(\eta) - \phi(\Rmap)|
\end{equation}
from Lemma \ref{gammaT_sim_a-b_ndim}.
Notice that in the proof of the Sobolev--Stein inequality, Theorems \ref{corol2a_l2sobolev} and \ref{corol2a_l2sobolev_n}, the equivalence \eqref{all_sim} was used to reduce the problem to one concerning the bilinear fractional integral $\mathcal{B}_{S, s}$, but could have equivalently reduced the problem to one of the above bilinear fractional integrals.
Explicitly, \eqref{all_sim} implies that
    \begin{equation}
        \mathcal{B}_{S, s}(g_1, g_2) \sim \int_I \frac{g_1(\eta)\overline{g_2(R_{\xi}\eta)}}{|e_1(\xi)\cdot (\phi(\eta)-\phi(R_{\xi}\eta))|^s}J(\xi, \eta)\rmd\eta
        \sim
        \int_I \frac{g_1(\eta)\overline{g_2(R_{\xi}\eta)}}{|\eta-R_{\xi}\eta|^s}J(\xi, \eta)\rmd\eta,
    \end{equation}
    with uniform implicit constants.  The choice to study $\mathcal{B}_{S, s}$ was made to align with the work of \cite{BGNO}.
\end{remark}

\section*{Acknowledgements}
\textbf{Acknowledgments.} The author profusely thanks her PhD advisor, Jonathan Bennett, for his unwavering support, guidance, and impactful conversations without which this paper would not have been possible. She thanks him for introducing her to the problem with such joy and contagious enthusiasm and for reviewing many earlier drafts of the paper. Throughout the nearly three years since beginning this project, the author has had insightful and motivating conversations with many individuals whom she thanks, including Hannah Cairo, Tony Carbery, Gianluca Giacchi, Zane Li, and Po-Lam Yung. She especially thanks Itamar Oliveira and Shohei Nakamura for numerous conversations which have influenced this work.

\textbf{Funding.} The author is
supported by the EPSRC Doctoral Training Partnership and the University
of Birmingham.

\textbf{Use of Generative AI.} The impact of generative AI on the main results of this work has been minimal. Claude was used to check some of the elementary lemmas and remarks for typographical and logical issues: explicitly, it suggested small improvements to earlier versions of Remark \ref{corol_injectivity_of_delta}, Lemma \ref{lemma_derivatives_frenet_vectors}, Lemma \ref{lem:Ck_control}, and Lemma \ref{jacobian_limits_n}.
The results, mathematical ideas, and writing of this paper are all due to the author.

\bibliographystyle{abbrv}
\bibliography{ref}

\end{document}